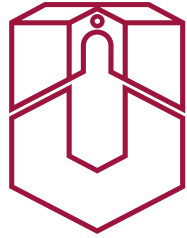

UNIVERSITÄT
OSNABRÜCK

# THE SYMMETRIC SIGNATURE

vorgelegt von

**Alessio Caminata**

Betreuer

Prof. Dr. **Holger Brenner**

To Elisa,

my love and inspiration.

# Contents

# Introduction

> In principle, you can define everything.
>
> ——————————————
>
> Gerd Faltings

A recurring procedure used in algebraic geometry and commutative algebra, and more in general in Mathematics, to understand a rather complicated object is to attach to it a number, a numerical invariant. This invariant should capture the properties we want to investigate, so isomorphic objects must have the same numerical invariant, and it should be relatively easy to compute, thus this is not a strict condition. Classical examples in algebraic geometry are the *dimension* of an algebraic variety or the *genus* of an algebraic curve.

Another important and classical example is the *Hilbert-Samuel multiplicity*. Given a local ring $(R, \mathfrak{m})$ of dimension $d$ (one can think of it as the ring of a point of an algebraic variety), we can associate to it the following number

$$e(R) := \lim_{n \to +\infty} \frac{d! \cdot l_R(R/\mathfrak{m}^{n+1})}{n^d},$$

which is always a positive integer and is called the *Hilbert-Samuel multiplicity* of $R$. This is an algebraic invariant which can be interpreted as a measure of the singularity of the point in the variety: the bigger the multiplicity is, the more complicated the singularity. Moreover one has that $e(R) = 1$ if and only if the ring is regular, that is the point is smooth.

The Hilbert-Samuel multiplicity is quite easy to compute in general, there exist algorithms that have been implemented on computer algebra softwares such as Macaulay2 [GS]. As a drawback, this multiplicity is not a fine invariant and is not able to distinguish different types of singularities. For example, the two-dimensional rings $k[\![x, y, z]\!]/(x^2 + y^3 + z^5)$ and $k[\![x, y, z]\!]/(x^2 + y^3 + z^7)$ have both multiplicity 2, but they have very different algebraic properties, the first one being a unique factorization domain, the latter not.

If we work over a field of prime characteristic $p$, we can use the Frobenius homomorphism given by $F(r) = r^p$ and build finer invariants to understand the singularities of the ring. A pioneer in this study was Kunz [Kun69]. In 1969 he introduced the function $q \mapsto l_R(R/\mathfrak{m}^{[q]})$ and the limit

$$e_{HK}(R) := \lim_{q \to +\infty} \frac{l_R(R/\mathfrak{m}^{[q]})}{q^d},$$

which are now known as *Hilbert-Kunz function* and *Hilbert-Kunz multiplicity* respectively. The difference with the classical setting is that $q = p^e$ varies only among powers of the characteristic and $\mathfrak{m}^{[q]}$ is the Frobenius power of the ideal, generated by all the $q$-th powers of elements of $\mathfrak{m}$.

In 1983, Monsky [Mon83] proved that the limit defining $e_{HK}(R)$ always exists and extended the definition to $R$-modules in general. In 2000, Watanabe and Yoshida [WY00] proved that $e_{HK}(R) = 1$ if and only if the ring is regular. So, the Hilbert-Kunz multiplicity detects regularity as the classical Hilbert-Samuel multiplicity does. Moreover, it is deeply connected with the theory of *tight closure* and *F-singularities* in prime characteristic (see for example the book of Huneke [Hun96] and the survey article of Schwede and Tucker [ST12]). On the other hand, the Hilbert-Kunz multiplicity is in general much more difficult to compute than the Hilbert-Samuel multiplicity, and is not necessarily an integer, and not even a rational number, as shown by Brenner [Bre13].

The main motivation for this thesis comes from a more recent invariant of prime characteristic rings: the *F-signature*. This was introduced in 2002 by Huneke and Leuschke [HL02], who continued the ideas of Smith and Van den Bergh [SVB96] of focusing on the splitting properties of the Frobenius homomorphism.

Consider $(R, \mathfrak{m}, k)$ a $d$-dimensional reduced Noetherian local ring of prime characteristic $p$ which is F-finite, and with algebraically closed residue field $k$. For every $e \in \mathbb{N}$, let $q = p^e$ and let $^eR$ be the $R$-module which is equal to $R$ as abelian group and has left multiplication twisted via Frobenius, that is $r \circ s := r^{p^e} s$. We decompose it as

$$^eR = R^{a_q} \oplus M_q,$$

where the module $M_q$ has no free direct summands. The limit

$$s(R) := \lim_{e \to +\infty} \frac{a_q}{q^d}$$

is called *F-signature* of $R$. The number $a_q$ is also called the *free rank* of $^eR$ and denoted by $\mathrm{frk}_R(^eR)$, while if $R$ is a domain, one has that $q^d = \mathrm{rank}_R(^eR)$, the usual rank of an $R$-module.

Huneke and Leuschke [HL02] proved that the limit defining the F-signature exists assuming that $R$ is Gorenstein. Then, other authors showed the existence in some other cases (see e.g. Yao [Yao05] for rings of finite F-representation type, or Singh [Sin05] for affine semigroup rings). Finally, in 2012 Tucker [Tuc12] proved that the F-signature exists for every reduced F-finite Noetherian local ring.

It is worth mentioning that, almost at the same time of Huneke and Leuschke, Watanabe and Yoshida [WY04] introduced a similar notion, called *minimal relative Hilbert-Kunz multiplicity*. Then, Yao [Yao06] proved it actually coincides with the F-signature.

The F-signature is a real number between 0 and 1 and provides delicate information about the singularities of $R$. Two principal results in this direction are the fact that $s(R) = 1$ if and only if the ring is regular, and that $s(R) > 0$ if and only if $R$ is strongly F-regular. The first one is a consequence of the analogous result of Watanabe and Yoshida [WY00] for the Hilbert-Kunz multiplicity, and the second one was proved by Aberbach and Leuschke [AL03]. There are other results which justify the statement that *F-signature measures singularities*. In general, the closer to 1 the F-signature is, the nicer the singularity.

The F-signature is an intrinsic concept of rings of prime characteristic, but it would be profitable to have an analogous invariant which is independent of the characteristic and could be used to study singularities also in characteristic zero. This problem is actually the main motivation of this thesis, and led us to the definition of the *symmetric signature*.

A first natural attempt to define a characteristic zero version of the F-signature is perhaps the following. If $R$ is a reduced $\mathbb{Z}$-algebra such that $\operatorname{Spec} R \to \operatorname{Spec} \mathbb{Z}$ is dominant, then for every prime number $p$ we consider its reduction mod $p$, $R_p := R \otimes_{\mathbb{Z}} (\mathbb{Z}/p\mathbb{Z})$, and compute the F-signature $s(R_p)$. One may ask whether the limit

$$\lim_{p \to +\infty} s(R_p)$$

exists, and use this limit to define a characteristic 0 version of the F-signature.

Unfortunately, it is difficult to find an appropriate meaning and compute the previous limit, or even determine whether it exists at all. In fact, if $p_1$ and $p_2$ are two distinct prime numbers, the two Frobenius homomorphisms $F_{p_1} : R_{p_1} \to R_{p_1}$ and $F_{p_2} : R_{p_2} \to R_{p_2}$ are quite different in general, and difficult to compare. For example, if $R$ is a domain then the rank of ${}^e R_{p_1}$ and the rank of ${}^e R_{p_2}$ are different, the first one being $p_1^e$ and the second one $p_2^e$.

For our definition of symmetric signature, we choose a different approach. We replace the module ${}^e R$ with another module $\mathscr{S}^q$ (here $q$ is simply a natural number), which is defined without using the Frobenius homomorphism or other properties of prime characteristic rings. This approach has the clear advantage that in a relative situation, when $R$ is a reduced $\mathbb{Z}$-algebra, we can give a meaning to the limit $\lim_{p \to +\infty} s(R_p)$. The existence of this limit is not clear also in this case, but at least now the modules $\mathscr{S}_p^q$ over the fibers $R_p$ are specializations of the same module $\mathscr{S}^q$ over $R$.

Consider a Noetherian local d-dimensional domain $(R, \mathfrak{m}, k)$, with algebraically closed residue field $k$ and consider $\operatorname{Syz}_R^d(k)$, the top-dimensional syzygy module of $k$. For every natural number $q$ we define $\mathscr{S}^q := \left(\operatorname{Sym}_R^q(\operatorname{Syz}_R^d(k))\right)^{**}$, where $(-)^*$ denotes the functor $\operatorname{Hom}_R(-, R)$. We decompose the last module as

$$\mathscr{S}^q = R^{a_q} \oplus M_q,$$

with the module $M_q$ containing no free direct summands, so that $a_q = \operatorname{frk}_R \mathscr{S}^q$. For ease of notation we fix also $b_q := \operatorname{rank}_R \mathscr{S}^q$, and we introduce the following (Definition 2.2.1).

**Definition.** The real number

$$s_\sigma(R) := \lim_{N \to +\infty} \frac{\sum_{q=0}^{N} a_q}{\sum_{q=0}^{N} b_q}$$

is called *symmetric signature* of $R$, provided the limit exists.

In our construction we looked at two important properties of the module $^eR$ which we would like to keep. First, if the ring $R$ is Cohen-Macaulay then $^eR$ is a maximal Cohen-Macaulay module. The second fact is an outstanding result of Kunz, which states that $^eR$ is a free module (for all or for some $e \in \mathbb{N}$) if and only if $R$ is a regular ring.

These properties are shared also by the top dimensional syzygy module of the residue field $\mathrm{Syz}_R^d(k)$, which is always maximal Cohen-Macaulay, and is $R$-free if and only if the ring $R$ is regular. These are consequences of the Depth Lemma (Lemma 1.4.9), of a famous result of Auslander, Buchsbaum, and Serre (Theorem 1.4.6), and of the Auslander-Buchsbaum formula (Theorem 1.4.3). In order to obtain an asymptotic behaviour, we apply the symmetric powers to $\mathrm{Syz}_R^d(k)$ and we take reflexive hull, obtaining the module $\mathscr{S}^q$. After Auslander, who introduced the reflexive tensor product, we call the functor $\mathrm{Sym}_R^q(-)^{**}$ *reflexive symmetric powers*.

Another interesting possibility would be to replace the module $\mathrm{Syz}_R^d(k)$ with the module $\Omega_{R/k}^{**}$ of *Zariski differentials of $R$ over $k$*, and then applying again reflexive symmetric powers. We consider also this case and we call the corresponding limit *differential symmetric signature* (Definition 2.2.15), but in most parts of the thesis we will rather concentrate on the symmetric signature defined with $\mathrm{Syz}_R^d(k)$.

As a fundamental test to sustain our definition of symmetric signature we focus on two important classes of examples: two-dimensional quotient singularities and coordinate rings of plane elliptic curves.

A two-dimensional quotient singularity $R$ over a field $k$ is the invariant ring under the linear action of a finite small group $G \subseteq \mathrm{GL}(2,k)$ on a power series ring $S = k[\![u,v]\!]$. Watanabe and Yoshida [WY04] proved that the F-signature of a quotient singularity $R$ over a field of prime characteristic $k$ is equal to $\frac{1}{|G|}$, provided that the characteristic of $k$ does not divide the order of the group.

We focus on two particular classes of quotient singularities $R = k[\![u,v]\!]^G$, with the characteristic of $k$ not dividing the order of $G$.

1. *Kleinian quotient singularities* or *ADE singularities*, that is when $G \subseteq \mathrm{SL}(2,k)$.
2. *Cyclic quotient singularities*, that is when $G$ is a cyclic subgroup of $\mathrm{GL}(2,k)$.

We compute the symmetric signature for these two classes of rings.

The methods we use to handle these cases live at the intersection of homological algebra and representation theory. In particular, our main tool is the so-called *Auslander correspondence* (Theorem 3.5.10). This states that under certain hypothesis there is a one to one correspondence between irreducible $k$-representations of $G$ and indecomposable maximal Cohen-Macaulay $R$-modules. Moreover, this correspondence is functorial, that is there exists a functor $\mathscr{A}$, called *Auslander functor*, which sends every $k$-representation $V$ to the maximal Cohen-Macaulay module $\mathscr{A}(V) := (S \otimes_k V)^G$.

The Auslander correspondence will be the main object of study of Section 3.4, where we will recall its proof and generalize it to the case where the ring $S$ is a normal $n$-dimensional domain, but not necessarily a power series ring.

Then, we will prove the following two fundamental results (Theorem 4.1.8 and Corollary 4.3.6), which will allow us to translate the problem of computing the symmetric signature of $R$ into a problem of representation theory of finite groups.

**Theorem.** *The Auslander functor $\mathscr{A}$ commutes with reflexive symmetric powers, that is for every $k$-representation $V$ of $G$ we have*

$$\mathscr{A}\left(\mathrm{Sym}_k^q(V)\right) \cong \left(\mathrm{Sym}_R^q(\mathscr{A}(V))\right)^{**}.$$

Here, $\mathrm{Sym}_R^q(-)$ denotes the $q$-th symmetric power of an $R$-module, and $\mathrm{Sym}_k^q(V)$ is the $q$-th symmetric power of the representation $V$, which is canonical isomorphic to its double dual.

**Theorem.** *Let $R = S^G$ be a two-dimensional Kleinian singularity, and let $V_1$ be the two-dimensional fundamental representation which defines the action of $G$ on $S$. Then, the second syzygy module of the residue field $k$ is isomorphic to the image of $V_1$ via Auslander functor, that is*

$$\mathrm{Syz}_R^2(k) \cong \mathscr{A}(V_1).$$

The second theorem extends a result of Yoshino and Kawamoto [YK88], who proved that in this situation $\mathscr{A}(V_1)$ is isomorphic to the third syzygy of the residue field, $\mathrm{Syz}_R^3(k)$.

As a consequence of the previous two theorems it turns out that for a two-dimensional Kleinian singularity the rank of $\mathscr{S}^q$ is just the dimension of the representation $\mathrm{Sym}_k^q(V_1)$, that is $q+1$, and the free rank of $\mathscr{S}^q$ is equal to the multiplicity of the trivial representation into $\mathrm{Sym}_k^q(V_1)$. This multiplicity can be more easily computed using tools from representation theory of finite groups, such as character theory. Thus, we obtain the following result.

**Theorem.** *Let $R = S^G$ be a Kleinian singularity over an algebraically closed field $k$ such that* char$\,k$ *does not divide the order of $G$. Then the symmetric signature of $R$ is*

$$s_\sigma(R) = \frac{1}{|G|}.$$

The same result is obtained also for the differential symmetric signature (Theorem 4.4.11).

For cyclic quotient singularities the situation is a little bit more complicated, since $\mathrm{Syz}_R^2(k)$ is not necessarily isomorphic to $\mathscr{A}(V_1)$, and it does not even have rank 2 in general. However it is a maximal Cohen-Macaulay module, so it is isomorphic to the image via Auslander functor $\mathscr{A}$ of a $k$-representation of $G$. We do not determine explicitly which representation it is, but we prove it is a faithful representation (Lemma 4.5.5). Then, using again techniques from representation theory we are able to prove the following result, which is obtained in collaboration with Lukas Katthän.

**Theorem.** *Let $R = S^G$ be a cyclic singularity over an algebraically closed field $k$ such that* char$k$ *does not divide the order of $G$. Then the symmetric signature of $R$ is*

$$s_\sigma(R) = \frac{1}{|G|}.$$

The second class of examples we are interested in are coordinate rings of plane projective elliptic curves over an algebraically closed field $k$, that is $R = k[x, y, z]/(f)$ with $f$ a homogeneous non-singular polynomial of degree 3. If $k$ has prime characteristic, then $R$ is not strongly F-regular, so by the result of Aberbach and Leuschke [AL03] the F-signature of $R$ is 0. We would like to prove that also the symmetric signature of $R$ is 0.

The methods we use to handle this case are geometric. We use the correspondence between graded MCM $R$-modules and vector bundles over the smooth projective curve $Y = \mathrm{Proj}R$, and we translate the problem of computing $s_\sigma(R)$ into an analogous problem in the category VB$(Y)$ of vector bundles over $Y$. The main advantage of this approach is that over an elliptic curve $Y$ the structure of the category VB$(Y)$ has been clarified by Atiyah [Ati57]. He gave a description of the indecomposable vector bundles on $Y$, and their behaviour under the tensor product operation.

Using Atiyah's classification we are able to prove the desired result for the differential symmetric signature defined with $\Omega^{**}_{R/k}$ (Theorem 5.2.1).

**Theorem.** *Let $Y$ be a plane elliptic curve over an algebraically closed field $k$ of characteristic $\geq 5$ with coordinate ring $R$. Then the differential symmetric signature of $R$ is $s_{d\sigma}(R) = 0$.*

Unfortunately we are not able at the moment to obtain the same result for the symmetric signature. However we present two strategies that allow us to get partial results and a better understanding of the decomposition of the module $\mathscr{S}^q$ involved in the definition of symmetric signature. In particular, we prove that over the field of the complex numbers $s_\sigma(R) \leq \frac{1}{2}$, provided the limit exists (Corollary 5.2.6).

## Summary of the thesis

In the first chapter we present and discuss some classical notions in commutative algebra. The contents of this chapter are not new, on the contrary they are important and established results that are necessary to define and understand the symmetric signature. In particular, in Section 1.1 we present the definition and examples of Krull-Schmidt categories. Roughly speaking, these are the categories, where each object can be uniquely decomposed as a direct sum of indecomposable objects. In Section 1.2 we give the definition of a syzygy of a module over a Noetherian commutative ring, and in the particular case of a local Noetherian ring. Then, in Section 1.3 we give a small survey on reflexive modules. In Section 1.4 we present some very important results in commutative algebra, which will appear often throughout the thesis. These are for example the Auslander-Buchsbaum formula (Theorem 1.4.3) and the Depth Lemma (Lemma 1.4.9). Finally, in Section 1.5 we introduce the class of maximal Cohen-Macaulay modules and present some of their most important properties.

Chapter 2 is dedicated to the F-signature and the symmetric signature. In Section 2.1 we present a short survey on F-signature and its properties. We give no proofs, but rather we focus on the examples and ideas that we consider more significant to motivate our choice to define the symmetric signature. This will be defined in Section 2.2, together with its variant, the differential symmetric signature. In this section we will present also some easy results and observations that follow almost directly from our definition.

The third chapter is dedicated to quotient singularities and to the Auslander correspondence. We give a short review of basic notions of representation theory of finite groups (Section 3.1), Brauer characters (Section 3.2), and skew-group ring (Section 3.3). Then, we will focus on the Auslander correspondence and we will present a proof (Section 3.4). Our proof follows mainly the original ideas of Auslander [Aus86b] (see also [Yos90], [LW12], and [IT13]), but it is more general, since we do not assume that the ring $S$ where the finite group acts on is regular and two-dimensional. Moreover, we present an alternative geometric proof of certain well-established facts (e.g. Lemma 3.4.11). Finally, in Section 3.5 we recall the main properties of quotient singularities, and we provide some classical examples of the Auslander correspondence applied to this setting.

Chapter 4 contains some of the main results of this thesis. We prove here that the Auslander functor commutes with reflexive symmetric powers (Section 4.1), and that the second syzygy of the residue field of a Kleinian singularity is isomorphic to the fundamental module of Auslander theory (Section 4.3). In Section 4.4 we compute the symmetric signature of the Kleinian singularities, and in Section 4.5 we compute the symmetric signature of a two-dimensional cyclic singularity. The differential variant is also considered and computed here.

In the last chapter we present some partial results concerning the symmetric signature of the coordinate ring of a plane elliptic curve over an algebraically closed field. In Section 5.1 we give a short review of the methods we use to attack this problem. We recall the definition and the basic properties of vector bundles over a smooth projective curve, and then we focus on Atiyah's classification of vector bundles over an elliptic curve (Section 5.1.2). In Section 5.2 we prove that the differential symmetric signature of the coordinate ring of a plane elliptic curve is 0 (Theorem 5.2.1), and we prove that over $\mathbb{C}$ the symmetric signature, if it exists, is at most $\frac{1}{2}$ (Corollary 5.2.6). Finally, in Section 5.2.1 and Section 5.2.2 we present two possible strategies that hopefully may be completed to prove that the symmetric signature is 0 in this case.

## Acknowledgements

First of all, I want to thank my advisor Holger Brenner. He has guided me with patience and experience to reach the results contained in this thesis. Without his help and his friendly approach this would have not been possible.

I would like to thank Winfried Bruns for his kindness and humanity. He has been always helpful and very nice to me and my girlfriend Elisa. He has done a lot to make my stay in Osnabrück much more confortable.

During these four years in Osnabrück, I have been in the middle of two generations of Ph.D. students and post docs. So I have been lucky enough to get to know, and benefit from the presence of both the *old* and the *new* generation. Let me thank both for the stimulating conversations and for all the fun we had. Their presence has made the math department in Osnabrück a warm and familiar working place. Among them, I thank Lukas Katthän, whom I had the pleasure to work with, and I thank Davide Alberelli for his great help with LATEX and for all the stylistic advices.

During my Ph.D. studies I attended many conferences, workshops, and schools around the world, where I met a lot of other young mathematicians. I shared many ideas with them, and part of my growing as a mathematician and as a man is due to these fruitful and frequent conversations. I am grateful to them for this reason, and I am looking forward to have this possibility again. In particular, I want to thank Yusuke Nakajima for all discussions about Auslander theory that we had in Tateyama, Barcelona, and during his visit in Osnabrück. I learned a lot from him, and I hope to meet him again in the future, and continue our conversation.

I want to thank Alberto Fernandez Boix for suggesting me the notation $s_\sigma(-)$ for the symmetric signature.

I want to thank my parents, Francesco and Nicoletta, and my sister Ilaria for all the support that they gave me during these years. It was not easy for me to be distant from my family, but they did everything they could to fill this gap and make me feel always serene.

The greatest thanks of all goes to my girlfriend Elisa. She has been by my side during all my Ph.D. studies, always supporting me and sustaining my choices, even if this was sometimes not easy for her. Words cannot express how grateful I am to her, but I am sure she knows.

Finally, even if he will probably never read these lines, I want to thank David Robert Jones. His amazing work for more than fourty years has inspired me a lot during the writing of this thesis.

# Notation

In this thesis we will use notations according mainly to the following textbooks.

1. For what concerns commutative algebra, we refer to the books of Atiyah and Macdonald [AMD69], and Bruns and Herzog [BH98].
2. For what concerns algebraic geometry, we refer to the books of Hartshorne [Har77], and Le Potier [LeP97].
3. For what concerns representation theory, we refer to the books of Fulton and Harris [FH91], and Feit [Fei82], and the notes of Webb [Web14].

In particular, we will adopt the following conventions.

- When we say that $R$ is a *ring*, we always mean a non-zero, associative, unitary ring. And by homomorphism of rings $f : R \to S$, we always assume that $f(1_R) = 1_S$.

- In almost all thesis, with the exception of some parts of Chapter 1 and Chapter 3, our rings will be also Noetherian and commutative.

- We denote by $\mathbb{N}$ the set of *natural numbers* including 0, that is $\mathbb{N} = \{0, 1, 2, 3, 4, \ldots\}$.

- We denote by $\mathbb{Z}$ the set of *integer numbers*, by $\mathbb{Q}$ the set of *rational numbers*, and by $\mathbb{C}$ the set of *complex numbers*.

- An ideal $I$ of a ring $R$ will be always different from $R$.

- We denote by $(R, \mathfrak{m}, k)$ or by $(R, \mathfrak{m})$ a *local ring*, that is a commutative unitary ring with a unique maximal ideal $\mathfrak{m}$, where $k = R/\mathfrak{m}$.

- Given a commutative ring $R$, an $R$-module $M$ and a multiplicative system $S \subseteq R$, we will denote by $S^{-1}M$ the *localization of $M$ at $S$*. If $\mathfrak{p}$ is a prime ideal of $R$, we will write $M_\mathfrak{p}$ instead of $(R \setminus \mathfrak{p})^{-1}M$.

- We say that a commutative ring $R$ is *graded* if it is $\mathbb{N}$-graded, that is if there is a decomposition $R = \bigoplus_{i \in \mathbb{N}} R_i$ such that each $R_i$ is an $R_0$-module and $R_i R_j \subseteq R_{i+j}$. The irrelevant ideal of $R$ is $R_+ := \bigoplus_{i > 0} R_i$. An element $x$ of $R$ has degree $i$ if $x \in R_i$. Therefore 0 has degree $i$ for all $i \in \mathbb{N}$.

- We say that an algebra $R$ is *standard graded* over a field $k$ if $R$ is graded, $R_0 = k$, and $R$ is generated by finitely many elements of degree 1, i.e. $R = k[R_1]$.

- The *spectrum* of a commutative ring $R$ is the topological space

$$\operatorname{Spec} R := \{\mathfrak{p} \subseteq R : \mathfrak{p} \text{ is a prime ideal}\}$$

intended with the Zariski topology.

- The *projective spectrum* of a commutative graded ring $R$ is the topological space

$$\operatorname{Proj} R := \{\mathfrak{p} \subseteq R : \mathfrak{p} \text{ is a homogeneous prime ideal not containing } R_+\}$$

intended with the Zariski topology.

Further references and notations will be explained throughout the thesis.

# 1. Preliminaries and definitions

> I don't know where I'm going from
> here, but I promise it won't be boring.
>
> ——————————————————————
>
> David Bowie

The main goal of this chapter is to present and discuss the notions and the results which are necessary to define and understand the symmetric signature. All theorems in this chapter are classical and well-known, so we skip proofs and we refer to appropriate sources. We reserve to give a proof in two cases: if we did not find a suitable reference or if we think that the proof is particularly enlightening to understand the key ideas.

Our main references for this chapter are the following: the book of Bruns and Herzog [BH98] for Section 1.2 on syzygies and Section 1.4 on depth and Serre's conditions, the books of Leuschke and Wiegand [LW12], and Yoshino [Yos90] for Section 1.1 on Krull-Schmidt categories and ranks, and Section 1.5 on maximal Cohen-Macaulay modules, and the book of Lam [Lam98] for Section 1.3 on reflexive modules. Other important references used will be outlined throughout the presentation.

## 1.1. Krull-Remak-Schmidt property and ranks

The name *Krull-Remak-Schmidt property* (KRS property, for short) refers to the so-called Krull-Remak-Schmidt Theorem. Roughly speaking, this states that a group, satisfying certain finiteness conditions, can be uniquely written as a finite direct product of indecomposable subgroups. It was first proved by Wedderburn [Wed09], and then generalized and extended by Remak [Rem11], Schmidt [Sch13], and Krull [Kru25]. The name Krull-Remak-Schmidt Theorem is nowadays used for results that concern a unique decomposition also in settings different from group theory, such as modules over Artinian and Noetherian rings and coherent sheaves.

It appears to us, that the KRS property is better explained using the language of category theory. We will use the terminology of Mac Lane's book [Mac71] and we assume that the reader is familiar with the concept of an abelian category. This is an additive category in which every map has a kernel and a cokernel, and in which every monomorphism (respectively epimorphism) is a kernel (respectively cokernel).

Fundamental examples of abelian categories are categories of modules over rings. We introduce now the notation for two of those categories. Even if we will work mainly with commutative rings, it is useful to define these categories also in the non-commutative setting.

**Definition 1.1.1.** Let $R$ be a ring, we denote by $\text{mod}(R)$ the category whose objects are finitely generated left modules and whose morphisms are $R$-linear maps between them. The category $\text{proj}(R)$ is the full subcategory of $\text{mod}(R)$ consisting of projective finitely generated left $R$-modules.

Given two objects $X$ and $Y$ in an abelian category $\mathscr{C}$, there exists a unique object $X \oplus Y$ in $\mathscr{C}$, called the direct sum of $X$ and $Y$. This motivates the following definitions.

**Definition 1.1.2.** Let $X$ be a non-zero object of $\mathscr{C}$. $X$ is called *indecomposable* if $X = X_1 \oplus X_2$ implies $X_1 = 0$ or $X_2 = 0$. $X$ is called *simple* if $X_1 \subseteq X$ implies $X_1 = 0$ or $X_1 = X$.

Clearly, simple objects are indecomposable, but the converse is in general not true.

**Example 1.1.3.** We consider the category $\text{mod}(\mathbb{Z})$ of $\mathbb{Z}$-modules. The simple objects of $\text{mod}(\mathbb{Z})$ are the cyclic modules $\mathbb{Z}/p\mathbb{Z}$ of prime order. The module $\mathbb{Z}$ is indecomposable, but is not simple.

**Definition 1.1.4.** Let $X$ be an object of an abelian category $\mathscr{C}$. We denote by $\text{Add}(X)$ the full subcategory consisting of all finite direct sums of copies of $X$ and their direct summands.

The category $\text{Add}(X)$ is the smallest additive subcategory of $\mathscr{C}$ which contains $X$ and is closed under taking direct summands.

**Definition 1.1.5.** An additive category $\mathscr{C}$ is said to be a *Krull-Schmidt category* or to have the *KRS property* if every object of $\mathscr{C}$ decomposes into a finite direct sum of objects having local endomorphism rings.

Notice that the endomorphism ring $\text{End}(X)$ of an object $X$ in $\mathscr{C}$ may be non-commutative. One says that $\text{End}(X)$ is local if the sum of non-units in $\text{End}(X)$ is a non-unit.

In Krull-Schmidt categories we have an analogue of the Krull-Remak-Schmidt Theorem.

**Theorem 1.1.6** (Krull-Remak-Schmidt)**.** *Let $\mathscr{C}$ be a Krull-Schmidt category, then the following facts hold.*

*1) An object $X$ is indecomposable if and only if $\text{End}(X)$ is local.*
*2) Every object is isomorphic to a finite direct sum of indecomposable objects.*
*3) If $X$ is an object of $\mathscr{C}$ and there are two decompositions*

$$X_1 \oplus \cdots \oplus X_r = X = Y_1 \oplus \cdots \oplus Y_s$$

*into indecomposable objects, then $r = s$ and there exists a permutation $\pi$ such that $X_i \cong Y_{\pi(i)}$ for $i = 1, \ldots, r$.*

**Remark 1.1.7.** Let $R$ be a commutative local ring such that $\mathrm{mod}(R)$ has the KRS property. Then the endomorphism ring $\mathrm{End}_R(R) \cong R$ is local, so $R$ is an indecomposable object in $\mathrm{mod}(R)$.

Our main interest is for categories of modules over some ring $R$, but the KRS property may fail in these categories, even if the ring is regular.

**Example 1.1.8.** The category $\mathrm{mod}(\mathbb{Z})$ does not have the KRS property. Although every finitely generated $\mathbb{Z}$-module is a direct sum of indecomposable objects, the indecomposable module $\mathbb{Z}$ does not have a local endomorphism ring.

**Example 1.1.9** (Leuschke-Wiegand [LW12])**.** Let $R = k[x,y]$, $\mathfrak{m} = (x,y)$ and $\mathfrak{n} = (x-1,y)$. We have $\mathfrak{m} + \mathfrak{n} = R$, so we get a short exact sequence

$$0 \to \mathfrak{m} \cap \mathfrak{n} \to \mathfrak{m} \oplus \mathfrak{n} \to R \to 0.$$

Since $R$ is projective, the sequence splits, so $\mathfrak{m} \oplus \mathfrak{n} \cong R \oplus (\mathfrak{m} \cap \mathfrak{n})$. Neither $\mathfrak{m}$ nor $\mathfrak{n}$ is isomorphic to $R$, therefore the category $\mathrm{mod}(R)$ is not KRS.

We present now some examples of Krull-Schmidt categories without proofs. If not otherwise stated, the reader should consult the book of Leuschke and Wiegand [LW12, Chapter 1] or the notes of Krause [Kra14] for proofs and further details.

**Theorem 1.1.10.** *Let $R$ be a commutative ring, then the following facts hold.*

*1) If $R$ is Artinian, then $\mathrm{mod}(R)$ has the KRS property.*
*2) If $R$ is local and complete, then $\mathrm{mod}(R)$ has the KRS property.*

**Remark 1.1.11.** The condition complete in part *2)* of Theorem 1.1.10 could be replaced with Henselian and the statement remains true. We recall that by Hensel's lemma complete local rings are Henselian.

**Theorem 1.1.12** (Atiyah [Ati56])**.** *Let $(X, \mathscr{O}_X)$ be a complete algebraic variety over an algebraically closed field. Then the KRS property holds for the category $\mathrm{Coh}(X)$ of coherent $\mathscr{O}_X$-modules. Moreover if $X$ is connected then KRS holds also for the category $\mathrm{VB}(X)$ of vector bundles over $X$.*

### 1.1.1. Rank and free rank

Let $R$ be a commutative ring, and let $S$ be the multiplicative system consisting of the non-zero-divisors of $R$. If $R$ is reduced, then $S$ is the complement of the union of the minimal primes of $R$. The *total ring of fractions* of $R$ is the ring $Q(R) := S^{-1}R$.

**Definition 1.1.13.** Let $M$ be an $R$-module, we say that *$M$ has rank*, if the module $S^{-1}M$ is a free $Q(R)$-module. In this case we define

$$\mathrm{rank}_R M := \mathrm{rank}_{Q(R)} S^{-1}M.$$

**Remark 1.1.14.** If $R$ is a domain, then $S = R \setminus \{0\}$ and $Q(R)$ is a field. So every finitely generated $R$-module has rank.

**Definition 1.1.15.** Let $M$ be a finitely generated $R$-module. We define the *free rank* of $M$ as

$$\mathrm{frk}_R(M) := \max\{n : \exists \text{ a split surjection } \varphi : M \twoheadrightarrow R^n\}.$$

If the module $M$ has finite rank, we can write it as $M \cong R^{\mathrm{frk}_R(M)} \oplus N$, where the module $N$ has no free direct $R$-summands. It follows that in general

$$\mathrm{frk}_R M \leq \mathrm{rank}_R M.$$

**Example 1.1.16.** Let $M$ be a finite free $R$-module, then $M$ has rank and $\mathrm{rank}_R M = \mathrm{frk}_R M$.

**Example 1.1.17.** Let $R$ be a principal ideal domain, then every finitely generated $R$-module $M$ can be decomposed as $M = R^n \oplus T(M)$, where $T(M)$ is the torsion sub-module of $M$. It follows that $n = \mathrm{frk}_R M = \mathrm{rank}_R M$.

**Remark 1.1.18.** Let $R$ be a domain and let $\mathscr{C}$ be a full subcategory of $\mathrm{mod}(R)$, such that $\mathscr{C}$ has the KRS property. Then for every object $M$ of $\mathscr{C}$ we have a unique decomposition into indecomposable objects

$$M = M_1^{\oplus n_1} \oplus \cdots \oplus M_r^{\oplus n_r},$$

for some integers $n_i$. Then, $\mathrm{rank}_R M = \sum_i n_i \, \mathrm{rank}_R M_i$.

When $R$ is a graded commutative ring, it is convenient to introduce a graded version of the free rank. We recall that by graded ring we always mean $\mathbb{N}$-graded ring.

**Definition 1.1.19.** Let $R$ be a graded ring and let $M$ be a finitely generated graded $R$-module. We define the *graded free rank* of $M$ as

$$\mathrm{frk}_R^{\mathrm{gr}}(M) := \max\{n : \exists \text{ a homogeneous of degree 0 split surjection } \varphi : M \twoheadrightarrow F,$$
$$\text{with } F \text{ free graded } R\text{-module of rank } n\}.$$

The main reason to introduce this definition is that for a graded ring $R$ the category $\mathrm{mod}(R)$ is rarely a Krull-Schmidt category, even if $R$ is complete or regular (see Example 1.1.9). To hope to recover a unique decomposition of an $R$-module into indecomposable objects, one should rather work with the category $\mathrm{mod}_{\mathbb{Z}}(R)$ of finitely generated graded $R$-modules, whose maps are $R$-linear homomorphisms of degree 0. Thus, Definition 1.1.19 looks more natural in this category.

One should remind that in $\mathrm{mod}_{\mathbb{Z}}(R)$ two $R$-modules $R(a)$ and $R(b)$ with $a \neq b$ are not isomorphic. However, all modules of the form $R(a)$ for some integer $a$ contribute to the graded free rank part of the module. For example, the modules $R(1) \oplus R(-2)$ and $R^2$ are two non-isomorphic graded $R$-modules of free rank 2. Therefore the free $R$-module $F$ of Definition 1.1.19 can be uniquely written as

$$F = \bigoplus_{i=1}^{n} R(a_i),$$

for some integers $a_i$.

Finally, we point out that for a finitely generated graded module $M$ over a graded ring $R$, the following inequality holds

$$\mathrm{frk}_R^{\mathrm{gr}} M \leq \mathrm{frk}_R M.$$

## 1.2. Syzygies

We consider a Noetherian ring $R$ and a left $R$-module $M$. We want to define the $n$-th syzygy module of $M$.

**Definition 1.2.1.** Two left $R$-modules $M_1$ and $M_2$ are *projective equivalent* if there exist two projective left $R$-modules $P_1$ and $P_2$ such that $M_1 \oplus P_1 \cong M_2 \oplus P_2$.

We point out that an $R$-module is equivalent to a projective module if and only if it is projective.

The following lemma is a generalization of Schanuel's Lemma [Rot09, Proposition 3.12].

**Lemma 1.2.2** (Schanuel)**.** *Let $M$ be a left $R$-module and let*

$$0 \to K_n \to P_{n-1} \to \cdots P_1 \to P_0 \to M \to 0$$
$$0 \to K_n' \to P_{n-1}' \to \cdots P_1' \to P_0' \to M \to 0,$$

*be two exact sequences in* $\mathrm{mod}(R)$*, where $P_i$ and $P_i'$ are projective $R$-modules, then $K_n$ and $K_n'$ are projective equivalent.*

This leads us to define for each natural number $n$ the *$n$-th syzygy module of $M$* as

$$\Omega_R^n(M) := K_n.$$

The module $\Omega_R^n(M)$ is defined up to projective equivalence. In fact, if we denote by $\underline{\mathrm{mod}}(R)$ the stable category of left $R$-modules, then $\Omega_R^n(-)$ is a functor $\underline{\mathrm{mod}}(R) \to \underline{\mathrm{mod}}(R)$. We recall that the category $\underline{\mathrm{mod}}(R)$ has the same objects of $\mathrm{mod}(R)$, and a morphism $f : M \to N$ in $\underline{\mathrm{mod}}(R)$ is an $R$-linear map modulo the relation that $f \sim g$ if $f - g$ factors through a projective module.

Using syzygy modules, we can define the *projective dimension* of a left $R$-module $M$ as

$$\mathrm{proj.dim}_R(M) := \inf\{n : \Omega_R^n(M) \text{ is projective}\}.$$

Schanuel's Lemma and its generalization ensure us that this is well-defined.

We observe that $\Omega_R^n(M) = 0$ if $\mathrm{proj.dim}_R(M) < n$ and $\Omega_R^0(M) = M$, so $M$ is projective if and only if $\mathrm{proj.dim}_R(M) = 0$.

### 1.2.1. Syzygies over local rings

Over a Noetherian commutative local ring $(R, \mathfrak{m}, k)$ working with syzygies is particularly nice. One can define the syzygy modules uniquely, without the projective equivalence condition, using the following construction.

Let $M$ be a finitely generated $R$-module. We choose $x_1, \ldots, x_{\beta_0}$ a minimal system of generators of $M$. By Nakayama's Lemma the number $\beta_0$ is unique, precisely it is the $k$-vector space dimension

$$\beta_0 = \mu(M) := \dim_k(M \otimes_R k).$$

We fix an epimorphism $\varphi_0 : R^{\beta_0} \to M$ sending the canonical basis $\{e_1, \ldots, e_{\beta_0}\}$ of $R^{\beta_0}$ to $\{x_1, \ldots, x_{\beta_0}\}$. The module $\ker \varphi_0$ is again finitely generated, so we set $\beta_1 := \mu(\ker \varphi_0)$ and we continue similarly defining an epimorphism $\varphi_1 : R^{\beta_1} \to \ker \varphi_0$. Proceding in this manner we construct a minimal free resolution

$$\mathscr{F} : \cdots R^{\beta_n} \xrightarrow{\varphi_n} R^{\beta_{n-1}} \to \cdots \to R^{\beta_1} \xrightarrow{\varphi_1} R^{\beta_0} \xrightarrow{\varphi_0} M \xrightarrow{\varphi_{-1}} 0, \quad (1.1)$$

which is uniquely determined by $M$ up to isomorphism. The numbers $\beta_i(M) := \beta_i$ are called the *Betti numbers* of $M$.

**Definition 1.2.3.** Let $(R, \mathfrak{m}, k)$ be a Noetherian local ring and $M$ a finitely generated $R$-module with minimal free resolution $\mathscr{F}$ as in (1.1). For every natural number $n$ the *n-th syzygy module of $M$* is

$$\mathrm{Syz}_R^n(M) := \ker \varphi_{n-1}.$$

Since finite free modules are projective this definition is consistent with the previous one, in other words $\mathrm{Syz}_R^n(M)$ is a representative of the projective equivalence class of $\Omega_R^n(M)$.

## 1.3. Reflexive modules

Let $R$ be a Noetherian ring. We do not assume that $R$ is commutative in this section, so by $R$-module we mean left $R$-module. We denote by $(-)^*$ the contravariant functor $\mathrm{Hom}_R(-, R) : \mathrm{mod}(R) \to \mathrm{mod}(R)$. Then, for every left $R$-module $M$ we have a canonical map

$$\lambda_M : M \to M^{**}$$

such that $\lambda_M(m)(f) = f(m)$ for every $f \in M^*$. The homomorphism $\lambda_M$ is in general not injective nor surjective. We say that

- $M$ is *torsionless* if $\lambda_M$ is injective;
- $M$ is *reflexive* if $\lambda_M$ is bijective.

**Example 1.3.1.** The ring $R$ is a reflexive $R$-module. Every finite free $R$-module is reflexive.

**Lemma 1.3.2.** *Let $R$ be a Noetherian domain and let $M$ be a finite $R$-module. Then $M$ is reflexive if and only if $M_{\mathfrak{m}}$ is a reflexive $R_{\mathfrak{m}}$-module for all maximal ideals $\mathfrak{m}$ of $R$.*

For every $R$-module $M$, the module

$$T(M) := \{m \in M : \exists r \in R \text{ non-zero-divisor, s.t. } r \cdot m = 0\}$$

is called the *torsion submodule* of $M$. We say that

- $M$ is *torsion-free* if $T(M) = 0$;
- $M$ is a *torsion module* if $T(M) = M$.

In other words $M$ is torsion-free if and only if all zero-divisors on $R$ are zero-divisors on $M$.

Torsionless and reflexive modules are very closely related to syzygies

**Proposition 1.3.3.** *Let $R$ be a Noetherian ring and let $M$ be a finitely generated $R$-module. Then the following facts hold.*

*1) $M$ is torsionless if and only if it is a first syzygy.*
*2) If $M$ is reflexive then it is a second syzygy.*
*3) Reflexive $\Rightarrow$ torsionless $\Rightarrow$ torsion-free.*
*4) If $R$ is a domain and $M$ is torsion-free, then $M$ is torsionless.*

*Proof: 1)* First assume that $M$ is torsionless and let $f_1, \ldots, f_n$ be a set of generators for $M^*$. We define a map $\varphi : M \to R^n$ by $\varphi(x) = (f_1(x), \ldots, f_n(x))$. If $\varphi(x) = 0$ then $f(x) = 0$ for every $f \in M^*$, so $\lambda_M(x) = 0$ which implies $x = 0$, since $\lambda_M$ is injective. Hence $\varphi$ is injective and $M$ is the first syzygy of the module $\mathrm{Coker}(\varphi)$. On the other hand, let $M$ be the first syzygy of an $R$-module $N$, that is

$$0 \to M \xrightarrow{\psi} R^m \to N \to 0.$$

Let $x$ be a non-zero element of $M$, then there exists an index $j$ such that $\pi_j \psi(x) \neq 0$, where $\pi_j : R^m \to R$ is the $j$-th canonical projection. Then $\pi_j \psi \in M^*$ does not vanish on $x$, so $x \notin \ker \lambda_M$.

*2)* We consider a free presentation of the dual $F_2 \to F_1 \to M^* \to 0$. We apply the left exact functor $\mathrm{Hom}_R(-, R)$ and we get

$$0 \to M^{**} \to F_1^* \to F_2^*.$$

Since $M^{**} \cong M$, this shows that $M$ is a second syzygy.

*3)* Reflexive modules are clearly torsionless from the definition. For the second implication, we have that a torsionless module $M$ is a submodule of a free module by *1)*. So non-zero-divisors on $R$ are zero-divisors on $M$.

*4)* Let $f_1, \ldots, f_m$ be a system of generators of $M$, and let $S = R \setminus \{0\}$. Then the localization $M_S$ is a finite dimensional vector space over the quotient ring $Q(R) = S^{-1}R$. Thus, it can be embedded in $Q(R)^n$ for some positive integer $n$. We obtain an injective homomorphism of $R$-modules $\psi : M \to Q(R)^n$ given by the composition

$$M \hookrightarrow M_S \hookrightarrow Q(R)^n.$$

Actually the map $\psi$ is given by $f_i \mapsto \frac{1}{s_i}(r_{i,1}, \ldots, r_{i,n})$ for some $s_i \in S$ and $r_{i,j} \in R$. Let $s = s_1 \cdots s_m$, and consider the map $\varphi$ given by the composition

$$M \xrightarrow{s} M \xrightarrow{\psi} Q(R)^n.$$

The image of $\varphi$ is contained in $R^n$, and since multiplication by $s$ is injective on $M$, $\varphi : M \to R^n$ is also injective. Therefore $M$ is a first syzygy, hence torsionless by *1)*. $\square$

**Example 1.3.4.** The last implication of Proposition 1.3.3 is no longer true if we remove the finitely generated hypothesis on $M$. For example $\mathbb{Q}$ is a torsion-free $\mathbb{Z}$-module which is not torsionless.

From now on, we will restrict our attention almost exclusively to finitely generated modules. We denote by $\mathrm{Ref}(R)$ the category whose objects are finitely generated reflexive left $R$-modules and whose maps are just $R$-linear homomorphisms between them. Then $\mathrm{Ref}(R)$ is a full subcategory of $\mathrm{mod}(R)$.

The converse of *2)* in Proposition 1.3.3 does not hold in general, however we have the following characterization of reflexive modules in terms of syzygies.

**Lemma 1.3.5.** *Let $R$ be a Noetherian domain and let $M$ be a finitely generated $R$-module. Then the following facts are equivalent.*

*1) $M$ is reflexive.*

*2) There exists a short exact sequence $0 \to M \to F \to N \to 0$ with $F$ finite free and $N$ torsion-free.*

*Proof:* We consider a short exact sequence as in *2)* and we apply the functor $(-)^{**}$. We get a commutative diagram

$$
\begin{array}{ccccccccc}
0 & \longrightarrow & M & \longrightarrow & F & \longrightarrow & N & \longrightarrow & 0 \\
& & \downarrow & & \downarrow & & \downarrow & & \\
0 & \longrightarrow & M^{**} & \longrightarrow & F^{**} & \longrightarrow & N^{**} & \longrightarrow & 0.
\end{array}
$$

The first row of the diagram is exact, while the second one may not be. Anyway the map $F \to F^{**}$ is an isomorphism, since $F$ is finite free, hence reflexive, and the map $N \to N^{**}$ is injective since $N$ is torsionless, being torsion-free and finitely generated over a domain. The module $M^{**}$ is torsion-free, so the map $M^{**} \to F^{**}$ is injective. Then by diagram chasing we obtain that $M$ is reflexive.

Conversely, assume that $M$ is reflexive. As in the proof of Proposition 1.3.3 we have an exact sequence $0 \to M^{**} \cong M \xrightarrow{\varphi} F_1 \to F_2$, with $F_1$ and $F_2$ finite free. It is clear that the module $\mathrm{Coker}\,\varphi$ is torsion-free. $\square$

**Lemma 1.3.6.** *Let $R$ be a Noetherian ring and let $M$ be a finitely generated $R$-module. Then $M^*$ is a reflexive module.*

*Proof:* Consider the canonical map $\lambda_{M^*} : M^* \to M^{***}$, and let $f \neq 0$ be an element of $M^*$. We prove that $f^{**} := \lambda_{M^*}(f)$ is not zero. Since $f \neq 0$, there exists $x \in M$ such that $f(x) \neq 0$. Let $\mathrm{Ev}_x \in M^{**}$ be the evaluation at $x$, we have $f^{**}(\mathrm{Ev}_x) = f(x) \neq 0$, so $f^{**} \neq 0$ and $\lambda_{M^*}$ is injective.

For surjectivity, let $g \in M^{***}$. We consider the map $h = g \circ \lambda_M : M \to R$. Then $h \in M^*$, and the following diagram

$$
\begin{array}{ccc}
M^{**} & \xrightarrow{\ g\ } & R \\
{\scriptstyle \lambda_M}\big\uparrow & \nearrow{\scriptstyle h} & \\
M & &
\end{array}
$$

shows that $\lambda_{M^*}(h) = g$. $\square$

**Proposition 1.3.7.** *Let $R$ be a Noetherian ring, let $M_1$ and $M_2$ be reflexive $R$-modules and let $N$ be an $R$-module. Then the following facts hold.*

*1)* $M_1 \oplus M_2$ *is reflexive.*
*2)* $\mathrm{Hom}_R(N, M_1)$ *is reflexive.*
*3)* *Let $N$ be a direct summand of $M_1$, then $N$ is reflexive.*

Since projective modules are direct summands of finite free modules, which are reflexive, we immediately get the following.

**Corollary 1.3.8.** *Let $R$ be a Noetherian ring and let $M$ be a finitely generated projective $R$-module. Then $M$ is reflexive.*

We recall that $\mathrm{proj}(R)$ denotes the category of finitely generated projective left $R$-modules. From the previous corollary we have the following inclusions of categories:

$$\mathrm{proj}(R) \subseteq \mathrm{Ref}(R) \subseteq \mathrm{mod}(R).$$

**Example 1.3.9.** Reflexive modules are in general not projective. Consider $R = k[x^2, xy, y^2]$, the module $M = (x^2, xy)$ is reflexive, but not projective.

**Definition 1.3.10.** Let $R$ be a Noetherian domain and let $M$ be a finitely generated $R$-module. The module $M^{**} = \mathrm{Hom}_R(\mathrm{Hom}_R(M, R), R)$ is called the *reflexive hull* of $M$.

The name reflexive hull makes sense, since by Lemma 1.3.6 $M^{**}$ is a reflexive module. Thus, we have a functor

$$(-)^{**} : \mathrm{mod}(R) \to \mathrm{Ref}(R),$$

which is a left adjoint to the inclusion functor $\mathrm{Ref}(R) \hookrightarrow \mathrm{mod}(R)$. To see this it is enough to see that if $N \in \mathrm{Ref}(R)$ every $R$-module map $f : M \to N$ factors through the reflexive hull

of $M$, that is

$$
\begin{array}{ccc}
M & \xrightarrow{\ f\ } & N \\
{\scriptstyle \lambda_M}\Big\downarrow & \nearrow & \\
M^{**}. & &
\end{array}
$$

Given two reflexive $R$-modules $M$ and $N$ the tensor product $M \otimes_R N$ does not need to be reflexive. For example consider $R = k[\![x, y]\!]/(xy)$ and $M = (x, y)$. Then $M$ is reflexive, but the tensor product $M \otimes_R M$ is not even torsionless.

For this reason, after Auslander [Aus86b] it is convenient to introduce the following definition.

**Definition 1.3.11.** Let $M, N \in \mathrm{Ref}(R)$. The *reflexive tensor product* of $M$ and $N$ is the module

$$
M \boxtimes_R N := (M \otimes_R N)^{**}.
$$

Now for each $M, N, L \in \mathrm{Ref}(R)$ we have a natural homomorphism $M \otimes_R N \to M \boxtimes_R N$, which induces an $R$-isomorphism $\mathrm{Hom}_R(M \boxtimes_R N, L) \xrightarrow{\simeq} \mathrm{Hom}_R(M \otimes_R N, L)$, which is functorial in $M$, $N$ and $L$. From the usual hom-tensor adjunction $\mathrm{Hom}_R(M \otimes_R N, L) \cong \mathrm{Hom}_R(M, \mathrm{Hom}_R(N, L))$ we get the isomorphism

$$
\mathrm{Hom}_R(M \boxtimes_R N, L) \xrightarrow{\simeq} \mathrm{Hom}_R(M, \mathrm{Hom}_R(N, L)).
$$

In fact we obtained an adjunction in the category $\mathrm{Ref}(R)$.

**Proposition 1.3.12.** *Let $M$ be a reflexive $R$-module. Then the reflexive tensor product functor $M \boxtimes_R -$ is left adjoint to the* $\mathrm{Hom}$*-functor* $\mathrm{Hom}_R(M, -)$*.*

**Proposition 1.3.13.** *Let $M_1, M_2, M_3 \in \mathrm{Ref}(R)$. Then*

$$
M_1 \boxtimes_R M_2 \cong M_2 \boxtimes_R M_1, \quad (M_1 \boxtimes_R M_2) \boxtimes_R M_3 \cong M_1 \boxtimes_R (M_2 \boxtimes_R M_3).
$$

As for the tensor product, also the $q$-th symmetric power $\mathrm{Sym}_R^q(M)$ of a reflexive module need not to be reflexive. With the same spirit of the reflexive tensor product, we give the following definition.

**Definition 1.3.14.** Let $M \in \mathrm{Ref}(R)$ and let $q$ be a natural number. The *$q$-th reflexive symmetric power* of $M$ is the $R$-module $\left(\mathrm{Sym}_R^q(M)\right)^{**}$.

## 1.4. Depth and Serre's conditions

Let $(R, \mathfrak{m}, k)$ be a Noetherian commutative local ring, we recall the definition of depth of an $R$-module.

**Definition 1.4.1.** A sequence $x_1, \ldots, x_n$ of elements in $\mathfrak{m}$ is called an *M-regular sequence* if $x_1$ is a non-zero-divisor in $M$, $x_{i+1}$ is a non-zero-divisor in $M/(x_1, \ldots, x_i)M$ for any $i \in \{1, \ldots, n-1\}$, and $(x_1, \ldots, x_n)M \neq M$. The *depth* of $M$ is denoted by $\mathrm{depth}_R(M)$ and is equal to the maximal lenght of an *M*-regular sequence in $\mathfrak{m}$.

**Remark 1.4.2.** The definition of regular sequence may depend a priori on the order of the elements, however in this setting this is not the case. In fact, if $x_1, \ldots, x_n$ is an *M*-regular sequence, then every permutation $x_{\pi(1)}, \ldots, x_{\pi(n)}$ is an *M*-sequence (cf. [BH98, Exercise 1.2.21]).

When the ring $R$ is clear from the context we will denote the depth of $M$ simply by $\mathrm{depth}\, M$.

The following well-known inequality holds for any finitely generated $R$-module $M$:

$$\mathrm{depth}\, M \leq \dim M \leq \dim R,$$

where $\dim R$ denotes the Krull dimension of $R$ and $\dim M$ the Krull dimension of $M$, defined as $\dim M = \dim(R/\mathrm{Ann}_R(M))$.

The projective dimension and the depth of a module are related by the Auslander-Buchsbaum formula (cf. [BH98, Theorem 1.3.3]).

**Theorem 1.4.3** (Auslander-Buchsbaum)**.** *Let $(R, \mathfrak{m}, k)$ be a Noetherian local ring and $M \neq 0$ a finitely generated $R$-module such that $\mathrm{proj.dim}_R(M) < \infty$. Then*

$$\mathrm{proj.dim}_R(M) + \mathrm{depth}_R(M) = \mathrm{depth}\, R.$$

Moreover, we recall that in this setting projective modules and free modules coincide, and that the finiteness of the projective dimension of the residue field $k$ characterizes regularity (cf. [BH98, Theorem 2.2.7]).

**Theorem 1.4.4.** *Let $(R, \mathfrak{m}, k)$ be a Noetherian local ring and let $M$ be a finitely generated $R$-module. Then the following facts are equivalent.*

*1) $M$ is free.*
*2) $M$ is projective.*
*3) $M$ is flat.*

**Remark 1.4.5.** The implications *1) $\Rightarrow$ 2) $\Rightarrow$ 3)* hold in every Noetherian ring.

**Theorem 1.4.6** (Auslander-Buchsbaum-Serre)**.** *Let $(R, \mathfrak{m}, k)$ be a Noetherian local ring. Then the following facts are equivalent.*

*1) $R$ is regular.*
*2) $\mathrm{proj.dim}_R(M) < \infty$ for every finitely generated $R$-module $M$.*
*3) $\mathrm{proj.dim}_R(k) < \infty$.*

The depth of an $R$-module can be characterized in terms of local cohomology and Ext-functor. Local cohomology is an important tool in commutative algebra and algebraic geometry. We recall here the definition and some properties in our particular setting. For a more general introduction and proofs of these facts we recommend the books of Brodmann and Sharp [BS98], and of Iyengar et al. [ILLMMSW07].

**Definition 1.4.7.** Let $(R, \mathfrak{m}, k)$ be a Noetherian local ring, then the functor $\Gamma_{\mathfrak{m}} : \mathrm{mod}(R) \to \mathrm{mod}(R)$ given by

$$M \mapsto \Gamma_{\mathfrak{m}}(M) := \{x \in M : \mathfrak{m}^t x = 0 \text{ for some } t > 0\}.$$

is left exact. Its derived functors are the local cohomology functors, that is $H_{\mathfrak{m}}^i(M) := R^i \Gamma_{\mathfrak{m}}(M)$ for all $R$-modules $M$.

**Lemma 1.4.8.** *Let $M$ be an $R$-module. Then*

$$\mathrm{depth}_R(M) = \inf\{i \in \mathbb{N} : \mathrm{Ext}_R^i(R/\mathfrak{m}, M) \neq 0\}$$
$$= \inf\{i \in \mathbb{N} : H_{\mathfrak{m}}^i(M) \neq 0\}.$$

From the long exact sequence of Ext-modules it is easy to prove the following well-known result.

**Lemma 1.4.9** (Depth Lemma)**.** *Let $0 \to L \to M \to N \to 0$ be a short exact sequence of finitely generated $R$-modules. Then the following facts hold.*

*1) If* $\mathrm{depth}\, N < \mathrm{depth}\, M$, *then* $\mathrm{depth}\, L = \mathrm{depth}\, N + 1$.
*2)* $\mathrm{depth}\, L \geq \min\{\mathrm{depth}\, M, \mathrm{depth}\, N\}$.
*3)* $\mathrm{depth}\, M \geq \min\{\mathrm{depth}\, L, \mathrm{depth}\, N\}$.

**Definition 1.4.10.** Let $M$ be a finitely generated module over a Noetherian ring $R$ and let $n$ be a natural number. We say that $M$ satisfies *Serre's condition* $(S_n)$ provided

$$\mathrm{depth}_{R_{\mathfrak{p}}}(M_{\mathfrak{p}}) \geq \min\{n, \dim R_{\mathfrak{p}}\},$$

for every $\mathfrak{p} \in \mathrm{Spec}\, R$.

**Remark 1.4.11.** We point out that this definition agrees with the one given by Leuschke and Wiegand [LW12, p. 310], but differs from the one given by Bruns and Herzog [BH98, p. 63], who require that $\mathrm{depth}_{R_{\mathfrak{p}}}(M_{\mathfrak{p}}) \geq \min\{n, \dim M_{\mathfrak{p}}\}$. Of course if the module $M$ is the ring $R$ then the two definitions coincide. The main motivation for our choice should be found in the upcoming Lemma 1.4.14.

The following theorem is known as *Serre's criterion for normality* (cf. [Ser00, Theorem IV.D.11]).

**Theorem 1.4.12** (Serre)**.** *Let $R$ be a Noetherian ring. Then the following facts are equivalent.*

*1)* *R is a normal domain.*
*2)* *R satisfies* $(S_2)$ *and* $R_{\mathfrak{p}}$ *is a regular local ring for each prime ideal* $\mathfrak{p}$ *such that* $\mathrm{ht}(\mathfrak{p}) \leq 1$.

**Lemma 1.4.13.** *Let M and N be finitely generated modules over a local ring* $(R, \mathfrak{m}, k)$. *Then* $\mathrm{depth}\,\mathrm{Hom}_R(M, N) \geq \min\{2, \mathrm{depth}\,N\}$

*Proof:* For depth $N = 0$ there is nothing to prove. If depth $N = 1$, then there exists $r \in \mathfrak{m}$ such that $r \cdot n \neq 0$ for all non-zero $n \in N$. It follows that if $\varphi \in \mathrm{Hom}_R(M, N)$ is non-zero, then also $r \cdot \varphi$ is non-zero. Thus $r$ is a regular element for $\mathrm{Hom}_R(M, N)$, and the claim is proved.

Now assume that depth $N \geq 2$ and let $x, y \in \mathfrak{m}$ be a regular sequence for $N$. The previous observation implies that $x$ is a regular element for $\mathrm{Hom}_R(M, N)$, so it remains to show that $y$ is non-zero-divisor for $\mathrm{Hom}_R(M, N)/x\,\mathrm{Hom}_R(M, N)$.

Let $\varphi \in \mathrm{Hom}_R(M, N)$ such that $y \cdot [\varphi] = [0]$ in $\mathrm{Hom}_R(M, N)/x\,\mathrm{Hom}_R(M, N)$, we prove that $[\varphi] = [0]$. Since $y \cdot [\varphi] = [0]$, there exists $\psi \in \mathrm{Hom}_R(M, N)$ such that $y \cdot \varphi = x \cdot \psi$. It follows that for every $m \in M$ we have $y \cdot \varphi(m) = x \cdot \psi(m)$, which is 0 in $N/xN$. But $y$ is regular for $N/xN$, then $[\varphi(m)] = [0]$ in $N/xN$. Therefore there exists $n_m \in M$ such that $\varphi(m) = n_m$. We define a function $f \in \mathrm{Hom}_R(M, N)$ by $f(m) = n_m$, then we have $\varphi = x \cdot f$, which implies that $[\varphi] = [0]$. □

The following lemma is due to Auslander and Buchsbaum [AB59], the reader may consult also [LW12, Lemma 5.11].

**Lemma 1.4.14** (Auslander-Buchsbaum)**.** *Let R be a commutative Noetherian ring and let* $f : M \to N$ *be a homomorphism of finitely generated R-modules such that M satisfies* $(S_2)$ *and N satisfies* $(S_1)$. *If the map* $f_{\mathfrak{p}} : M_{\mathfrak{p}} \to N_{\mathfrak{p}}$ *is*

* *injective for every minimal prime* $\mathfrak{p}$*;*
* *surjective for every prime ideal* $\mathfrak{p}$ *of height* $\leq 1$*;*

*then f is an isomorphism.*

The following two results show how the reflexive property can be characterized using the depth and Serre's conditions.

**Lemma 1.4.15.** *Let R be a Noetherian domain and let M be a finitely generated R-module. Then the following facts are equivalent.*

*1)* *M is reflexive.*
*2)* *For every prime ideal* $\mathfrak{p}$ *of R one of the following happens*
    *a)* $M_{\mathfrak{p}}$ *is a reflexive* $R_{\mathfrak{p}}$*-module, or*
    *b)* $\mathrm{depth}(R_{\mathfrak{p}}) \geq 2$ *and* $\mathrm{depth}(M_{\mathfrak{p}}) \geq 2$.

*Proof:* If *1)* is true, then *2)* follows from Lemma 1.3.2 and Lemma 1.4.13 applied to $R_{\mathfrak{p}}^{**} \cong R_{\mathfrak{p}}$ and $M_{\mathfrak{p}}^{**} \cong M_{\mathfrak{p}}$.

Conversely, assume that *2)* holds. We fix $N = M^{**}$ and we consider the canonical map $\lambda : M \to N$. We want to apply Lemma 1.4.14 to $\lambda$. It is easy to check that conditions *a)* and *b)* imply that $M$ and $R$ are $(S_2)$. Let $\mathfrak{p}$ be a prime ideal of $R$. If *a)* holds, then the map $\lambda_{\mathfrak{p}} : M_{\mathfrak{p}} \to N_{\mathfrak{p}}$ is an isomorphism, otherwise depth $R_{\mathfrak{p}} \geq 2$, which implies that also the height of $\mathfrak{p}$ is $\geq 2$. Therefore by Lemma 1.4.14, $\lambda$ is an isomorphism, that is $M$ is reflexive. $\square$

**Theorem 1.4.16.** *Let $R$ be a Noetherian normal domain and let $M$ be a finitely generated $R$-module. Then the following facts are equivalent.*

*1)  M is reflexive.*
*2)  M is torsion-free and has property $(S_2)$.*

*Proof:* We observe that in both cases $M$ is torsion-free. Let $\mathfrak{p}$ be a prime ideal of $R$. If the height of $\mathfrak{p}$ is $\leq 1$, then by Serre's criterion for normality (Theorem 1.4.12) $R_{\mathfrak{p}}$ is a regular local ring of dimension $\leq 1$. Thus, the torsion-free module $M_{\mathfrak{p}}$ is finite free, hence reflexive. If the height of $\mathfrak{p}$ is $\geq 2$, then again by Theorem 1.4.12 we have depth $R_{\mathfrak{p}} \geq 2$. It follows by Lemma 1.4.15 that $M$ is reflexive if and only if depth $M_{\mathfrak{p}} \geq 2$ for primes of height $\geq 2$, which in this case is equivalent to the $(S_2)$ condition.  $\square$

## 1.5.  Maximal Cohen-Macaulay modules

We introduce another important class of modules: maximal Cohen-Macaulay modules (MCM for short). We will see how these modules are related to reflexive modules and to syzygies. Then we will investigate some of their properties, which will be useful later on. Every ring in this section is commutative and Noetherian.

**Definition 1.5.1.** Let $(R, \mathfrak{m})$ be a Noetherian local ring. An $R$-module $M$ is *Cohen-Macaulay* (abbr. CM) if depth$_R(M) = \dim M$.
An $R$-module $M$ is *maximal Cohen-Macaulay* (abbr. MCM) if depth$_R(M) = \dim R$.
The ring $R$ is *Cohen-Macaulay ring* if it is a Cohen-Macaulay $R$-module.

Even if we will work mainly with local rings, CM modules may be defined also in the non-local setting thanks to the following proposition (cf. [Eis94, Proposition 18.8]).

**Proposition 1.5.2.** *Let $(R, \mathfrak{m})$ be a Noetherian local ring and $M$ a CM $R$-module, then $M_{\mathfrak{p}}$ is CM over $R_{\mathfrak{p}}$ for all prime ideals $\mathfrak{p}$ of $R$.*

**Definition 1.5.3.** Let $R$ be a Noetherian ring and $M$ an $R$-module. $M$ is CM (resp. MCM) if $M_{\mathfrak{m}}$ is CM (MCM) over $R_{\mathfrak{m}}$ for all maximal ideals $\mathfrak{m}$ of $R$.  $R$ is a Cohen-Macaulay ring if it is a CM module.

**Example 1.5.4.** Some easy examples of CM rings and modules are the following.

- Regular rings are Cohen-Macaulay.

- An Artinian ring $R$ (i.e. $\dim R = 0$) is Cohen-Macaulay.
- A domain of Krull dimension 1 is Cohen-Macaulay.
- The two-dimensional ring $k[\![x, y, z]\!]/(x^4 - yz)$ is Cohen-Macaulay.
- A finite free $R$-module is MCM.

Given a Noetherian ring $R$ we denote by $\mathrm{MCM}(R)$ the full-subcategory of $\mathrm{mod}(R)$ consisting of all maximal Cohen-Macaulay modules.

Thanks to the following result of Grothendieck (cf. [BH98, Theorem 3.5.7]) we can give a homological characterization of the Cohen-Macaulay property. Notice that the vanishing of $H^i_{\mathfrak{m}}(M)$ for $i < \mathrm{depth}_R(M)$ is a consequence of Lemma 1.4.8.

**Theorem 1.5.5** (Grothendieck)**.** *Let $(R, \mathfrak{m}, k)$ be a Noetherian local ring and $M$ a finitely generated $R$-module of depth $t$ and dimension $d$. Then the following facts hold.*

*1)* $H^i_{\mathfrak{m}}(M) = 0$ *for $i < t$ and $i > d$.*
*2)* $H^t_{\mathfrak{m}}(M) \neq 0$ *and $H^d_{\mathfrak{m}}(M) \neq 0$.*

**Corollary 1.5.6.** *Let $(R, \mathfrak{m}, k)$ be a Noetherian local ring and $M$ an $R$-module. Then the following facts are equivalent.*

*1)* $M$ *is MCM over $R$.*
*2)* $\mathrm{Ext}^i_R(k, M) = 0$ *for $i < \dim R$.*
*3)* $H^i_{\mathfrak{m}}(M) = 0$ *for $i \neq \dim R$.*

We use Corollary 1.5.6 to prove some basic results on MCM modules.

**Proposition 1.5.7.** *Let $(R, \mathfrak{m}) \subseteq (T, \mathfrak{n})$ be a finite extension of local Noetherian rings and $M$ a Noetherian $T$-module. Then $M$ is MCM over $R$ if and only if it is MCM over $T$.*

*Proof:* The forgetful functor $\mathrm{mod}(T) \to \mathrm{mod}(R)$ is exact, therefore $\Gamma_{\mathfrak{n}}(M) \cong \Gamma_{\mathfrak{m}}(M)$ as $R$-modules. Hence also $H^i_{\mathfrak{n}}(M) \cong H^i_{\mathfrak{m}}(M)$ for all $i \geq 0$, which concludes the proof. $\square$

Under the assumptions of Proposition 1.5.7, if the ring $R$ is regular the map $R \to T$ is called a *Noether normalization* of $T$. So Proposition 1.5.7 states that the property of being MCM is invariant under Noether normalization.

The following is an immediate consequence of the Depth Lemma

**Proposition 1.5.8.** *Let $(R, \mathfrak{m}, k)$ be a Noetherian local ring and let $0 \to L \to M \to N \to 0$ be a short exact sequence of MCM $R$-modules. Then the following facts hold.*

*1)* *If $L$ and $N$ are MCM, then so is $M$.*
*2)* *If $M$ and $N$ are MCM, then so is $L$.*

**Remark 1.5.9.** It is not necessary that if $L$ and $M$ are MCM then so is $N$. As a counterexample consider the power series ring in one variable $R = k[\![x]\!]$ and the following short exact sequence of $R$-modules

$$0 \to R \xrightarrow{\cdot x^2} R \to R/(x^2) \to 0.$$

As a consequence of the Auslander-Buchsbaum formula (Theorem 1.4.3) the theory of MCM modules over regular local rings is trivial.

**Proposition 1.5.10.** *Let $(R, \mathfrak{m}, k)$ be a Noetherian local regular ring and $M$ a finitely generated $R$-module. Then $M$ is maximal Cohen-Macaulay if and only if it is free.*

*Proof:* Since $M$ is a finitely generated module over a regular ring it has finite projective dimension by Theorem 1.4.6. Then the Auslander-Buchsbaum formula

$$\operatorname{proj.dim}_R(M) + \operatorname{depth}_R(M) = \operatorname{depth} R = \dim R$$

implies that $M$ is free if and only if $\operatorname{proj.dim}_R(M) = 0$ if and only if $\dim R - \operatorname{depth}_R(M) = 0$, which is the definition of MCM. $\square$

Combining this result with Corollary 1.5.7 we get the following.

**Corollary 1.5.11.** *Let $R \to T$ be a Noether normalization of the local ring $T$. Then a $T$-module $M$ is MCM over $T$ if and only if it is $R$-free.*

We consider now maximal Cohen-Macaulay modules over rings of low dimension. The zero-dimensional case is trivial, since every finitely generated module is automatically MCM. In dimension one and two the situation is clarified by the following proposition.

**Proposition 1.5.12.** *Let $(R, \mathfrak{m}, k)$ be a Noetherian local ring and let $M$ be a finitely generated $R$-module. Then the following facts hold.*

1) *If $R$ is reduced of dimension one, then $M$ is MCM if and only if it is torsion-free.*
2) *If $R$ is a normal domain of dimension two, then $M$ is MCM if and only if it is reflexive, that is $\operatorname{Ref}(R) = \operatorname{MCM}(R)$.*
3) *If $R$ is a normal domain of dimension $\geq 3$, then every MCM $R$-module is reflexive. However it is not true in general that reflexive modules are MCM.*

*Proof:* It follows from the definition that $M$ is torsion-free if and only if $\operatorname{Hom}_R(R/\mathfrak{m}, M) = 0$, which is equivalent to $\operatorname{depth}_R(M) \geq 1$. Since $\operatorname{depth}_R(M) \leq \dim M \leq \dim R$, if the ring has dimension 1 the condition $\operatorname{depth}_R(M) \geq 1$ is equivalent to the module being maximal Cohen-Macaulay. This proves *1)*.

Now we prove that for normal domains of dimension $\geq 2$ MCM modules are reflexive. Let $M$ be an MCM $R$-module, then $M_\mathfrak{p}$ is MCM over $R_\mathfrak{p}$ for every prime ideal $\mathfrak{p}$. If $\operatorname{ht}(\mathfrak{p}) \leq 1$, then by Serre's criterion (Theorem 1.4.12) $R_\mathfrak{p}$ is regular, so $M_\mathfrak{p}$ is free, in particular it is reflexive. If $\operatorname{ht}(\mathfrak{p}) \geq 2$ then again by Serre's criterion $\operatorname{depth} R_\mathfrak{p} \geq 2$ so $M_\mathfrak{p}$ is reflexive. Then Lemma 1.3.2 implies that $M$ is reflexive.

Conversely, assume that $R$ has dimension 2 and that $M$ is reflexive. Let $\mathfrak{p}$ be a prime ideal of $R$ and apply Lemma 1.4.15. If $\operatorname{depth} R_\mathfrak{p} = 2$ then $M_\mathfrak{p}$ is MCM. If $\operatorname{depth} R_\mathfrak{p} \leq 1$ then $M_\mathfrak{p}$ is reflexive, hence MCM by *1)*. Then by Proposition 1.5.2 $M$ is MCM. $\square$

One of the most challenging problems in the theory of maximal Cohen-Macaulay modules is to determine whether the category $\operatorname{MCM}(R)$ has only finitely many isomorphism classes of indecomposable MCM modules or not.

**Definition 1.5.13.** Let $(R, \mathfrak{m}, k)$ be a Noetherian local ring, and let $\mathscr{M} = \{M_1, \ldots, M_r\}$ be a set of MCM $R$-modules. We say that $R$ has *finite Cohen-Macaulay type* given by $\mathscr{M}$ if every MCM $R$-module is isomorphic to a direct sum $M_1^{\oplus n_1} \oplus \cdots \oplus M_r^{\oplus n_r}$ for some natural numbers $n_1, \ldots, n_r$.

**Example 1.5.14.** If $R$ is regular, then the MCM modules are exactly the free modules by Proposition 1.5.10, so $R$ has finite CM type given by $\mathscr{M} = \{R\}$.

**Remark 1.5.15.** If $\mathrm{MCM}(R)$ is a Krull-Schmidt category, then to say that $R$ is of finite CM type is equivalent to the condition that $\mathrm{MCM}(R)$ has only a finite number of isomorphism classes of indecomposable objects.

Proposition 1.5.16 below furnish a way to construct examples of MCM modules: it is enough to take $n$-th syzygies, where $n$ is bigger or equal to the dimension of the ring. One of the modules arising in this way is particularly interesting. It is the top-dimensional syzygy module of the residue field $k$, that is $\mathrm{Syz}_R^{\dim R}(k)$. The main reason is that this module characterizes the regularity of the ring, as shown in Corollary 1.5.17. This is actually a consequence of the Auslander-Buchsbaum-Serre Theorem.

**Proposition 1.5.16.** *Let $(R, \mathfrak{m}, k)$ be a Cohen-Macaulay local ring of Krull dimension $d$. Then for any finitely generated $R$-module $M$ and for any integer $n \geq d$ the module $\mathrm{Syz}_R^n(M)$ is either $0$ or maximal Cohen-Macaulay.*

*Proof:* Assume that $\mathrm{Syz}_R^n(M) \neq 0$, then a minimal free resolution of $M$ has length at least $n$, so we can write the exact sequence

$$0 \to \mathrm{Syz}_R^n(M) \to F_{n-1} \xrightarrow{\varphi_{n-1}} \cdots \to F_1 \xrightarrow{\varphi_1} F_0, \qquad (1.2)$$

where the modules $F_i$ are finite free. From (1.2) we get the short exact sequence

$$0 \to \mathrm{Syz}_R^n(M) \to F_{n-1} \to \mathrm{Syz}_R^{n-1}(M) \to 0,$$

since $\mathrm{Syz}_R^{n-1}(M) = \mathrm{Im}\,\varphi_{n-1} = \ker \varphi_{n-2}$. The Depth Lemma yields the recursive formula

$$\mathrm{depth}_R(\mathrm{Syz}_R^n(M)) = \min\{d, \mathrm{depth}_R(\mathrm{Syz}_R^{n-1}(M)) + 1\}.$$

Thus if $n \geq d$ then $\mathrm{depth}_R(\mathrm{Syz}_R^n(M)) \geq d$, so $\mathrm{Syz}_R^n(M)$ is MCM. $\square$

**Corollary 1.5.17.** *Let $(R, \mathfrak{m}, k)$ be a CM ring of Krull dimension $d$. Then $R$ is regular if and only if $\mathrm{Syz}_R^d(k)$ is free.*

*Proof:* Theorem 1.4.3 and Theorem 1.4.6 combined give the following chain of equivalences:

$$\begin{aligned}
R \text{ is regular} \quad &\Leftrightarrow \quad \mathrm{proj.dim}_R(k) < \infty \\
&\Leftrightarrow \quad \mathrm{proj.dim}_R(k) = d - \mathrm{depth}_R(k) \leq d \\
&\Leftrightarrow \quad \mathrm{Syz}_R^d(k) \text{ is free.}
\end{aligned}$$

$\square$

### 1.5.1. Geometry of MCM modules

From a geometric point of view, maximal Cohen-Macaulay modules can be seen as a measure of the singularities of the ring. Over a regular ring MCM modules are exactly the free modules, so in a certain sense there exists no non-trivial MCM module. The next step in this direction is a ring of finite CM type, where only finitely many indecomposable MCM modules appear up to isomorphism. A beautiful result of Auslander says that these rings are necessarily isolated singularities, so they are not too much complicated in geometric terms.

Moreover, it turns out that the category of MCM modules over isolated singularites has interesting properties. In particular, since MCM modules are locally free over the punctured spectrum $U$ of an isolated singularity, they correspond to vector bundles over $U$. This is a powerful and important tool to translate algebraic information into geometric information and viceversa.

**Definition 1.5.18.** Let $(R, \mathfrak{m}, k)$ be a Noetherian local ring. $R$ is called an *isolated singularity* if for any $\mathfrak{p} \in \mathrm{Spec}(R)$, $\mathfrak{p} \neq \mathfrak{m}$, the ring $R_\mathfrak{p}$ is regular.

Note that we include also regular rings in this definition. Sometimes we will say that a ring has at most an isolated singularity to stress this fact. Rephrasing Serre's criterion for normality (Theorem 1.4.12) we obtain the following important characterization of isolated singularity in dimension two.

**Theorem 1.5.19** (Serre)**.** *A two-dimensional local domain is normal if and only if it is Cohen-Macaulay and isolated.*

**Definition 1.5.20.** Let $(R, \mathfrak{m}, k)$ be a Noetherian local ring. A finitely generated module $M$ is called *locally free on the punctured spectrum* if the $R_\mathfrak{p}$-module $M_\mathfrak{p}$ is free for any prime ideal $\mathfrak{p} \neq \mathfrak{m}$.

Lemma 1.5.21 below is a direct consequence of Proposition 1.5.2 and Proposition 1.5.10.

**Lemma 1.5.21.** *Let $(R, \mathfrak{m}, k)$ be a Cohen-Macaulay isolated singularity and let $M$ be a MCM $R$-module. Then $M$ is locally-free on the punctured spectrum.*

The following result was first proved by Auslander in [Aus86a] in the complete setting using the theory of almost split sequences. Huneke and Leuschke [HL02] gave a shorter proof and extended the result also to non-complete CM rings.

**Theorem 1.5.22** (Auslander)**.** *Let $(R, \mathfrak{m}, k)$ be a local Cohen-Macaulay ring of finite Cohen-Macaulay type. Then $R$ is an isolated singularity.*

**Lemma 1.5.23.** *Let $(R, \mathfrak{m}, k)$ be a normal two-dimensional local domain, let $X = \mathrm{Spec}\, R$, and let $U = X \setminus \{\mathfrak{m}\}$. For a torsion-free finitely generated $R$-module $M$ we have an isomorphism*

$$M^{**} \cong \Gamma(U, \widetilde{M}).$$

*Proof:* First of all, we assume that $M$ is MCM. From the exact sequence

$$0 \to H^0_{\mathfrak{m}}(M) \to M \to \Gamma(U, \widetilde{M}) \to H^1_{\mathfrak{m}}(M) \to 0,$$

we obtain that $M \cong \Gamma(U, \widetilde{M})$, since $H^0_{\mathfrak{m}}(M) = H^1_{\mathfrak{m}}(M) = 0$ for a MCM module $M$.

Now, let $M$ be torsion-free and consider the following short exact sequence

$$0 \to M \to M^{**} \to T \to 0.$$

Since $M$ and $M^{**}$ have the same Krull dimension, $T$ is zero-dimensional, hence it has finite length. For every prime ideal $\mathfrak{p} \neq \mathfrak{m}$ we have $T_{\mathfrak{p}} = 0$, so we obtain an isomorphism of sheaves $\widetilde{M}|_U \cong \widetilde{M^{**}}|_U$. Moreover $M^{**}$ is MCM, so we have

$$M^{**} \cong \Gamma(U, \widetilde{M^{**}}) \cong \Gamma(U, \widetilde{M}),$$

which concludes the proof. $\square$

**Remark 1.5.24.** We point out that in the proof of Lemma 1.5.23 we proved actually that the restriction of the sheaf $\widetilde{M^{**}}$ on the punctured spectrum $U$ coincide with the restriction to $U$ of the sheaf $\widetilde{M}$, that is

$$\widetilde{M^{**}}|_U = \widetilde{M}|_U.$$

**Corollary 1.5.25.** *Let $(R, \mathfrak{m}, k)$ be a normal two-dimensional local ring, and let $\mathrm{VB}(U)$ denote the category of vector bundles over $U = \mathrm{Spec}\, R \setminus \{\mathfrak{m}\}$. Then the functor $\mathrm{MCM}(R) \to \mathrm{VB}(U)$ mapping a MCM $R$-module $M$ to the locally free sheaf $\widetilde{M}|_U$ is an equivalence of categories. Moreover, for any two objects $M_1, M_2 \in \mathrm{MCM}(R)$ we have*

$$\widehat{M_1 \boxtimes_R M_2}\Big|_U \cong \widetilde{M_1}|_U \otimes \widetilde{M_2}|_U.$$

**Remark 1.5.26.** If the ring $R$ has a canonical module $K_R$, then many authors consider the canonical dual $(-)^{\vee} := \mathrm{Hom}_R(-, K_R)$, instead of the standard dual $(-)^* = \mathrm{Hom}_R(-, R)$. Lemma 1.5.23 is true also for the canonical dual, see for example [BD08, Proposition 3.10].

# 2. Signatures

## 2.1. F-signature

We want to give a short review of the *F-signature*, which is the main inspiration that lead us to the definition of the symmetric signature. This is by no mean intended to be a complete survey on F-signature, rather we will focus on those properties and results that we consider important to understand the connection between these two numerical invariants.

We recall some general facts about rings of prime characteristic. Let $R$ be a commutative ring containing a field of prime characteristic $p$. With the letter $e$ we always denote a natural number, and with $q = p^e$ a power of the characteristic. The *Frobenius homomorphism* is the ring homomorphism $F : R \to R$, $F(r) = r^p$, we often consider also its iterates $F^e : R \to R$, $F^e(r) = r^q$. For any finitely generated $R$ module $M$, we denote by $^eM$ the $R$-module $M$, whose multiplicative structure is pulled back via $F^e$. The scalar multiplication on $^eM$ is given by $r \cdot m := r^q m$, for any $r \in R$, $m \in {}^eM$. If $R$ is reduced, we can identify the Frobenius map $F^e : R \to {}^eR$ with the inclusion $R \hookrightarrow R^{1/q}$, where $R^{1/q}$ is the over-ring of $q$-th roots of elements of $R$.

We say that $R$ is *F-finite* if $^1R$ is a finitely generated $R$-module, and we say that $R$ is *F-split* if the Frobenius map $F : R \to {}^1R$ splits. If $(R, \mathfrak{m}, k)$ is a complete Noetherian local ring, then $R$ is F-finite if and only if $[k^{1/p} : k]$ is finite. We will always assume that our rings are F-finite, and that $k$ is perfect, that is $[k^{1/p} : k] = 1$. The latter is not a strong restriction, but it is done simply to avoid that numbers like $[k^{1/p} : k]$ appear in the definitions and in the results.

**Remark 2.1.1.** If $R$ is F-finite and Cohen-Macaulay, then $^eR$ is a maximal Cohen-Macaulay module for every $e \in \mathbb{N}$. Actually, if $x_1, \ldots, x_n$ is a maximal $R$-regular sequence, then it is also an $^eR$-regular sequence for every $e \in \mathbb{N}$.

Let $(R, \mathfrak{m}, k)$ be a Noetherian local ring of characteristic $p > 0$ and let $I = (f_1, \ldots, f_r)$ be an $\mathfrak{m}$-primary ideal. We denote by $I^{[q]}$ the ideal $(f_1^q, \ldots, f_r^q)$, which does not depend on the choice of the generators $f_1, \ldots, f_r$.

**Definition 2.1.2.** The numerical function $HK(I, -)$ given by

$$HK(I, q) := l_R\left(R/I^{[q]}\right),$$

where $q = p^e$, is called *Hilbert-Kunz function of R with respect to I*.

**Theorem 2.1.3** (Monsky). *Let $(R, \mathfrak{m}, k)$ be a Noetherian local ring of characteristic $p > 0$ and dimension $d$ and let $I$ be an $\mathfrak{m}$-primary ideal. Then there exists a positive real number $e_{HK}(I)$ such that*

$$HK(I, q) = e_{HK}(I)q^d + O(q^{d-1}).$$

The number $e_{HK}(I)$ is called *Hilbert-Kunz multiplicity* of $I$. If $I = \mathfrak{m}$, then we denote it also by $e_{HK}(R)$ and we speak of the Hilbert-Kunz multiplicity of the ring $R$.

The Hilbert-Kunz multiplicity is a characteristic $p$ analogue of the classical Hilbert-Samuel multiplicity of local rings (cf. [Ser00, Chapter V]). As the classical multiplicity, also the Hilbert-Kunz multiplicity encodes information about the singularites of the ring. Watanabe and Yoshida [WY00] proved that if $R$ is unmixed, then $R$ is regular if and only if $e_{HK}(R) = 1$. Moreover the Hilbert-Kunz multiplicity is also related with F-singularities arising from tight closure theory. For a survey on Hilbert-Kunz function and multiplicity, and their connection with F-singularities, the reader may consult the excellent book of Huneke [Hun96], or his notes [Hun13].

The F-signature is a more recent development, and was defined by Huneke and Leuschke [HL02]. They continued the work, started by Smith and Van Den Bergh [SVB96], of looking at the free $R$-summands of the module ${}^eR$ for $e \to +\infty$.

Let $(R, \mathfrak{m}, k)$ be a $d$-dimensional Noetherian reduced local ring of prime characteristic $p$, which is F-finite and such that $k$ is perfect. For every $q = p^e$, we denote by $a_q = \text{frk}_R({}^eR)$ the free rank of ${}^eR$. In other words, $a_q$ is choosen such that

$${}^eR \cong R^{a_q} \oplus M_q,$$

and the module $M_q$ contains no free $R$-direct summands.

**Definition 2.1.4.** The *F-signature* of $R$, denoted by $s(R)$ is the limit

$$s(R) := \lim_{e \to +\infty} \frac{a_q}{q^d}.$$

It is not easy to prove that the limit defining the F-signature exists. It was proved by Huneke and Leuschke [HL02] assuming that $R$ is Gorenstein, and then in some other special cases by other authors (see e.g. [Yao05] for rings of finite F-representation type, or [Sin05] for affine semigroup rings). Finally, Tucker [Tuc12] proved that the F-signature exists for every reduced F-finite Noetherian local ring.

**Remark 2.1.5.** If $R$ is a domain of dimension $d$, then $\operatorname{rank}_R({}^eR) = q^d$. So we can write the limit defining the F-signature as

$$s(R) = \lim_{e \to +\infty} \frac{\operatorname{frk}_R({}^eR)}{\operatorname{rank}_R({}^eR)}.$$

**Remark 2.1.6.** Almost at the same time of Huneke and Leuschke, Watanabe and Yoshida [WY04] introduced the notion of *minimal relative Hilbert-Kunz multiplicity*. This is defined as

$$me_{HK}(R) := \liminf_{e \to +\infty} \frac{l_R\left(R/\operatorname{ann}_R z^q\right)}{q^d},$$

where $z$ is a generator of the socle of the injective hull of the residue field $k$. If $R$ is Gorenstein this is in fact

$$me_{HK}(R) = e_{HK}(J) - e_{HK}(J : \mathfrak{m}), \tag{2.1}$$

for any parameter ideal $J$ of $R$. Yao [Yao06] proved that the minimal relative Hilbert-Kunz multiplicity coincide with the F-signature.

The F-signature is always a real number in the interval $[0, 1]$. The extreme value 1 is obtained if and only if the ring is regular. This is a consequence of the results of Watanabe and Yoshida [WY00] and Proposition 2.1.7 below.

**Proposition 2.1.7.** *Let $(R, \mathfrak{m}, k)$ be a reduced F-finite Cohen-Macaulay local ring such that $k$ is infinite, Then*

$$(e(R) - 1)(1 - s(R)) \geq e_{HK}(R) - 1,$$

*where $e_{HK}(R)$ is the Hilbert-Kunz multiplicity and $e(R)$ is the Hilbert-Samuel multiplicity.*

**Theorem 2.1.8** (Watanabe-Yoshida)**.** *Let $(R, \mathfrak{m}, k)$ be a reduced F-finite local ring such that $k$ is infinite and perfect. Then $s(R) = 1$ if and only if $R$ is regular.*

Also the value $s(R) = 0$ has a special meaning in terms of the singularities of the ring, it is equivalent to the ring being not strongly F-regular, an important notion in tight closure theory.

**Definition 2.1.9.** The ring $R$ is *strongly F-regular* if for every $c \in R$ not in any minimal prime of $R$ the inclusion $Rc^{1/q} \subseteq R^{1/q}$ splits for $q \gg 0$. $R$ is called *weakly F-regular* provided every ideal of $R$ is tightly closed.

The following implications hold

$$\text{regular} \implies \text{strongly F-regular} \implies \text{weakly F-regular}, \tag{2.2}$$

and the last arrow can be reversed in some special cases (e.g. for Gorenstein rings [Hun96] or for F-finite $\mathbb{N}$-graded rings [LS99]), but it is not known in general. Moreover, we mention the fact that weakly F-regular rings are always Cohen-Macaulay and normal.

**Theorem 2.1.10** (Aberbach-Leuschke, [AL03])**.** *Let $(R, \mathfrak{m}, k)$ be a reduced excellent F-finite local ring of prime characteristic such that $k$ is perfect. Then $s(R) > 0$ if and only if $R$ is strongly F-regular.*

**Remark 2.1.11.** If the ring $R$ has dimension 0 or 1, then it is normal if and only if it is regular. Therefore, from the implications in (2.2) it follows that strongly F-regular is equivalent to regular. So by Theorem 2.1.8 and Theorem 2.1.10 we have only two possibilities for the F-signature of $R$. We have $s(R) = 1$ if $R$ is regular, and $s(R) = 0$ otherwise.

Theorem 2.1.8 and Theorem 2.1.10 justify the statement that *F-signature measures singularities*. Roughly speaking, the closer the F-signature to 1 is, the nicer the singularity. Further examples are: the F-signature of quotient singularities computed by Watanabe and Yoshida [WY04], the F-signature of normal monomial rings computed by Singh [Sin05], and the F-signature of coordinate rings of affine toric varieties computed by Von Korff [VKo11].

**Theorem 2.1.12** (Watanabe-Yoshida)**.** *Let $R = k[\![x_1, \dots x_n]\!]^G$ be a quotient singularity over a field $k$ of prime characteristic $p$. Assume that the acting group $G \subseteq \mathrm{GL}(n, k)$ is a small finite group such that that $(p, |G|) = 1$. Then the F-signature of $R$ is*

$$s(R) = \frac{1}{|G|}.$$

For the definition of small group and related properties look at Definition 3.5.3 and the corresponding Section 3.5 on quotient singularities.

**Theorem 2.1.13** (Singh)**.** *Let $k$ be a perfect field of positive characteristic and let $R$ be a normal subring of the polynomial ring $k[x_1, \dots, x_d]$ which is generated, as $k$-algebra, by monomials in the variables $x_1, \dots, x_d$. Then the F-signature $s(R)$ exists and is a positive rational number which does not depend on the characteristic of $k$.*

**Example 2.1.14** ([Sin05])**.** Let $k$ be a perfect field of characteristic $p > 0$ and let $R$ be the $n$-th Veronese subring of $S = k[x_1, \dots, x_d]$, where $d \geq 2$ and $n \geq 1$. In other words, $R$ is the subring of $S$ which is generated by the monomials of degree $n$. For $d = 2$ we have the Veronese ring of further Example 4.5.1. The F-signature of $R$ is $s(R) = \frac{1}{n}$.

**Theorem 2.1.15** (Von Korff)**.** *Let $N$ be a lattice and let $M = N^*$ be its dual. Let $\sigma \subseteq N \otimes_{\mathbb{Z}} \mathbb{R}$ be a full-dimensional strongly convex rational polyhedral cone with primitive generators $v_1, \dots, v_r$. Let $k$ be a perfect field of prime characteristic, let $S = M \cap \sigma^{\vee}$, and let $R = k[S]$. Then $R$ is the coordinate ring of an affine toric variety, and the F-signature of $R$ is*

$$s(R) = \mathrm{Vol}(P_{\sigma}),$$

*where $P_{\sigma} = \{w \in M \otimes_{\mathbb{Z}} \mathbb{R} : 0 \leq w \cdot v_i < 1, \ \forall i = 1, \dots, r\}$ is a polytope.*

In all examples known so far, the F-signature is a rational number. Thus, it is natural to ask whether this is always true.

**Question 2.1.16.** Is the F-signature always a rational number?

The same question has been asked for the Hilbert-Kunz multiplicity as well. In this case for a long time a positive answer seemed to be the most reasonable, however Monsky [Mon08] proposed the following conjecture, which is supported by a huge amount of numerical evidence.

**Conjecture 2.1.17** (Monsky). Let $k$ be a finite field of characteristic 2 and $h = x^3 + y^3 + xyz \in k[x, y, z]$. Then the Hilbert-Kunz multiplicity of the hypersurface $R/f$, where $f = uv + h$ and $R = k[x, y, z, u, v]$ is

$$e_{HK}(R/f) = \frac{3}{4} + \frac{5}{14\sqrt{7}}.$$

**Remark 2.1.18.** The ring $R/f$ of Monsky's example is Gorenstein, so if Conjecture 2.1.17 is true, then also the F-signature $s(R/f)$ is an irrational number. To see this, consider a minimal reduction $J = (v, x, y, z)$ of the maximal ideal $\mathfrak{m}$ of $R/f$. This reduction has the property that $\mathfrak{m}/J$ is a vector space of dimension 1 generated by the residue class of the element $u$. This forces the equality $J : \mathfrak{m} = \mathfrak{m}$. Thus we have $e_{HK}(J) - e_{HK}(\mathfrak{m}) = s(R/f)$ by (2.1). Since $J$ is generated by a regular sequence and is a reduction of $\mathfrak{m}$, $e_{HK}(J) = e(J) = e(\mathfrak{m})$ the Hilbert-Samuel multiplicity, which is always an integer. Therefore we may write $e(\mathfrak{m}) - e_{HK}(R/f) = s(R/f)$, which explains the claim.

Question 2.1.16 and Conjecture 2.1.17 are still open as of the writing of this chapter. However Brenner [Bre13] exhibits an example of a ring with irrational Hilbert-Kunz multiplicity.

We conclude with the notion of generalized F-signature, introduced by Hashimoto and Nakajima [HN15]. Following the ideas and notations of Section 1.1 we will give a definition which is more general than the one contained in [HN15].

Let $(R, \mathfrak{m}, k)$ be a Noetherian reduced F-finite local ring of characteristic $p > 0$ such that $k$ is perfect. Let $\mathscr{C}$ be a full subcategory of $\mathrm{mod}(R)$ such that $\mathscr{C}$ has the KRS property and ${}^e R \in \mathscr{C}$ for every $e \in \mathbb{N}$. For example, if $R$ is complete we may choose $\mathscr{C} = \mathrm{mod}(R)$ by Theorem 1.1.10. Let $M$ be an indecomposable object in $\mathscr{C}$, for each $e \in \mathbb{N}$ and $q = p^e$ we can consider the multiplicity $a_q$ of $M$ into ${}^e R$. In other words, we can write ${}^e R = M^{a_q} \oplus N_q$, with the module $N_q$ containing no copies of $M$ as direct summands.

**Definition 2.1.19.** The *generalized F-signature of R with respect to M* is

$$s(R, M) = \lim_{e \to +\infty} \frac{a_q}{q^d},$$

provided the limit exists.

**Remark 2.1.20.** If the ring $R$ has finite F-representation type, then the generalized F-signature exists. This is a consequence of [SVB96, Proposition 3.3.1] and of [Yao05, Theorem 3.11]. Moreover, Seibert [Sei97] proved that in this situation the F-signature and the Hilbert-Kunz multiplicity of $R$ are rational numbers.

**Theorem 2.1.21** (Hashimoto-Nakajima)**.** *Let $R = k[\![x_1, \ldots x_n]\!]^G$ be a quotient singularity over an algebraically closed field $k$ of prime characteristic $p$. Assume that the acting group $G \subseteq \mathrm{GL}(n, k)$ is a small finite group such that that $(p, |G|) = 1$. Let $V_t$ be an irreducible $k$-representation of $G$ and let $M_t = \mathscr{A}(V_t) = (S \otimes_k V_t)^G$ be the corresponding $R$-module via Auslander functor. Then the generalized F-signature of $R$ with respect to $M_t$ is*

$$s(R, M_t) = \frac{\mathrm{rank}_R M_t}{|G|}.$$

The definition and the properties of quotient singularities and Auslander functor $\mathscr{A}$ will be studied in details in Chapter 3 (e.g. see Theorem 3.5.10).

**Remark 2.1.22.** In the situation of Theorem 2.1.21, Smith and Van Den Bergh [SVB96, Proposition 3.2.1] proved that $R$ has finite F-representation type given by $\{M_0, \ldots, M_{r-1}\}$, where $M_t = \mathscr{A}(V_t) = (S \otimes_k V_t)^G$ and $\{V_0, \ldots, V_{r-1}\}$ is a complete set of irreducible $k$-representations of $G$. In other words, the $R$-modules $M_t$ are the only indecomposable modules appearing in the decomposition of ${}^e R$ into indecomposable modules.

### 2.1.1. F-signature in characteristic zero

The F-signature is an intrinsic concept of rings of prime characteristic, but it would be profitable to have an analogous invariant which is independent of the characteristic and could be used to study singularities also in characteristic zero. The first approach that one may think of is perhaps the following.

Let $R$ be a reduced $\mathbb{Z}$-algebra such that $\mathrm{Spec} R \to \mathrm{Spec} \mathbb{Z}$ is dominant. For every prime number $p$ we consider the reduction mod $p$ $R_p := R \otimes_{\mathbb{Z}} (\mathbb{Z}/p\mathbb{Z})$ and compute the F-signature $s(R_p)$. One may ask whether the limit

$$\lim_{p \to +\infty} s(R_p)$$

exists, and use this limit to define a characteristic 0 version of the F-signature.

Despite the fact that it seems quite natural, the previous definition presents a fundamental problem. If $p_1$ and $p_2$ are two distinct prime numbers, the two Frobenius homomorphisms $F_{p_1} : R_{p_1} \to R_{p_1}$ and $F_{p_2} : R_{p_2} \to R_{p_2}$ are quite different in general, and difficult to compare. For example we have that $\mathrm{rank}_R({}^e R_{p_1}) = p_1^e$, and $\mathrm{rank}_R({}^e R_{p_2}) = p_2^e$. So it is difficult to find an appropriate meaning to the previous limit and compute it, and even to determine whether it exists or not.

On the other hand, results like those of Watanabe and Yoshida (Theorem 2.1.12), Singh (Theorem 2.1.13) and Von Korff (Theorem 2.1.15) show that often the F-signature does not depend on the characteristic, so a meaningfull extension to the characteristic zero setting could be done, at least in these cases.

In order to give a characteristic free definition, we choose to replace the role played by the Frobenius homomorphism with another construction which is characteristic independent, but recovers some of the properties of the Frobenius. In particular we look at two of them. The first one is the famous result of Kunz [Kun69] that the Frobenius homomorphism characterizes regularity.

**Theorem 2.1.23** (Kunz)**.** *Let $(R, \mathfrak{m}, k)$ be a reduced local ring of prime characteristic $p$. Then the following facts are equivalent.*

*1) $R$ is regular.*
*2) ${}^e R$ is a free $R$-module for all $e \in \mathbb{N}$.*
*3) ${}^1 R$ is a free $R$-module.*

The second property is the fact that ${}^e R$ is a MCM $R$-module, if the ring is Cohen-Macaulay (Remark 2.1.1).

Our idea is to use the top dimensional syzygy module of the residue field, which is a MCM module (Proposition 1.5.16) and characterizes regularity (Corollary 1.5.17), and apply reflexive symmetric powers (Definition 1.3.14) to obtain an asymptotic behaviour. This leads us to the definition of *symmetric signature*, which will be explained in the next section.

## 2.2. Symmetric signature

Keeping in mind the observations of the previous Section 2.1.1, we are going to define the main object of interest of this thesis: *the symmetric signature*. We will give its definition and explain some basic properties. Further computations and examples will be given in Chapter 4 and Chapter 5.

Let $(R, \mathfrak{m}, k)$ be a local Noetherian domain of dimension $d$. We consider a minimal free resolution of $k$

$$0 \longrightarrow \mathrm{Syz}_R^d(k) \longrightarrow R^{\beta_{d-1}} \longrightarrow \cdots \longrightarrow R^{\beta_1} \longrightarrow R \longrightarrow k \longrightarrow 0,$$

and the top-dimensional syzygy module $\mathrm{Syz}_R^d(k)$. For every natural number $q$, we consider the reflexive symmetric powers of the last module, for ease of notation we define

$$\mathscr{S}^q := \left( \mathrm{Sym}_R^q \left( \mathrm{Syz}_R^d(k) \right) \right)^{**}.$$

The module $\mathscr{S}^q$ will play almost the same role played by the module ${}^e R$ in the definition of F-signature. Here the advantage is that the construction used to define $\mathscr{S}^q$ does not depend explicitly from the characteristic of the ring. At the same time, $\mathscr{S}^q$ encodes some homological properties of the ring $R$ through the top-dimensional syzygy module $\mathrm{Syz}_R^d(k)$.

**Definition 2.2.1.** The number

$$s_\sigma(R) := \lim_{N \to +\infty} \frac{\sum_{q=0}^N \mathrm{frk}_R \mathscr{S}^q}{\sum_{q=0}^N \mathrm{rank}_R \mathscr{S}^q} \tag{2.3}$$

is called *symmetric signature* of $R$, provided the limit exists.

We point out that for $q = 0$, $\mathscr{S}^q = k^{**} \cong R$. So one has $\mathrm{frk}_R \mathscr{S}^0 = \mathrm{rank}_R \mathscr{S}^0 = 1$, and the sums in the previous limit start both from 1.

We don't know if the limit defining the symmetric signature always exists. To obtain an invariant which is for sure well defined and exists, one may replace the limit in (2.3) with a limit inferior.

Moreover one may ask why should we consider the limit of (2.3), insted of the simpler limit

$$\lim_{q \to +\infty} \frac{\mathrm{frk}_R \mathscr{S}^q}{\mathrm{rank}_R \mathscr{S}^q}. \tag{2.4}$$

The main reason is that the limit (2.4) does not exist even in simple cases, as it will be shown in Example 4.4.9. However, when the simpler limit (2.4) exists, then also the symmetric signature exists and they coincide, as shown by the following elementary statement.

**Lemma 2.2.2.** *Let $(a_n)_{n \in \mathbb{N}}$, $(b_n)_{n \in \mathbb{N}}$ be sequences of real numbers such that $b_n > 0$ for all $n \in \mathbb{N}$. Assume that the infinite series $\sum_{n \in \mathbb{N}} b_n$ diverges. Then, if the limit $\lim_{n \to +\infty} \frac{a_n}{b_n}$ exists, then also the limit*

$$\lim_{n \to +\infty} \frac{\sum_{k=0}^n a_k}{\sum_{k=0}^n b_k}$$

*exists and the two limits coincide.*

*Proof:* Let $a := \lim_{n \to +\infty} \frac{a_n}{b_n}$. Then $c_n := \frac{a_n}{b_n} - a$ converges to 0 and we have $a_n = (a + c_n) b_n$ for all $n$. Let $\varepsilon > 0$, then there exists a $n_0 \in \mathbb{N}$ such that for all natural numbers $n \geq n_0$ we have $|c_n| \leq \varepsilon$. Then

$$\left| \frac{\sum_{k=0}^n a_k}{\sum_{k=0}^n b_k} - a \right| = \left| \frac{\sum_{k=0}^n c_k b_k}{\sum_{k=0}^n b_k} \right|$$

$$\leq \left| \frac{\sum_{k=0}^{n_0} c_k b_k}{\sum_{k=0}^n b_k} + \frac{|c_n| \sum_{k=n_0+1}^n b_k}{\sum_{k=0}^n b_k} \right|$$

$$\leq \left| \frac{\sum_{k=0}^{n_0} c_k b_k}{\sum_{k=0}^n b_k} + \frac{\varepsilon \sum_{k=n_0+1}^n b_k}{\sum_{k=0}^n b_k} \right|$$

$$\leq \varepsilon + \left| \frac{\sum_{k=0}^{n_0} c_k b_k}{\sum_{k=0}^n b_k} \right|.$$

The last term tends to zero because the series $\sum_{n \in \mathbb{N}} b_n$ diverges, so the claim follows. $\square$

**Example 2.2.3.** The assumption that $\sum_{n \in \mathbb{N}} b_n$ diverges is necessary. Consider

$$a_n = \frac{n}{2^n} \text{ and } b_n = \frac{n+1}{2^n}.$$

Then $\lim_{n \to +\infty} \frac{a_n}{b_n} = 1$, but

$$\lim_{n \to +\infty} \frac{\sum_{k=0}^n a_k}{\sum_{k=0}^n b_k} = \frac{1}{2}.$$

**Remark 2.2.4.** The module $\mathscr{S}_q$ is reflexive for every $q \geq 0$. If $R$ is CM of dimension $\leq 2$, then $\mathscr{S}_q$ is also MCM.

**Example 2.2.5.** If $R$ is a regular ring of dimension $d$, then $\mathrm{Syz}_R^d(k)$ is a free module by Corollary 1.5.17. Then $\mathscr{S}^q$ is also free, therefore $\mathrm{frk}_R \mathscr{S}^q = \mathrm{rank}_R \mathscr{S}^q$ for every $q$. It follows that the symmetric signature is $s_\sigma(R) = 1$.

**Example 2.2.6.** If $R$ has dimension 0, then is a field since it is a domain. In particular $R$ is regular, so by the previous Example 2.2.5 we have $s_\sigma(R) = 1$.

**Remark 2.2.7.** If $(R, \mathfrak{m}, k)$ has dimension 1, then the first syzygy of the residue field $\mathrm{Syz}_R^1(k)$ is just the maximal ideal $\mathfrak{m}$, as explained by the following short exact sequence

$$0 \to \mathfrak{m} \to R \to R/\mathfrak{m} \to 0.$$

Since $\mathfrak{m}$ is an ideal it has rank 1, therefore also all its symmetric powers have rank 1, and the same holds for the reflexive hull. In other words, we have $\mathrm{rank}_R \left( \mathrm{Sym}_R^q(\mathrm{Syz}_R^1(k)) \right)^{**} = 1$ for all $q \in \mathbb{N}$. Therefore, for each $q$ we have only two possibilities: either $\mathscr{S}^q \cong R$ and $\mathrm{frk}_R \mathscr{S}^q = \mathrm{rank}_R \mathscr{S}^q = 1$ or $\mathscr{S}^q$ is not free and $\mathrm{frk}_R \mathscr{S}^q = 0$.

Now assume in addition that the category $\mathrm{Ref}(R)$ has the KRS property. We fix an indecomposable object $M$ in $\mathrm{Ref}(R)$, then for every $q \in \mathbb{N}$ we have a unique decomposition

$$\mathscr{S}^q = M^{a_q} \oplus N_q,$$

where the module $N_q$ contains no copy of $M$ as direct summand. The natural number $a_q$ is called *the multiplicity of $M$ in $\mathscr{S}^q$*.

**Definition 2.2.8.** The number

$$s_\sigma(R, M) := \lim_{N \to +\infty} \frac{\sum_{q=0}^N a_q}{\sum_{q=0}^N \mathrm{rank}_R \mathscr{S}^q}$$

is called *generalized symmetric signature of $R$ with respect to $M$*, provided the limit exists.

We can define the symmetric signature also for graded rings.

Let $R$ be a standard graded Noetherian $k$-domain of dimension $d$. We work in the category $\mathrm{mod}_{\mathbb{Z}}(R)$ of finitely generated graded $R$-modules. In this category we consider the graded module $\mathrm{Syz}_R^d(k)$ coming from a minimal graded free resolution of $k$ as $R$-module, and its reflexive symmetric powers $\mathscr{S}^q := \left( \mathrm{Sym}_R^q \left( \mathrm{Syz}_R^d(k) \right) \right)^{**}$, which are again naturally graded. The symmetric signature of $R$ is defined by the limit

$$s_\sigma(R) := \lim_{N \to +\infty} \frac{\sum_{q=0}^N \mathrm{frk}_R^{\mathrm{gr}} \mathscr{S}^q}{\sum_{q=0}^N \mathrm{rank}_R \mathscr{S}^q},$$

provided it exists.

**Remark 2.2.9.** Notice that in the definition of $s_\sigma(R)$ for graded rings we consider the graded free rank of Definition 1.1.19 instead of the free rank as in the local setting. Therefore some differences may occur. For example, one should remind that all direct summands of the form $R(a)$ for some integer $a$ contribute to the graded free rank part of the module.

### 2.2.1. Differential symmetric signature

We present an alternative version of the symmetric signature, which is based on the reflexive hull of the differential module $\Omega_{R/k}$ instead of the top-dimensional syzygy module $\mathrm{Syz}_R^d(k)$. We recall briefly the definition of the module $\Omega_{R/k}$ and some basic properties. For a more detailed description the reader may consult [Eis94, Chapter 16].

Let $R$ be a commutative ring, and let $M$ be an $R$-module. A *derivation* d from $R$ to $M$ is a map of abelian groups $\mathrm{d} : R \to M$ which satisfies the Leibniz rule $\mathrm{d}(fg) = f\mathrm{d}(g) + g\mathrm{d}(f)$ for all $f, g \in R$. If $R$ is an algebra over a field $k$, and the derivation d is a $k$-linear map, then we say that d is a $k$-linear derivation. A classical example is given by $R = k[x]$ a polynomial ring in one variable, and $\mathrm{d} : R \to R$ the formal derivative of a polynomial.

**Definition 2.2.10.** Let $R$ be an algebra over a field $k$. The *cotangent module* or the *module of Kähler differentials* of $R$ over $k$ is denoted by $\Omega_{R/k}$, and is the $R$-module generated by all elements of the form $\mathrm{d}f$ with $f \in R$ subject to the relations

$$\mathrm{d}(fg) = f\mathrm{d}g + g\mathrm{d}f,$$
$$\mathrm{d}(af + bg) = a\mathrm{d}f + b\mathrm{d}g,$$

for all $f, g \in R$ and $a, b \in k$. The map $\mathrm{d} : R \to \Omega_{R/k}$ defined by $f \mapsto \mathrm{d}f$ is called the *universal $k$-linear derivation*.

The cotangent module can be also defined by the following universal property. For any $R$-module $M$ and $k$-linear derivation $\mathrm{e} : R \to M$, there exists a unique $k$-linear homomorphism $\varphi : \Omega_{R/k} \to M$ such that $\mathrm{e} = \varphi\mathrm{d}$. In other words, we have the following commutative diagram.

$$
\begin{array}{ccc}
R & \xrightarrow{\ \mathrm{d}\ } & \Omega_{R/k} \\
& \searrow{\scriptstyle \mathrm{e}} & \downarrow{\scriptstyle \exists!\, \varphi} \\
& & M.
\end{array}
$$

**Example 2.2.11.** Let $R = k[x_1, \ldots, x_n]$ be the polynomial ring in $n$ variables. Then $\Omega_{R/k}$ is the free module generated by the elements $\mathrm{d}x_1, \ldots, \mathrm{d}x_n$.

**Proposition 2.2.12** (Conormal Sequence). *Let $A$ be a $k$-algebra, let $I$ be an ideal of $A$, and let $R = A/I$. Then there is a natural exact sequence of $R$-modules*

$$I/I^2 \to \Omega_{A/k} \otimes_A R \to \Omega_{R/k} \to 0.$$

As a corollary to the last proposition, one obtains the remarkable fact that if $R$ is a finitely generated $k$-algebra, then $\Omega_{R/k}$ is a finitely generated $R$-module.

In the local setting, the cotangent module characterizes regularity as well as the top-dimensional syzygy module of the residue field does.

**Theorem 2.2.13.** *Let $(R, \mathfrak{m}, k)$ be a Noetherian local ring. Assume that $k$ is perfect, and $R$ is a localization of a finitely generated $k$-algebra. Then $\Omega_{R/k}$ is a free $R$-module of rank $\dim R$ if and only if $R$ is a regular local ring.*

*Proof:* See [Har77, Theorem 8.8]. □

**Remark 2.2.14.** Differently from $\mathrm{Syz}_R^{\dim R}(k)$, the cotangent module $\Omega_{R/k}$ is not necessarily maximal Cohen-Macaulay, and not even torsion-free. In fact, Herzog [Her78a] proved that if $(R, \mathfrak{m}, k)$ is a Noetherian local $k$-algebra of dimension 1 with $\mu(\mathfrak{m}) = 3$, then the torsion submodule of $\Omega_{R/k}$ is not zero.

Since in general $\Omega_{R/k}$ is not reflexive, it is convenient to work with its reflexive hull $\Omega_{R/k}^{**}$. The module $\Omega_{R/k}^{**}$ is called *module of Zariski (or regular) differentials of $R$ over $k$*, and has been studied intensively by Auslander [Aus86b], and Martsinkovsky [Mar87], [Mar90], [Mar92].

**Definition 2.2.15.** Let $(R, \mathfrak{m}, k)$ be a local Noetherian $k$-domain of dimension $d$, and for every natural number $q$ consider the $R$-module

$$\mathscr{C}^q := \left( \mathrm{Sym}_R^q \left( \Omega_{R/k}^{**} \right) \right)^{**}.$$

The number

$$s_{d\sigma}(R) := \lim_{N \to +\infty} \frac{\sum_{q=0}^N \mathrm{frk}_R \mathscr{C}^q}{\sum_{q=0}^N \mathrm{rank}_R \mathscr{C}^q}$$

is called *differential symmetric signature* of $R$, provided the limit exists.

**Remark 2.2.16.** If $R$ is two-dimensional and normal, then by Lemma 1.5.23 and Remark 1.5.24 the first double dual in the definition of $\mathscr{C}^q$ is unnecessary, that is

$$\mathscr{C}^q = \left( \mathrm{Sym}_R^q \left( \Omega_{R/k} \right) \right)^{**}.$$

If $R$ is regular, $\Omega_{R/k}$ is free by Theorem 2.2.13, hence also the module $\Omega_{R/k}^{**}$ is free. It follows that for regular rings the differential symmetric signature is 1, as for the symmetric signature.

In the rest of the thesis, and in particular in Chapter 4 and Chapter 5 we will concentrate on the symmetric signature of Definition 2.2.1. However, we will mention each time which results extend also to the differential symmetric signature, and which differences may occur.

# 3. Quotient Singularites

Let $G$ be a finite group and let $A$ be a commutative ring. For an $A$-module $V$ we denote by $\mathrm{GL}(V)$ the group of $A$-linear automorphisms of $V$. If $V$ is free of rank $n$, then this group is isomorphic to to the group of all non-singular $n \times n$ matrices over $A$ and we denote it also by $\mathrm{GL}(n, A)$.

**Definition 3.0.17.** An *A-representation* of a finite group $G$ is a couple $(V, \rho)$, where $V$ is a finitely generated $A$-module and $\rho : G \to \mathrm{GL}(V)$ is a homomorphism of groups. If the module $V$ is free, then the rank of $V$ is called the *dimension* or the *degree* of the representation.

Sometimes by abuse of notation we will call the module $V$ the representation, but it should more properly be called the representation module or representation space. If $V$ is a free $A$-module we may write each element of $\mathrm{GL}(V)$ as a matrix, so we obtain for each $g \in G$ a matrix $\rho(g)$. These matrices act on the elements of $V$ by the usual matrix-vector multiplication, sometimes we will write $gv$ or $g(v)$ instead of $\rho(g)(v)$ if the map $\rho$ is clear.

**Example 3.0.18.** We take $V = A$ and $\rho(g) = \mathrm{id}_A$ for all $g \in G$. This representation is called *trivial representation* and it is sometimes denoted simply by $A$.

An important way to think about representations is in terms of group actions. In fact an $A$-representation $(V, \rho)$ of $G$ defines an action of $G$ on $V$ given by

$$G \times V \to V, \ (g, v) \mapsto \rho(g)(v).$$

This action is linear, in the sense that $g(v_1 + v_2) = gv_1 + gv_2$ and $g(\alpha v) = \alpha g(v)$ for every $v, v_1, v_2 \in V$, $\alpha \in A$, and $g \in G$. Viceversa, given a linear action $\psi : G \times V \to V$, we obtain a representation

$$G \to \mathrm{GL}(V), \ g \mapsto \psi(g, -),$$

which is called the *fundamental representation* of the group action.

### 3. Quotient Singularites

If the homomorphism $\rho$ is injective, then the representation is said to be *faithful*. In this case one may identify the abstract group $G$ with its image $\rho(G)$ inside $\mathrm{GL}(V)$. With this identification in mind we will write sometimes $G \subset \mathrm{GL}(V)$ meaning that the group is acting on $V$ through a faithful representation.

An important situation is when the representation module $S$ has also the structure of an $A$-algebra, and the $\rho(g) : S \to S$ are $A$-algebra-homomorphisms. Then, the representation $\rho : G \to \mathrm{GL}(S)$ induces an $A$-linear action on $S$. In particular if $A = \mathbb{Z}$, then $\mathrm{GL}(S) = \mathrm{Aut}(S)$ the automorphism group of $S$ and one recovers the classical theory of group actions on rings. In this setting we define the invariant subring of $S$ as

$$S^G := \{x \in S : g(x) = x \ \forall g \in G\} \subseteq S.$$

Moreover if the order of the group $|G|$ is invertible in $S$ we can define the *Reynolds operator*

$$\rho : S \to S^G, \ \rho(s) = \frac{1}{|G|} \sum_{g \in G} g(s),$$

which is a left splitting for the inclusion $i : S^G \hookrightarrow S$, that is $\rho \circ i = \mathrm{id}_{S^G}$. It follows that $S^G$ is a $S^G$-direct summand of $S$.

If the ring $S$ is a regular ring, then the invariant ring $S^G$ is called *quotient singularity*. The main goal of this chapter is to recall some fundamental facts and properties of quotient singularities. In particular we will review and prove the so-called *Auslander correspondence* (Theorem 3.5.10). This states that, under certain assumptions, there is a one-one correspondence between representations of $G$ and MCM $S^G$-modules. It follows that the study of the invariant ring $S^G$ and its modules is deeply connected with the study of the representations of the group $G$.

We will procede as follows. In Sections 3.1 and 3.2 we will review some basic facts concerning finite group representations over a field and their characters. In Section 3.3 we will introduce the skew group ring, which is a twisted version of the classical group ring, and is necessary to construct the Auslander correspondence. The latter will be the main object of study in Section 3.4. There, we will work over a Noetherian normal local domain equipped with a group action and we will explicitly construct a functor between the category of group representations and the category of reflexive modules over the invariant ring, using modules over the skew group ring as an intermediate step. Finally in Section 3.5 we will apply this correspondence to the quotient singularites and we will also present other important results on quotient singularites.

We claim no originality for the results of this chapter. However the theorems of Section 3.4 are stated and proved for a normal Noetherian local ring, instead of a power series ring as it is classically done in the literature ([Yos90, Chapter 10] and [LW12, Chapter 5]).

## 3.1. Representation theory of finite groups

We consider a finite group $G$ and we work over an arbitrary field $k$. Many definitions and results from this section can be extended to representations over any (not necessarily commutative) ring $A$. However for our purpose representations over a field are enough, so we restrict already to this situation to avoid many complications.

We will skip proofs of the main results of this section and the following Section 3.2. The interested reader may consult the books of Etingof et al. [EGHLSVY11], Feit [Fei82], Fulton and Harris [FH91] (for representations over $\mathbb{C}$), and Serre [Ser77], or the notes of Webb [Web14].

We begin with some definitions.

**Definition 3.1.1.** Let $(V, \rho)$ be a $k$-representation. A subspace $W$ of $V$ is said to be *stable under $G$* or a *$G$-invariant subspace* of $V$ if $gw \in W$ for all $w \in W$ and $g \in G$. Such an invariant submodule $W$ gives rise to a representation $\rho_W : G \to \mathrm{GL}(W)$, which is called a *subrepresentation* of $V$. A representation whose unique $G$-invariant subspaces are $0$ and $V$ itself is called *irreducible*.

**Remark 3.1.2.** One should not confuse a $G$-invariant subspace $W$ of $V$ with the subspace of fixed points $V^G = \{v \in V : gv = v \ \forall g \in G\}$.

**Definition 3.1.3.** A *$G$-linear map* between two representations $(V, \rho)$ and $(W, \rho')$ is a $k$-linear homomorphism $V \to W$ such that the diagram

$$
\begin{array}{ccc}
V & \longrightarrow & W \\
\rho(g) \downarrow & & \rho'(g) \downarrow \\
V & \longrightarrow & W
\end{array}
$$

commutes for every $g \in G$. We denote the space of $G$-linear maps between $V$ and $W$ by $\mathrm{Hom}_G(V, W)$. Two representations $V$ and $W$ are said to be *equivalent* or *isomorphic* if there exists $\varphi \in \mathrm{Hom}_G(V, W)$ bijective. In particular $\varphi$ is an isomorphism of vector spaces, so two equivalent representations have the same dimension.

Given two representations $V$ and $W$ we can construct the following representations:

- the *direct sum $V \oplus W$*, with action $g(v + w) = gv + gw$ for $g \in G$, $v, w \in V$;

- the *tensor product $V \otimes_k W$*, with action $g(v \otimes w) = gv \otimes gw$ for $g \in G$, $v, w \in V$;

- the *homomorphism representation* $\mathrm{Hom}_k(V, W)$, with action $g f(v) = f(g^{-1}v)$ for $g \in G$, $v \in W$ and $f : V \to W$ $k$-linear. In particular the representation $V^* = \mathrm{Hom}_k(V, k)$ is called *dual representation* of $V$.

With the same rule, if $V$ is a representation then the *$n$-th tensor product $V^{\otimes n}$* is again a representation and the *exterior powers $\bigwedge^n V$* and the *symmetric powers $\mathrm{Sym}^n(V)$* are subrepresentations of it.

It is not hard to see that the $G$-linear maps between two representations $V$ and $W$ coincide with the fixed subspace of $\mathrm{Hom}_k(V, W)$, that is

$$\mathrm{Hom}_G(V, W) = \mathrm{Hom}_k(V, W)^G.$$

**Example 3.1.4.** Let $X$ be a finite set with an action of $G$, then we associate a representation to this action. Let $V$ be a vector space with basis $\{e_x : x \in X\}$. $G$ acts naturally on this basis by $g e_x = e_{gx}$, and we can extend the action to $V$ by linearity. The representation associated with this action is called *permutation representation*. When $X = G$ the group itself, it is also called *regular representation*. The regular representation plays a key role in the study of irreducible representations.

Another important way to look at representations is in terms of the group algebra $k[G]$. In fact every statement about the representations of $G$ can be rephrased in terms of the group algebra.

**Definition 3.1.5.** The *group algebra* or *group ring* of $G$ over $k$ is denoted by $k[G]$ (or $kG$) and is the $k$-vector space with basis $\{e_g : g \in G\}$ and with multiplication given on the basis elements by group multiplication $e_g e_h = e_{gh}$ and extended linearly to arbitrary elements.

More explicitely given $\sum_{g \in G} a_g e_g$ and $\sum_{g \in G} b_g e_g$, where $a_g$ and $b_g$ are elements of $k$ we have

$$\left( \sum_{g \in G} a_g e_g \right) \left( \sum_{g \in G} b_g e_g \right) = \sum_{g \in G} \left( \sum_{fh=g} a_f b_h \right) e_g.$$

The group algebra is commutative if and only if the group $G$ is abelian. So when we speak of $k[G]$ modules we will always mean left $k[G]$-module.

**Proposition 3.1.6.** *We have a one-one correspondence between the following objects:*

1. *$k$-representations of $G$;*
2. *finitely generated $k[G]$-modules.*

*Moreover the dimension of a representation is equal to the $k$-dimension of the corresponding $k[G]$-module.*

*Proof:* Given a $k$-representation $(V, \rho)$ of $G$ we define a $k[G]$-module action on $V$ by $e_g \cdot v = \rho(g)(v)$ for every $g \in G$, $v \in V$ and we extend by linearity.

Conversely given a $k[G]$-module $V$, then for every $g \in G$ the multiplication by the basis element $e_g$ is a $k$-linear map $\rho_g : V \to V$. Since $\rho_g \rho_h = \rho_{gh}$ we obtain a representation $\rho : G \to \mathrm{GL}(V)$. The statement on the dimension is clear. □

**Remark 3.1.7.** The $k[G]$-module $k[G]$ corresponds to the regular representation from Example 3.1.4.

We can translate the definitions we have given so far in the language of module theory. So a *G*-linear map between representations is just a homomorphism of $k[G]$-modules, that is

$$\mathrm{Hom}_{k[G]}(V, W) = \mathrm{Hom}_G(V, W).$$

Similarly, a subrepresentation is nothing but a $k[G]$-submodule and irreducible representations correspond to simple $k[G]$-modules. In fact we have an equivalence of categories between the category $\mathrm{mod}(k[G])$ of $k[G]$-modules and the category of representations, whose object are representations and whose morphisms are the *G*-linear maps. Keeping in mind this identification we will use indifferently the language of representations or the language of $k[G]$-modules.

When the characteristic of the field *k* does not divide the order of the group we say that *G* is *non-modular* over *k*. In particular this happens if *k* has characteristic 0. The structure of the representations of *G* is quite nice in the non-modular case, in fact every representation is isomorphic to a direct sum of irreducible ones. In other words, the category $\mathrm{mod}(k[G])$ has the KRS property. This fact is a consequence of the Theorem of Maschke, which ensures the existence of a decomposition, and of Schur's Lemma, which provides the unicity of such a decomposition.

**Theorem 3.1.8** (Maschke)**.** *Let k be a field and let G be a finite group such that $|G|$ is invertible in k. Let V be a representation of G and let W be an invariant subspace of V. Then there exists an invariant subspace $W_1$ of G such that $V = W \oplus W_1$ as representations.*

Using the language of $k[G]$-modules, Maschke's Theorem affirms that in the non-modular case the group algebra $k[G]$ is semisimple, which means that it can be written as direct sum of simple modules.

**Lemma 3.1.9** (Schur)**.** *Let G be a finite group, let k be field, let V and W be two irreducible k-representations of G and let $\varphi : V \to W$ be a G-linear homomorphism. Then the following facts hold.*

1) *Either $\varphi$ is an isomorphism or $\varphi = 0$.*
2) *If $V = W$ and k is algebraically closed, then $\varphi = \lambda \cdot \mathrm{id}_V$ for some $\lambda \in k$, i.e. $\mathrm{Hom}_G(V, V) \cong k$.*

In many results that will follow we will make the assumption that *k* is algebraically closed to make sure that the second condition of Schur's Lemma holds. Sometimes this is a significant point, but most of the times this is just a convenience to simplify the statements so that numbers such as $\dim_k \mathrm{Hom}_G(V, V)$ do not occur.

**Corollary 3.1.10.** *Let k be a field, let G be a finite group such that $|G|$ is invertible in k and let V be a representation of G. Then there is a unique decomposition*

$$V = V_1^{a_1} \oplus \cdots \oplus V_m^{a_m},$$

*where the $V_i$'s are distinct irreducible representations and the $a_i$'s are natural numbers. In other words the category* $\mathrm{mod}(k[G])$ *has the KRS property and its indecomposable objects are precisely the irreducible representations.*

**Remark 3.1.11.** Maschke's Theorem and Corollary 3.1.10 are not true in the modular case (if char$k$ divides $|G|$). In this case, it turns out that indecomposable objects and irreducible representations do not coincide. In other words, $k[G]$-modules no longer need to be semisimple, that is direct sums of simple modules.

Another consequence of Schur's Lemma is that irreducible representations of an abelian group are one-dimensional.

**Corollary 3.1.12.** *Let $k$ be an algebraically closed field and $G$ an abelian finite group. Then every irreducible $k$-representation is one-dimensional.*

*Proof:* Let $\rho : G \to \mathrm{GL}(V)$ be an irreducible $k$-representation. Let $g \in G$, then the linear map $\rho(g) : V \to V$ is $G$-linear. In fact given $h \in G$ we have $\rho(g)\rho(h) = \rho(gh) = \rho(hg) = \rho(h)\rho(g)$, since $G$ is abelian. Then by Schur's Lemma $\rho(g) = \lambda_g \mathrm{id}_V$ for some $\lambda_g \in k$. Thus every linear subspace of $V$ is $G$-invariant, which forces $V$ to be one-dimensional. $\square$

The important role of the regular representation is expressed in the following Theorem 3.1.13. It says that the regular representation contains at least one copy of every irreducible representation. As a consequence, there exists only a finite number of irreducible representations.

**Theorem 3.1.13.** *Let $k$ be a field and $G$ a finite group such that $|G|$ is invertible in $k$. Let $k[G]$ be the regular representation of $G$ then*

$$k[G] \cong V_1^{n_1} \oplus \cdots \oplus V_r^{n_r},$$

*where the $V_i$ are pairwise non-isomorphic simple $k[G]$-modules. Moreover we have that $V_1, \ldots, V_r$ are a complete set of representatives of the isomorphism classes of simple $k[G]$-modules.*

*If $k$ is algebraically closed, then*

$$n_i = \dim_k V_i, \quad |G| = \dim_k k[G] = n_1^2 + \cdots + n_r^2,$$

*and $r$ is equal to the number of conjugacy classes of $G$.*

The previous theorem provides a good numerical criterion to determine when we have constructed all the irreducible $k$-representations of a group $G$. However one would still like to solve the following problem: given a representation $V$, how can one decompose $V$ as a direct sum of irreducible representations? In the next section we will see that character theory provides a great tool to solve this question in many cases.

**Example 3.1.14.** Let $C_n$ be the cyclic group of order $n \geq 1$, with generator $g$. We wix a primitive $n$-th root of unity $\xi$ in the field $k$. Assume that char$k$ does not divide $n$, then the irreducible $k$-representations of $C_n$ are one-dimensional and given by

$$\rho_j(g) = \xi^j, \quad \text{for } j = 0, 1, \ldots, n-1.$$

**Example 3.1.15.** We consider the binary dihedral group $BD_2$ given by the presentation

$$BD_2 = <a, b\colon a^4 = e,\ b^2 = a^2,\ bab^{-1} = a^{-1}>,$$

where $e$ is the identity. This group has 8 elements and 5 conjugacy classes. The terminology $BD_2$ and the name binary dihedral group will be explained in Chapter 4. It is also called *quaternion group*, in fact it may be identified with the group of the quaternion units $\{\pm 1, \pm i, \pm j, \pm k\}$. One may take, for instance, $i = a$, $j = b$, $k = ab$.

We look for the 5 irreducible representations of $BD_2$ over $\mathbb{C}$ of dimensions $n_0, \ldots, n_4$. For sure we have the trivial representation $(V_0, \rho_0)$, then from Theorem 3.1.13 we have

$$8 = 1 + n_1^2 + n_2^2 + n_3^2 + n_4^2.$$

The only solution in positive integers, up to permutation, is given by $n_2 = n_3 = n_4 = 1$ and $n_1 = 2$. A small search shows that the irreducible one-dimensional representations $(V_j, \rho_j)$ are given by

$$\rho_0(a) = 1,\ \rho_0(b) = 1;$$
$$\rho_2(a) = 1,\ \rho_2(b) = -1;$$
$$\rho_3(a) = -1,\ \rho_3(b) = 1;$$
$$\rho_4(a) = -1,\ \rho_4(b) = -1.$$

The two-dimensional irreducible representation $(V_1, \rho_1)$ is

$$\rho_1(a) = \begin{pmatrix} i & 0 \\ 0 & -i \end{pmatrix},\ \rho_1(b) = \begin{pmatrix} 0 & i \\ i & 0 \end{pmatrix}.$$

## 3.2. Character theory

As remarked in the previous section, characters are an important tool for handling the irreducible representations of a group. We will work with Brauer characters, that may be defined over any algebraically closed field. They are very important to understand representations of a group in the modular situation. However we are only interested in the non-modular case, and in this situation Brauer characters satisfy almost the same properties of ordinary characters, that is characters over $\mathbb{C}$.

Let $G$ be a finite group of order $n$ and let $k$ be an algebraically closed field such that char$k$ does not divide $n$. Let $\mu_n(k)$ be the group of $n$-th roots of unity in $k$ and let $\mu_n(\mathbb{C})$ be the group of complex $n$-th roots of unity. Both $\mu_n(k)$ and $\mu_n(\mathbb{C})$ are cyclic of groups of order $n$, so we can fix an isomorphism

$$\phi : \mu_n(k) \to \mu_n(\mathbb{C}),$$

which we name a *lift*. In the same way, we say that a complex root of unity $z \in \mu_n(\mathbb{C})$ is a lift of $a \in \mu_n(k)$ if $z = \phi(a)$.

Let $(V, \rho)$ be a $k$-representation of $G$ of dimension $r \geq 1$ and let $g$ be an element of $G$. Then the matrix $\rho(g)$ has $r$ eigenvalues in $k$ and in particular they are elements of $\mu_n(k)$, since the order of $g$ divides $n$. Let $\lambda_1, \ldots, \lambda_r$ be these eigenvalues, counted with multiplicity.

**Definition 3.2.1.** The *Brauer character* or simply the *character* of $(V, \rho)$ is the function $\chi : G \to \mathbb{C}$ given by

$$\chi_V(g) = \phi(\lambda_1) + \cdots + \phi(\lambda_r).$$

This definition depend on the choice of the isomorphism $\phi$. Since there are in general many choices for $\phi$, we have a certain degree of arbitrariness. In fact, some of the definitions and results below may depend on the choice of $\phi$. However once chosen $\phi$, it will never be changed and sometimes we will simply say that we lift the eigenvalues to $\mathbb{C}$, meaning that the isomorphism $\phi$ is fixed. The reader should be aware of this fact and not be confused.

**Remark 3.2.2.** If $k = \mathbb{C}$, then we take $\phi$ as the identity map $\mu_n(\mathbb{C}) \to \mu_n(\mathbb{C})$. So the character $\chi$ of a $\mathbb{C}$-representation $(V, \rho)$ can be written as

$$\chi(g) = \lambda_1 + \cdots + \lambda_r = \mathrm{Tr}\rho(g),$$

the trace of the matrix $\rho(g)$. These characters over $\mathbb{C}$ are called *ordinary characters*.

**Example 3.2.3.** Let $k = \overline{\mathbb{F}}_2$ and let $G$ be the cyclic group of order 3. We consider the representation $(V, \rho)$ of $G$ given by

$$\rho(g) = \begin{pmatrix} 0 & 1 \\ 1 & 1 \end{pmatrix},$$

where $g$ is a generator of $G$. The characteristic polynomial of $\rho(g)$ is $t^2 + t + 1$ (recall that we are in characteristic 2), and its roots are the primitive cube roots of unity in $k$. These lift to $\exp(\frac{2\pi i}{3})$ and $\exp(\frac{4\pi i}{3})$, the primitive third roots of unity in $\mathbb{C}$. Therefore the Brauer character of $(V, \rho)$ is

$$\chi(g) = \exp\left(\frac{2\pi i}{3}\right) + \exp\left(\frac{4\pi i}{3}\right) = -1.$$

Notice that the trace of the matrix $\rho(g)$ is 1, which can be lifted to $1 \in \mathbb{C}$. However this does not give the correct Brauer character of $(V, \rho)$.

**Proposition 3.2.4.** *Let $V$ and $W$ be $k$-representations of $G$ with characters $\chi_V$ and $\chi_W$. Then the following facts hold:*

1. *$\chi_V(e) = \dim_k V$, where $e$ is the identity of $G$;*

2. *$\chi_V(g^{-1}) = \overline{\chi_V(g)}$, the complex conjugate of $\chi_V(g)$, for every $g \in G$;*

3. *$\chi_V(hgh^{-1}) = \chi_V(g)$ for every $g, h \in G$;*

4. *$\chi_{V \oplus W}(g) = \chi_V(g) + \chi_W(g)$ for every $g \in G$;*

5. *$\chi_{V \otimes W}(g) = \chi_V(g) \cdot \chi_W(g)$ for every $g \in G$.*

Since the character of the direct sum of two representations is the sum of the characters, it is important to know the characters of the irreducible representations. Moreover characters are constant on conjugacy class, so to determine a character it is sufficient to know its values on conjugacy classes.

These observations suggest to construct a table containing these fundamental data. This is called *character table* and is build as follows. In the top row we list the conjugacy classes of $G$, usually taking a representative from each class, and in the second row we write the number of elements of each class. Sometimes in the third row we may write also the order of the elements in the class. Then the rows are labelled by the irreducible representations of $G$ and in each entry of the table we write the value of the corresponding character on the conjugacy class.

**Example 3.2.5.** We consider the cyclic group $C_n$ of order $n \geq 1$ with generator $g$ and we fix a primitive $n$-th root of unity $\xi \in k$ as in Example 3.1.14. Furthermore we consider a lifting $\phi(\xi) = w := \exp\left(\frac{2\pi i}{n}\right)$. Then the character $\chi_j$ of the irreducible representation $\rho_j(g) = \xi^j$ is given by $\chi_j(g) = w^j$. With this information one can easily write the character table of $C_n$. For example for $n = 3$ we have

| representative | $e$ | $g$ | $g^2$ |
|---|---|---|---|
| $\lvert\,class\,\rvert$ | 1 | 1 | 1 |
| trivial $\rho_0$ | 1 | 1 | 1 |
| $\rho_1$ | 1 | $w$ | $w^2$ |
| $\rho_2$ | 1 | $w^2$ | $w$ |

**Example 3.2.6.** We consider the binary dihedral group $BD_2$ of Example 3.1.15. Its character table over $\mathbb{C}$ is given by

| representative | $e$ | $a^2$ | $a$ | $b$ | $ab$ |
|---|---|---|---|---|---|
| $\mid$ class $\mid$ | 1 | 1 | 2 | 2 | 2 |
| $\rho_0$ | 1 | 1 | 1 | 1 | 1 |
| $\rho_1$ | 2 | $-2$ | 0 | 0 | 0 |
| $\rho_2$ | 1 | 1 | 1 | $-1$ | $-1$ |
| $\rho_3$ | 1 | 1 | $-1$ | 1 | $-1$ |
| $\rho_4$ | 1 | 1 | $-1$ | $-1$ | 1 |

**Remark 3.2.7.** Notice that the first row of the character table, which corresponds to the trivial representation, consists of 1's. On the other hand, in the first column one can read the dimensions of the irreducible representations.

We consider the set $\mathbb{C}_{\text{class}}(G)$ of *class functions* on $G$. These are the functions $G \to \mathbb{C}$ which are constant on conjugacy classes. It follows from Proposition 3.2.4 that characters are class functions. The set $\mathbb{C}_{\text{class}}(G)$ equipped with pointwise sum and multiplication is a $\mathbb{C}$-algebra, isomorphic to $\mathbb{C}^r$, where $r$ is the number of conjugacy classes of $G$.

We define an Hermitian inner product on $\mathbb{C}_{\text{class}}(G)$ by

$$\langle \varphi, \psi \rangle := \frac{1}{|G|} \sum_{g \in G} \overline{\varphi(g)} \psi(g),$$

for every $\varphi, \psi \in \mathbb{C}_{\text{class}}(G)$. This bilinear form satisfies

$$\langle \chi \varphi, \psi \rangle = \langle \chi, \varphi^* \psi \rangle,$$

where $\varphi^*(g) = \overline{\varphi(g)}$ is the class function obtained by complex conjugation.

**Remark 3.2.8.** If we restrict the bilinear form to characters, it is easy to check that

$$\langle \chi, \psi \rangle = \frac{1}{|G|} \sum_{g \in G} \chi(g) \overline{\varphi(g)} = \langle \psi, \chi \rangle,$$

for every $\chi$, $\psi$ characters.

The characters of the irreducible representations of $G$ are orthonormal with respect to this inner product. In fact we have the following theorem.

**Theorem 3.2.9.** *Let $G$ be a finite group and let $k$ be an algebraically closed field such that* char$k$ *does not divide $|G|$. Then the following facts hold.*

1. *A representation $V$ over $k$ is irreducible if and only if $\langle \chi_V, \chi_V \rangle = 1$.*

2. *If $\chi$ and $\psi$ are the characters of two non-isomorphic irreducible $k$-representations then $\langle \chi, \psi \rangle = 0$.*

**Corollary 3.2.10.** *Let $k$ and $G$ be as in Theorem 3.2.9. Let $V$ be a $k$-representation of $G$ and let*

$$V = V_1^{n_1} \oplus \cdots \oplus V_r^{n_r}$$

*be its decomposition into irreducible representations $V_i$. If $\chi_V$ is the character of $V$ and $\chi_{V_i}$ is the character of $V_i$ then*

$$n_i = \langle \chi_V, \chi_{V_i} \rangle.$$

*In particular, two representations are isomorphic if and only if they have the same character.*

**Example 3.2.11.** The orthonormality of the irreducible representations of the groups $C_3$ and $BD_2$ can be read from their character tables in Example 3.2.5 and Example 3.2.6. The numbers over each conjugacy class tell how many times to count entries in that column.

**Remark 3.2.12.** In the modular case, when $\mathrm{char}\,k$ divides $|G|$, the orthonormal relations of Theorem 3.2.9 are no longer true. In this setting, Brauer characters can be defined only for $p$-regular elements, that is elements whose order is not an integer multiple of $p = \mathrm{char}\,k$. However one can still use Brauer characters to investigate the representations of the group, and obtain similar results.

## 3.3. Skew group ring

Now we come back to group actions on rings to investigate the relation between group representations and modules over the invariant ring. As an intermediate step in this correspondence it is convenient to consider modules over the skew group ring, which is a twisted version of the ordinary group ring.

**Definition 3.3.1.** Let $S$ be a (not necessarily commutative) ring and let $G$ a finite subgroup of $\mathrm{Aut}(S)$ such that $|G|$ is invertible in $S$. The *skew group ring of $G$ and $S$* is denoted by $S * G$ and it is

- $S * G = \bigoplus_{g \in G} Sg$, as $S$-module it is free on the elements of $G$;

- the multiplication is given by $(sg)(th) := sg(t) \cdot gh$ for all $s, t \in S$, $g, h \in G$.

Like the ordinary group ring, also the skew group ring is in general non-commutative, even if $S$ is commutative. So, by an $S * G$-module $M$ we will always mean a left $S * G$-module. According to the previous definition, this is just an $S$-module with a compatible action of $G$: $g(sm) = g(s)g(m)$ for all $g \in G$, $s \in S$ and $m \in M$. The ring $S$ is clearly an $S * G$-module, but not every $S$-module has a natural $S * G$-module structure.

Given two $S * G$-modules $M$ and $N$, an $S * G$-module homomorphism $\varphi : M \to N$ is an $S$-module homomorphism respecting the action of $G$, that is $\varphi(g(m)) = g(\varphi(m))$ for all $m \in M$ and $g \in G$. As for the group ring, also in this case it turns out that the $S * G$-module homomorphisms are just the $S$-module homomorphisms that are invariant under the action of $G$. First, we define an $S * G$-module structure on $\operatorname{Hom}_S(M, N)$ by $g(\varphi) := g \circ \varphi \circ g^{-1}$ for all $\varphi \in \operatorname{Hom}_S(M, N)$, $g \in G$. In other words the action of $G$ is expressed by the following commutative diagram

$$
\begin{array}{ccc}
M & \xrightarrow{\varphi} & N \\
g^{-1} \uparrow & & \downarrow g \\
M & \xrightarrow{g(\varphi)} & N.
\end{array}
$$

**Proposition 3.3.2.** *Let $M$ and $N$ be $S * G$-modules. Then*

$$
\operatorname{Hom}_{S*G}(M, N) = \operatorname{Hom}_S(M, N)^G
$$

*Proof:* Let $\varphi \in \operatorname{Hom}_S(M, N)$. If $\varphi$ is $S * G$-linear then for every $m \in M$ and every $g \in G$ we have $\varphi(g(m)) = g(\varphi(m))$. In particular $g(\varphi)(m) = g(\varphi(g^{-1}(m)) = g g^{-1} \varphi(m) = \varphi(m)$, so $\varphi$ is $G$-invariant.

Viceversa let $\varphi$ be $G$-invariant, then $g(\varphi)(m) = \varphi(m)$ for all $m \in M$ and $g \in G$, that is $g(\varphi(g^{-1}(m))) = \varphi(m)$. Applying $g^{-1}$ we get $\varphi(g^{-1}(m)) = g^{-1}\varphi(m)$. $\square$

Also the tensor product $M \otimes_S N$ of two $S * G$-modules $M$ and $N$ has a natural structure of $S * G$-module. This is given by $g(m \otimes n) = g(m) \otimes g(n)$.

**Lemma 3.3.3.** *Let $0 \to M \xrightarrow{\varphi} N \xrightarrow{\psi} L \to 0$ be an exact sequence of $S * G$-modules. Then the sequence*

$$
0 \to M^G \xrightarrow{\overline{\varphi}} N^G \xrightarrow{\overline{\psi}} L^G \to 0
$$

*is exact.*

*Proof:* Clearly $\overline{\varphi}$ is injective since it is a restriction of an injective map. We prove that $\overline{\psi}$ is surjective.

Let $\ell \in L^G$, since $\psi$ is surjective there exists $n \in N$ such that $\psi(n) = \ell$. The order of $G$ is invertible in $S$, so we consider the Reynolds operator and the element

$$
\bar{n} := \frac{1}{|G|} \sum_{g \in G} g(n) \in N^G.
$$

Then by $S * G$-linearity of $\psi$ we have

$$
\begin{aligned}
\overline{\psi}(\bar{n}) &= \psi(\bar{n}) \\
&= \frac{1}{|G|} \sum_{g \in G} \psi(g(n)) \\
&= \frac{1}{|G|} \sum_{g \in G} g(\psi(n)) \\
&= \frac{1}{|G|} \sum_{g \in G} \ell = \ell.
\end{aligned}
$$

The fact that $\mathrm{Im}\,\overline{\varphi} = \ker\overline{\psi}$ follows by a similar argument. □

We immediately get two corollaries.

**Corollary 3.3.4.** *Let $M$ and $N$ be $S * G$-modules. Then for every $i \in \mathbb{N}$ we have*

$$\mathrm{Ext}^i_{S*G}(M, N) = \mathrm{Ext}^i_S(M, N)^G.$$

**Corollary 3.3.5.** *Let $M$ be an $S * G$-module. Then the following facts are are equivalent.*

- *$M$ is projective as $S * G$-module.*
- *$M$ is projective as $S$-module.*

*Proof:* If $M$ is $S$-projective then $\mathrm{Ext}^i_S(M, N) = 0$ for all $i > 0$ and for every $S$-module $N$. In particular if $N$ is an $S * G$-module then $\mathrm{Ext}^i_S(M, N) = 0$ and also the invariant submodule $\mathrm{Ext}^i_S(M, N)^G = 0$ and this coincides with $\mathrm{Ext}^i_{S*G}(M, N)$ by Corollary 3.3.4. So $M$ is $S * G$-projective.

Viceversa let $M$ be $S * G$-projective, then $M$ is a direct summand of a free $S * G$-module $F$. Then $F$ is also $S$-free, so $M$ is a direct summand of a free $S$-module, hence it is $S$-projective. □

Projective $S * G$-modules have a special role. In fact we will see in Section 3.4 that when the ring $S$ is a normal Noetherian local commutative domain over an algebraically closed field $k$, then projective $S * G$-modules correspond to $k[G]$-modules, the group representations.

**Remark 3.3.6.** The ring $S$ can be embedded into $S * G$ in two different ways:

- $S \hookrightarrow S * G$, $s \mapsto s \cdot 1_G$;
- $S \hookrightarrow S * G$ via Reynolds operator, $s \mapsto \bar{\rho}(s) := \frac{1}{|G|} \sum_{g \in G} g(s) g$. The homomorphism $\bar{\rho}$ is injective and its image coincide with the fixed points of $S * G$, namely

$$S \cong \bar{\rho}(S) = (S * G)^G.$$

From now on we will identify $S$ with its image $\bar{\rho}(S)$ inside $S * G$.

Every element of the skew group ring $S * G$ can be viewed as a morphism $S \to S$, in particular this morphism is clearly $R$-linear, where $R = S^G$. In fact we obtain a map

$$\gamma : S * G \to \mathrm{End}_R(S)$$
$$sg \mapsto \gamma(sg),$$

such that $\gamma(sg)(t) := sg(t)$ for every $t \in S$. The map $\gamma$ is a ring homomorphism which extends the group homomorphism $G \to \mathrm{Aut}(S)$ that defines the action of $G$ on $S$.

Even though the map $G \to \mathrm{Aut}(S)$ is injective (if the group action is faithful), the map $\gamma$ is neither injective nor surjective in general. However Auslander [Aus62] found a sufficient condition for the homomorphism $\gamma$ to be bijective.

**Definition 3.3.7.** Let $(A, \mathfrak{m}, K) \subseteq (B, \mathfrak{n}, L)$ be an extension of commutative Noetherian local rings with $\mathfrak{m} = \mathfrak{n} \cap A$ such that $B$ is finitely generated as an $A$-algebra. The extension $A \subseteq B$ is called *unramified* if $\mathfrak{n} = \mathfrak{m}B$ and $K \subseteq L$ is a separable extension of fields.

**Definition 3.3.8.** Let $\varphi : A \hookrightarrow B$ be a ring extension, we say that $\varphi$ is *unramified in codimension one* if $\varphi$ is essentially of finite type and for every prime ideal $\mathfrak{p}$ of height one in $B$ the extension $A_{\mathfrak{p} \cap A} \to B_{\mathfrak{p}}$ is unramified.

**Theorem 3.3.9** (Auslander). *Let $(S, \mathfrak{m})$ be a normal commutative local domain and let $G$ be a finite subgroup of* Aut*$(S)$ such that $|G|$ is invertible in $S$. Let $R = S^G$ be the invariant ring and assume that the inclusion $R \hookrightarrow S$ is unramified in codimension one. Then the ring homomorphism $\gamma : S * G \to \text{End}_R(S)$ given by $\gamma(sg)(t) := sg(t)$ is an isomorphism.*

In addition to the original paper of Auslander [Aus62], the reader may find a proof of Theorem 3.3.9 also in the article of Iyama and Takahashi [IT13, Proposition 4.2], or in Yoshino's book [Yos90, Lemma 10.8] for the two-dimensional case.

## 3.4. The Auslander correspondence

We come now to the heart of the chapter: the proof of the Auslander correspondence between group representations and modules over the invariant ring. This result was first proved by Auslander in his remarkable paper [Aus86b]. Auslander worked over a power series ring in two variables over an algebraically closed field and proved Theorem 3.4.18 below in this setting. His results can be generalized to an $n$-dimensional power series ring $S$, as shown in [LW12, Chapter 5]. However one has to pay a little price for this, namely the equivalence between reflexive modules over the invariant ring $R$ and $R$-direct summands of $S$ of Theorem 3.4.19 is no longer true if $n > 2$.

We will relax further these hypothesis, and require our ring $S$ to be normal, but not necessarily regular. More precisely, we fix the following setting for this section. The condition that the extension $R \hookrightarrow S$ is unramified in codimension one is justified by Theorem 3.3.9 above.

**Situation 3.1.** Let $k$ be an algebraically closed field and let $S$ be a local Noetherian normal commutative $k$-domain with maximal ideal $\mathfrak{m}$ such that $S/\mathfrak{m} = k$. Let $G$ be a finite subgroup of $\text{Aut}_k(S)$, the $k$-linear algebra-automorphisms of $S$, such that $|G|$ is invertible in $k$ and let $R := S^G$ be the invariant ring. We assume that the extension $R \hookrightarrow S$ is unramified in codimension one.

Before proceeding with the Auslander correspondence, in the next theorem we collect some classical results from invariant theory, which will be useful later on. Its proof is taken mainly from [BD08, Theorem 4.1].

**Theorem 3.4.1.** *Let $(S, \mathfrak{m})$ be a Noetherian local normal commutative domain. Let $G$ be a finite group acting on $S$ such that $|G|$ is invertible in $S$ and let $R = S^G$ be the invariant ring. Then the following facts hold.*

1) *R is a Noetherian local domain.*
2) *The extension $R \hookrightarrow S$ is integral, in particular $R$ is normal of Krull dimension* $\dim R = \dim S$.
3) *If $S$ is Cohen-Macaulay then $R$ is Cohen-Macaulay.*
4) *If $S$ is complete then $R$ is complete.*
5) *$S$ is a finitely generated $R$-module of rank $|G|$.*

*Proof: 1)* The fact that $R$ is a domain is obvious. We show that $R$ is local with maximal ideal $\mathfrak{n} := \mathfrak{m} \cap R$. Let $x \in R \setminus \mathfrak{n}$, then $x \notin \mathfrak{m}$ so $x$ is invertible in $S$. Its inverse $x^{-1}$ is clearly $G$-invariant, so $x^{-1} \in R$ and $x$ is invertible in $R$. It follows that $\mathfrak{n}$ is the unique maximal ideal of $R$.

To prove that $R$ is Noetherian it is sufficient to prove the following claim: for every ideal $I$ in $R$ we have $(IS) \cap R = I$. In fact, for any chain of ideals in $R$

$$I_1 \subseteq I_2 \subseteq \cdots \subseteq I_n \subseteq \cdots \subseteq R \tag{3.1}$$

we have an induced chain of ideals in $S$: $I_1 S \subseteq I_2 S \subseteq \cdots \subseteq I_n S \subseteq \cdots \subseteq S$. Since $S$ is Noetherian this chain is stationary, that is there exists a natural number $c$ such that $I_n S = I_m S$ for all $n, m \geq c$. The claim implies that $I_m = (I_m S) \cap R = (I_n S) \cap R = I_n$, so the chain (3.1) is stationary too and $R$ is Noetherian.

Now we prove the claim. The inclusion $I \subseteq (IS) \cap R$ is clear. Let $x \in (IS) \cap R$, so $x = \sum_{i=1}^{m} x_i r_i$ with $x_i \in I$ and $r_i \in S$. Then we have:

$$|G| x = \sum_{g \in G} g(x) = \sum_{i=1}^{m} \sum_{g \in G} g(x_i r_i) = |G| \sum_{i=1}^{m} x_i \frac{1}{|G|} \sum_{g \in G} g(r_i) = |G| \sum_{i=1}^{m} x_i \bar{r}_i,$$

where $\bar{r}_i = \frac{1}{|G|} \sum_{g \in G} g(r_i) \in R$. Hence, $x = \sum_{i=1}^{m} x_i \bar{r}_i \in I$, so $I = (IS) \cap R$.

*2)* Let $a \in S$ and consider the monic polynomial $p(X) := \prod_{g \in G} (X - ag)$. It is easy to check that the coefficients of $p$ are $G$-invariants, hence elements of $R$, and that $p(a) = 0$. It follows that $S$ is integral over $R$ and therefore $\dim R = \dim S$.

Now we prove that $R$ is normal, that is $R$ is integrally closed in its field of fractions $Q(R)$. First we observe that since $|G|$ is invertible in $S$ we have $Q(R) = Q(S)^G$. Let $q \in Q(R)$ be integral over $R$. Then $q$ is also integral over $S$, which is normal, so $q \in S$. Since $Q(R) \cap S = Q(S)^G \cap S = S^G = R$ we obtain that $q \in R$, so $R$ is normal.

*3)* Since $|G|$ is invertible in $S$, the Reynolds operator gives a splitting for the inclusion $R \hookrightarrow S$, as shown in the proof of *1)*. Then $R$ is a direct summand of a Cohen-Macaulay ring, so it is Cohen-Macaulay too.

*4)* We know that $R$ is local with maximal ideal $\mathfrak{n} = \mathfrak{m} \cap R$. Let $(a_n)_{n \geq 1}$ be a Cauchy sequence of elements of $R$ with respect to the $\mathfrak{n}$-adic topology. Since $S$ is complete, there exists a unique element $a \in S$ such that $a \cong a_n \mod \mathfrak{m}^l$ for all $n$ bigger or equal than certain natural number $N_l$. Since for all $g \in G$ we have $g(a_n) = a_n$, for $n \geq N_l$ we get

$$a - g(a) = a - g(a) + a_n - g(a_n) = a - a_n - g(a - a_n) \in \mathfrak{m}^l.$$

By Krull's intersection theorem $\bigcap_{l \geq 1} \mathfrak{m}^l = 0$, therefore $a = g(a)$, so $a \in R$ and $R$ is complete.

*5)* We have already seen that $S$ is integral over $R$, this implies that $S$ is finitely generated as $R$-module. For the rank consider the following commutative diagram

$$
\begin{array}{ccc}
R & \lhook\joinrel\longrightarrow & S \\
\downarrow & & \downarrow \\
Q(R) = Q(S)^G & \lhook\joinrel\longrightarrow & Q(S).
\end{array}
$$

Then we have

$$
\begin{aligned}
\operatorname{rank}_R S &= \dim_{Q(R)} Q(R) \otimes_R S \\
&= \dim_{Q(S)^G} Q(S) \\
&= [Q(S) : Q(S)^G] \\
&= |G|,
\end{aligned}
$$

since $G$ is the Galois group of the field extension $Q(S)^G \hookrightarrow Q(S)$. □

**Remark 3.4.2.** The condition that $|G|$ is invertible in $S$ is necessary, otherwise the ring could be not Noetherian or not Cohen-Macaulay. Nagata [Nag65] constructed an example of a non Noetherian invariant ring under the action of a group $G$ on a regular ring $S$ with $|G|$ non invertible in $S$. Fogarty [Fog81] gave an example of a finite group acting on a Cohen-Macaulay ring $S$ such that the ring of invariants is not Cohen-Macaulay, but it is Noetherian.

### 3.4.1. The functor $\mathscr{F}$

Now assume the setting of Situation 3.1. The first step in the construction of the Auslander correspondence is to establish a relation between modules over the group ring $k[G]$ and modules over the skew group ring $S * G$. It turns out that the category $\operatorname{mod}(S * G)$ is too big in general, so one should rather look at $\operatorname{proj}(S * G)$, the category of finitely generated projective $S * G$-modules.

We recall that by Corollary 3.3.5 an $S * G$-module is $S * G$-projective if and only if it is projective as $S$-module, which is equivalent to be $S$-free, since $S$ is local (Theorem 1.4.4). In particular, it follows that $\operatorname{Hom}_S(P, Q)$ and $P \otimes_S Q$ are $S * G$-projective if $P$ and $Q$ are $S * G$-projective.

We define two functors. First, consider

$$
\begin{aligned}
\mathscr{F} : \operatorname{mod}(k[G]) &\longrightarrow \operatorname{proj}(S * G) \\
V &\longmapsto S \otimes_k V \cong S \otimes_k k^{\dim_k V} \\
(\varphi : V \to W) &\longmapsto \operatorname{id}_S \otimes \varphi,
\end{aligned}
$$

where $\mathrm{id}_S$ is the identity map on $S$. The $S * G$-module structure on $S \otimes_k V$ is given by diagonal action

$$sg(t \otimes v) = (sg(t)) \otimes g(v).$$

Moreover we observe that $sg(\mathrm{id}_S \otimes \varphi) = s \cdot \mathrm{id}_S \otimes g(\varphi)$ for all $g \in G$, $s, t \in S$, $v \in V$, and that the module $S \otimes_k V$ is $S * G$-projective, since it is clearly $S$-free. It follows that the functor $\mathscr{F}$ is well defined.

Furthermore we define

$$\mathscr{F}' : \mathrm{proj}(S * G) \to \mathrm{mod}(k[G])$$
$$P \mapsto P \otimes_S k \cong P/\mathfrak{m}P$$
$$(\psi : P \to Q) \mapsto \psi \mod \mathfrak{m}.$$

The quotient $P/\mathfrak{m}P$ is a finite-dimensional $k$-vector space with a linear action of $G$, so it is a $k[G]$-module. In particular, we point out that the image of the skew group ring is the ordinary group ring, that is $\mathscr{F}'(S * G) = (S * G) \otimes_S k \cong k[G]$.

**Proposition 3.4.3.** *The functors $\mathscr{F}$ and $\mathscr{F}'$ form an adjoint pair, that is*

$$\mathrm{Hom}_{k[G]}(\mathscr{F}'(P), V) = \mathrm{Hom}_{S*G}(P, \mathscr{F}(V))$$

*and they are inverses of each other on objects.*

*Proof:* Let $V$ be a $k[G]$-module, then we have

$$\mathscr{F}'\mathscr{F}(V) = (S \otimes_k V) \otimes_S k \cong k \otimes_k V \cong V.$$

We prove that the other composition is also the identity. We denote by $\pi : S \to k$ the canonical projection modulo $\mathfrak{m}$. Let $P$ be a projective $S * G$-module and consider $\mathscr{F}(\mathscr{F}'(P)) = S \otimes_k P/\mathfrak{m}P$, which is also a projective $S * G$-module. We have the following diagram where the vertical map is given by $\pi \otimes \mathrm{id}_{P/\mathfrak{m}P}$ and the existence of the map $\varphi$ is guaranteed by the fact that $S \otimes_k P/\mathfrak{m}P$ is projective.

We claim that $\varphi$ is surjective. To prove the claim, consider the following commutative diagram with exact rows where $\pi' : P \to P/\mathfrak{m}P$ is the projection modulo $\mathfrak{m}$. The map $P/\mathfrak{m}P \to P/\mathfrak{m}P$ is surjective, then by Nakayama's Lemma also its lifting $\varphi$ is surjective. Since $P$ is projective, the map $\varphi$ splits. In particular $P$ is a direct summand of $S \otimes_k P/\mathfrak{m}P$. Swapping the role of the two modules we get that $S \otimes_k P/\mathfrak{m}P$ is a direct summand of $P$, hence we have the required isomorphism $S \otimes_k P/\mathfrak{m}P \cong P$. $\quad\square$

**Remark 3.4.4.** Note that the categories $\mathrm{proj}(S * G)$ and $\mathrm{mod}(k[G])$ are not equivalent. In fact Hom-sets in $\mathrm{mod}(k[G])$ are $k$-vector spaces, but this is not true in $\mathrm{proj}(S * G)$.

**Corollary 3.4.5.** *Let $S$ and $G$ be as in Situation 3.1. Then the following facts hold.*

*1) The category $\mathrm{proj}(S * G)$ has the Krull-Remak-Schmidt property.*
*2) Let $V_0, V_1, \ldots, V_r$ be a complete set of non-isomorphic simple $k[G]$-modules and let $P_i := \mathscr{F}(V_i)$. Then $P_0, P_1, \ldots, P_r$ is a complete set of non-isomorphic indecomposable finitely generated projective $S * G$-modules.*

The fact that $\mathrm{proj}(S * G)$ is a Krull-Schmidt category allows us to give the following useful characterization.

**Lemma 3.4.6.** *The finitely generated indecomposable projective $S * G$-modules are exactly the indecomposable $S * G$-direct summands of $S * G$.*

*Proof:* The $S * G$-summands of $S * G$ are clearly projective. Conversely, let $P$ be an indecomposable object in $\mathrm{proj}(S * G)$, then $P$ is a direct summand of a finite free $S * G$-module, that is $P \oplus M = (S * G)^n$ for some $S * G$-module $M$. Since $P$ is indecomposable and the category $\mathrm{proj}(S * G)$ has the KRS property, $P$ must be a direct summand of $S * G$ as well, which concludes the proof. $\quad\square$

**Proposition 3.4.7.** *Let $V$ and $W$ be two objects in $\mathrm{mod}(k[G])$. Then the following facts hold.*

*1) $\mathscr{F}(\mathrm{Hom}_k(V, W)) \cong \mathrm{Hom}_S(\mathscr{F}(V), \mathscr{F}(W))$.*
*2) $\mathscr{F}(V \otimes_k W) \cong \mathscr{F}(V) \otimes_S \mathscr{F}(W)$.*
*3) $\mathrm{rank}_S \mathscr{F}(V) = \dim_k V$.*

*Proof: 1)* Consider $\mathrm{Hom}_S(\mathscr{F}(V), \mathscr{F}(W))$, the Hom-tensor adjunction formula yields

$$\mathrm{Hom}_S(S \otimes_k V, S \otimes_k W) \cong \mathrm{Hom}_k(V, \mathrm{Hom}_S(S, S \otimes_k W)).$$

It follows that $\mathrm{Hom}_k(V, \mathrm{Hom}_S(S, S \otimes_k W)) \cong \mathrm{Hom}_k(V, S \otimes_k W)$, and this is isomorphic to $\mathrm{Hom}_k(V, W) \otimes_k S$.

*2)* From ordinary properties of tensor product we have

$$\mathscr{F}(V \otimes_k W) = S \otimes_k (V \otimes_k W) \cong S \otimes_k V \otimes_S S \otimes_k W = \mathscr{F}(V) \otimes_S \mathscr{F}(W),$$

and it is easy to check that the second isomorphism is an isomorphism of $S * G$-modules.

*3)* We have

$$\text{rank}_S\, S \otimes_k V = \dim_{Q(S)} Q(S) \otimes_S S \otimes_k V$$
$$= \dim_{Q(S)} Q(S) \otimes_k V$$
$$= \dim_k V,$$

since $V$ is a free $k$-module. $\square$

### 3.4.2. The functor $\mathscr{G}$

Let $M$ be an $S * G$-module, we say that $M$ is *reflexive* if it is reflexive as $S$-module, that is the map $\lambda_M : M \to M^{**} = \text{Hom}_S(\text{Hom}_S(M,S),S)$ is an isomorphism of $S$-modules. So, in the category $\text{mod}(S * G)$ we consider the dual $(-)^* = \text{Hom}_S(-,S)$ of $\text{mod}(S)$, instead of $\text{Hom}_{S*G}(-,S*G)$. With a little abuse of notation we denote by $\text{Ref}(S * G)$ the full subcategory of $\text{mod}(S * G)$ consisting of reflexive $S * G$-modules. Since an $S * G$-module is $S * G$-projective if and only if it is $S$-projective, and finitely generated projective $S$-modules are reflexive by Corollary 1.3.8, we obtain that $\text{proj}(S * G)$ is a subcategory of $\text{Ref}(S * G)$.

We define two functors:

$$\mathscr{G} : \text{Ref}(S * G) \to \text{Ref}(R)$$
$$M \mapsto M^G,$$
$$(\varphi : M_1 \to M_2) \mapsto (\varphi|_{M_1^G} : M_1^G \to M_2^G),$$

and

$$\mathscr{G}' : \text{Ref}(R) \to \text{Ref}(S * G)$$
$$N \mapsto (S \otimes_R N)^{**},$$
$$(\psi : N_1 \to N_2) \mapsto (\text{id}_S \otimes \psi)^{**},$$

where $(-)^*$ denotes the dual $\text{Hom}_S(-,S)$.

The module $(S \otimes_R N)^{**}$ is reflexive, so the functor $\mathscr{G}'$ is well defined. For $\mathscr{G}$ we need to prove that $M^G$ is reflexive if $M \in \text{Ref}(S * G)$. We consider the presentation

$$(S * G)^m \to (S * G)^\ell \to M^* \to 0,$$

and we apply $\text{Hom}_S(-,S)$ and take $G$-invariants. Since $(-)^G$ is an exact functor by Lemma 3.3.3, and $(S * G)^G = S$ by Remark 3.3.6, we get an exact sequence

$$0 \to M^G \to S^\ell \to S^m.$$

The ring $S$ is normal, in particular it satisfies the $(S_2)$-condition by Theorem 1.4.12. Since it is integral over $R$, $S$ satisfies condition $(S_2)$ also as $R$-module. Applying the Depth Lemma to the previous sequence we get that $M^G$ satisfies $(S_2)$ too, hence it is reflexive by Theorem 1.4.16.

**Lemma 3.4.8.** *Let $M$ be a reflexive $S * G$-module, then we have an isomorphism of $S * G$-modules*

$$M \cong (S \otimes_R M^G)^{**}.$$

*Proof:* We define the map $\psi : (S \otimes_R M^G)^{**} \to M$ to be the composition

$$(S \otimes_R M^G)^{**} \xrightarrow{\delta^{**}} M^{**} \xrightarrow{\lambda^{-1}} M,$$

where $\delta : S \otimes_R M^G \to M$ is given by $\delta(s \otimes m) = sm$ and $\lambda$ is the canonical isomorphism $\lambda : M \to M^{**}$. We prove that $\psi$ is an isomorphism. It is enough to show that $\delta^{**}$ is an isomorphism, and since $(S \otimes_R M^G)^{**}$ and $M$ are reflexive, it suffices to show that $\delta^* : M^* \to (S \otimes_R M^G)^*$ is an isomorphism. We consider a presentation of the $S * G$-module $M^*$ and the following commutative diagram

$$
\begin{array}{ccccccc}
F_1 & \longrightarrow & F_0 & \longrightarrow & M^* & \longrightarrow & 0 \\
\downarrow{\scriptstyle \delta^*} & & \downarrow{\scriptstyle \delta^*} & & \downarrow{\scriptstyle \delta^*} & & \\
\left(S \otimes_R (F_1)^G\right)^* & \longrightarrow & \left(S \otimes_R (F_0)^G\right)^* & \longrightarrow & (S \otimes_R M^G)^* & \longrightarrow & 0,
\end{array}
$$

where $F_0$ and $F_1$ are projective $S * G$-modules. By the Five Lemma and additivity we can assume that $M = S * G$. In this case, we can compose $\delta^*$ with the isomorphisms

$$\left(S \otimes_R (S * G)^G\right)^* \cong (S \otimes_R S)^* \cong \left(S^* \otimes_R S\right)^* \cong (\mathrm{End}_R(S))^*,$$

to obtain a map $(S * G)^* \to (\mathrm{End}_R(S))^*$, which we call again $\delta^*$. By construction the last map is nothing but the dual of the inverse map of the isomorphism $\gamma : S * G \to \mathrm{End}_R(S)$ of Theorem 3.3.9, that is $\delta^* = (\gamma^{-1})^*$. Therefore $\delta^*$ is an isomorphism. □

If the ring $S$ has dimension 2 we can give a more geometric proof of the previous result using the following generalization of the classical Lemma of Speiser [SS88, §92, Aufgabe 16.c., page 763] due to Auslander and Goldman [AG60, Theorem 3.1]. The reader may consult also the paper of Chase, Harrison, and Rosenberg [CHR65, Theorem 1.3].

**Lemma 3.4.9** (Speiser, Auslander-Goldman). *Let $S$ be a commutative domain, $G$ a finite group of ring automorphisms of $S$ and denote by $R = S^G$ the invariant ring. Then the following facts are equivalent.*

*1) $S$ is a separable $R$-algebra.*

*2) For every $g \neq \mathrm{id}_G$ and every maximal ideal $\mathfrak{p}$ in $S$, there exists $s \in S$ such that $s - g(s) \notin \mathfrak{p}$.*

*3) The map $\delta : S \otimes_R M^G \to M$ given by $\delta(s \otimes m) = sm$ is an isomorphism for every $S * G$-module $M$.*

*4) There is a canonical isomorphism $\psi : S \otimes_R S \to \bigoplus_{g \in G} S_g$, where $S_g = S$ for all $g \in G$. The map $\psi$ is given by $\psi(s \otimes t)_g = sg(t) \in S_g$.*

**Remark 3.4.10.** Under the assumptions of Lemma 3.4.9, Chase, Harrison, and Rosenberg [CHR65] proved a ring version of the fundamental theorem of Galois theory. For this reason, one often says that the ring extension $R \subseteq S$ is *Galois with Galois group $G$.*

**Lemma 3.4.11.** *Let S,k, G and R be as in Situation 3.1 and assume that S has Krull dimension* 2*. Let M be a reflexive S ∗ G-module and consider the canonical map*

$$\delta : S \otimes_R M^G \to M,$$

*given by $\delta(s \otimes m) = sm$. Then the following facts hold.*

*1) $\delta$ induces an isomorphism of coherent sheaves $\widetilde{S \otimes_R M^G} \to \widetilde{M}$ on the punctured spectrum $U' = \operatorname{Spec} S \setminus \{\mathfrak{m}\}$.*

*2) $\delta$ induces an isomorphism of S ∗ G-modules*

$$(S \otimes_R M^G)^{**} \cong M.$$

*Proof:* By Theorem 1.5.19, $S$ is a Cohen-Macaulay isolated singularity, so coherent sheaves associated to MCM $S$-modules are locally free on the punctured spectrum. Moreover we recall that reflexive $S$-modules are MCM in this case. Therefore it follows from Lemma 1.5.23 and Corollary 1.5.25 that property *1)* implies property *2)*.

We prove *1)*. Let $f_1, \dots, f_\mu$ be elements in $R$ which generate $\mathfrak{m}$ up to radical, that is $\sqrt{(f_1, \dots, f_\mu)} = \mathfrak{m}$ as ideals in $S$. These elements exist, since the ring extension $R \hookrightarrow S$ is finite (Theorem 3.4.1). Thus, we have an open covering $U' = \bigcup_i D(f_i)$.

We claim that for every $f = f_i$ the induced map

$$\delta_f : S_f \otimes_{R_f} M_f^G \to M_f$$

is an isomorphism. We prove the claim, then the Lemma will follow from it. In fact, because we have a global homomorphism, we get a sheaf homomorphism defined on $U'$ which is locally an isomorphism, so it is forced to be an isomorphism on $U'$.

To prove the claim, let $f$ be one of the $f_i$'s. We check that condition *2)* of Lemma 3.4.9 is true for $S_f$ and $R_f$. Let $\mathfrak{p}$ be a maximal ideal of $S_f$. If $s - g(s) \in \mathfrak{p}$ for every $s \in S_f$ and every $g \in G$, then we have $s \in \mathfrak{p}$ if and only if $g(s) \in \mathfrak{p}$, which is equivalent to say that $\mathfrak{p}$ is a fix-point for the action of $G$ on on $S_f$. On the other hand, since $G$ is small, its action on $U'$ is fix-point-free, so we get a contradiction. Therefore $\delta_f$ is an isomorphism by Lemma 3.4.9 above. □

**Remark 3.4.12.** Let $\mathscr{C}$ and $\mathscr{C}'$ be two abelian categories, and let $F : \mathscr{C} \to \mathscr{C}'$ and $G : \mathscr{C}' \to \mathscr{C}$ be two functors which are quasi-inverse on the objects, that is $G(F(A)) \cong A$ and $F(G(B)) \cong B$ for all $A \in \mathscr{C}$ and $B \in \mathscr{C}'$. If $F \circ G$ and $G \circ F$ are surjective on Hom's, that is $\operatorname{Hom}_{\mathscr{C}}(A, B) \twoheadrightarrow \operatorname{Hom}_{\mathscr{C}}(G(F(A)), G(F(B)))$, and $\operatorname{Hom}_{\mathscr{C}'}(C, D) \twoheadrightarrow \operatorname{Hom}_{\mathscr{C}'}(F(G(C)), F(G(D)))$, then $F$ and $G$ are an equivalence of categories. This is explained by the following commutative diagram

$$\operatorname{Hom}_{\mathscr{C}}(A, B) \longrightarrow \operatorname{Hom}_{\mathscr{C}'}(F(A), F(B)) \longrightarrow \operatorname{Hom}_{\mathscr{C}}(G(F(A)), G(F(B)))$$
$$\Big\downarrow{\scriptstyle\cong}$$
$$\operatorname{Hom}_{\mathscr{C}}(A, B),$$

and by the analogous diagram for $F \circ G$.

**Theorem 3.4.13.** *Let $S$, $G$, $k$ and $R$ be as in Situation 3.1. Then the functors $\mathscr{G} : \mathrm{Ref}(S * G) \to \mathrm{Ref}(R)$, and $\mathscr{G}' : \mathrm{Ref}(R) \to \mathrm{Ref}(S * G)$ give an equivalence of categories*

$$\mathrm{Ref}(S * G) \cong \mathrm{Ref}(R).$$

*Proof:* First, we show that $\mathscr{G}$ and $\mathscr{G}'$ are quasi-inverse on the objects. From Lemma 3.4.8, or Lemma 3.4.11 in dimension 2, we know that $\mathscr{G}'(\mathscr{G}(M)) \cong M$ for every reflexive $S * G$-module $M$.

Let $N$ be an $R$-module, we define the following isomorphism $\varphi : N \to ((S \otimes_R N)^{**})^G$ such that $\varphi(n) = \lambda(1 \otimes n)$, where $\lambda : S \otimes_R N \to (S \otimes_R N)^{**}$ is the canonical map. By Lemma 1.4.14 to prove that $\varphi$ is an isomorphism, it is enough to show that

$$\varphi_\mathfrak{p} : N_\mathfrak{p} \to \left( \left( (S \otimes_R N)^{**} \right)^G \right)_\mathfrak{p} \cong \left( \left( S_\mathfrak{p} \otimes_{R_\mathfrak{p}} N_\mathfrak{p} \right)^{**} \right)^G$$

is an isomorphism for every $\mathfrak{p} \in \mathrm{Spec} R$ with $\dim R_\mathfrak{p} \le 1$. In this case $N_\mathfrak{p}$ is a free $R_\mathfrak{p}$-module, so we may assume $N_\mathfrak{p} = R_\mathfrak{p}$ by additivity. Then the isomorphism is clear.

Now we prove that $\mathscr{G}$ and $\mathscr{G}'$ are surjective on Hom's. Let $f : N_1 \to N_2$ be an arrow in $\mathrm{Ref}(R)$, then from the following commutative diagram

$$
\begin{array}{ccccc}
N_1 & \xrightarrow{\varphi} & ((S \otimes_R N_1)^{**})^G & \xrightarrow{\varphi^{-1}} & N_1 \\
\downarrow{\scriptstyle f} & & \downarrow & & \downarrow{\scriptstyle f'} \\
N_2 & \xrightarrow{\varphi} & ((S \otimes_R N_2)^{**})^G & \xrightarrow{\varphi^{-1}} & N_2
\end{array}
$$

we obtain the existence of an arrow $f'$ in $\mathrm{Ref}(R)$ such that $\mathscr{G}(\mathscr{G}'(f')) = f$. In the same way, if $g : M_1 \to M_2$ is an arrow in $\mathrm{Ref}(S * G)$ and $\psi$ is the isomorphism of Lemma 3.4.8, then the commutative diagram

$$
\begin{array}{ccccc}
M_1 & \xrightarrow{\psi^{-1}} & (S \otimes_R M_1^G)^{**} & \xrightarrow{\psi} & M_1 \\
\downarrow{\scriptstyle g} & & \downarrow & & \downarrow{\scriptstyle g'} \\
M_2 & \xrightarrow{\psi^{-1}} & (S \otimes_R M_2^G)^{**} & \xrightarrow{\psi} & M_2
\end{array}
$$

shows the surjectivity of $\mathscr{G}' \circ \mathscr{G}$ on Hom's. Therefore $\mathscr{G}$ and $\mathscr{G}'$ are an equivalence of categories by the previous Remark 3.4.12. $\square$

**Proposition 3.4.14.** *Let $P$ and $Q$ be reflexive $S * G$-modules. Then the following facts hold.*

*1)* $\mathscr{G} \left( \mathrm{Hom}_S(P, Q) \right) \cong \mathrm{Hom}_R \left( \mathscr{G}(P), \mathscr{G}(Q) \right)$.

*2)* $\mathscr{G}(P \otimes_S Q) \cong \mathscr{G}(P) \boxtimes_R \mathscr{G}(Q)$.

*3)* $\mathrm{rank}_R \mathscr{G}(P) = \mathrm{rank}_S P$.

*Proof: 1)* The fixed point functor $\mathscr{G}$ is an equivalence of categories, so we have

$$\mathrm{Hom}_S(P,Q)^G = \mathrm{Hom}_{S*G}(P,Q) \cong \mathrm{Hom}_R(P^G, Q^G).$$

as required.

*2)* We recall that if $M$ and $N$ are reflexive $R$-modules, then the reflexive tensor product of $M$ and $N$ is defined as $M \boxtimes_R N = (M \otimes_R N)^{**}$ and by Proposition 1.3.12 it satisfies the adjunction $\mathrm{Hom}_R(M \boxtimes_R N, L) \cong \mathrm{Hom}_R(M, \mathrm{Hom}_R(N, L))$ for every $L \in \mathrm{Ref}(R)$.

Let $L$ be a reflexive $R$-module and let $U \in \mathrm{Ref}(S*G)$ such that $U^G = L$. We have

$$\begin{aligned}
\mathrm{Hom}_R\big((P \otimes_S Q)^G, L\big) &= \mathrm{Hom}_R\big((P \otimes_S Q)^G, U^G\big) \\
&= \mathrm{Hom}_{S*G}(P \otimes_S Q, U) \\
&\cong \mathrm{Hom}_{S*G}(P, \mathrm{Hom}_S(Q, U)) \\
&\cong \mathrm{Hom}_R\big(P^G, \mathrm{Hom}_{S*G}(Q, U)\big) \\
&\cong \mathrm{Hom}_R\big(P^G, \mathrm{Hom}_R(Q^G, U^G)\big) \\
&\cong \mathrm{Hom}_R\big(P^G, \mathrm{Hom}_R(Q^G, L)\big) \\
&\cong \mathrm{Hom}_R(P^G \boxtimes_R Q^G, L).
\end{aligned}$$

So we get an isomorphism of $R$-modules

$$\mathrm{Hom}_R\big((P \otimes_S Q)^G, L\big) \cong \mathrm{Hom}_R(P^G \boxtimes_R Q^G, L),$$

which is functorial in $L$. It follows that $(P \otimes_S Q)^G \cong P^G \boxtimes_R Q^G$, as required.

*3)* It is enough to show that $\mathrm{rank}_S \mathscr{G}'(N) = \mathrm{rank}_R N$ for every $N \in \mathrm{Ref}(R)$. Since $\mathscr{G}'(N) = (S \otimes_R N)^{**}$, the last statement is clear. $\square$

We concentrate again on the category of finitely generated projective $S*G$-modules. These are reflexive modules, so we can restrict the functor $\mathscr{G}$ to the full subcategory $\mathrm{proj}(S*G) \subseteq \mathrm{Ref}(S*G)$. We denote this restriction and its quasi inverse again by $\mathscr{G}$ and $\mathscr{G}'$ respectively. It turns out that the image of $\mathrm{proj}(S*G)$ under $\mathscr{G}$ consists exactly of the reflexive $R$-modules which are direct summands of finite free $S$-modules, that is the category $\mathrm{Add}_R(S)$. More precisely we have the following result.

**Proposition 3.4.15.** *Let $S$, $G$ and $R$ as in Situation 3.1. The functors $\mathscr{G} : \mathrm{proj}(S*G) \to \mathrm{Add}_R(S)$ and $\mathscr{G}' : \mathrm{Add}_R(S) \to \mathrm{proj}(S*G)$ induce an equivalence of categories*

$$\mathrm{proj}(S*G) \cong \mathrm{Add}_R(S).$$

*Proof:* We know that $\mathscr{G}$ and $\mathscr{G}'$ are quasi inverse on the objects by Theorem 3.4.13, and that $\mathrm{proj}(S*G)$ is a Krull-Schmidt category by Corollary 3.4.5. Therefore it is enough to show that the indecomposable projective $S * G$-modules correspond exactly to the indecomposable $R$-direct summands of $S$.

Let $M$ be an indecomposable finitely generated projective $S * G$-module, then by Lemma 3.4.6 $M$ is a direct summand of $S * G$, so we have a splitting inclusion $0 \to M \hookrightarrow S * G$. Taking $G$-invariants we obtain an inclusion

$$0 \to M^G \hookrightarrow S,$$

which has a right splitting given by the Reynolds operator. It follows that $M^G \in \mathrm{Add}_R(S)$, which concludes the proof. $\square$

**Corollary 3.4.16.** *Let $S$, $G$ and $R$ be as in Situation 3.1. Then there is a one-one correspondence between*

- *indecomposable objects of $\mathrm{Add}_R(S)$;*
- *indecomposable objects of $\mathrm{proj}(\mathrm{End}_R(S))$;*
- *indecomposable objects of $\mathrm{proj}(S * G)$.*

### 3.4.3. The Auslander functor

We can compose the functors we have considered so far to obtain the desired correspondence between $k$-representations of $G$ and $R$-modules.

**Definition 3.4.17.** Let $S$, $G$, $k$ and $R$ be as in Situation 3.1. We define the *Auslander functor* to be

$$\mathscr{A} = \mathscr{G} \circ \mathscr{F} : \mathrm{mod}(k[G]) \to \mathrm{Add}_R(S)$$
$$V \mapsto (S \otimes_k V)^G,$$

and its right-adjoint

$$\mathscr{A}' = \mathscr{F}' \circ \mathscr{G}' : \mathrm{Add}_R(S) \to \mathrm{mod}(k[G])$$
$$N \mapsto (S \otimes_R N)^{**} \otimes_S k.$$

**Theorem 3.4.18** (Auslander correspondence). *Let $S$, $G$, $k$ and $R$ be as in Situation 3.1. The functors $\mathscr{A} : \mathrm{mod}(k[G]) \to \mathrm{Add}_R(S)$ and $\mathscr{A}' : \mathrm{Add}_R(S) \to \mathrm{mod}(k[G])$ have the following properties.*

1) *$\mathscr{A}(V) \cong \mathscr{A}(W)$ if and only if $V \cong W$.*
2) *$\mathscr{A}(V)$ is indecomposable in $\mathrm{Add}_R(S)$ if and only $V$ is an irreducible representation.*
3) *$\mathscr{A}(\mathrm{Hom}_k(V, W)) \cong \mathrm{Hom}_R(\mathscr{A}(V), \mathscr{A}(W))$ for every $V, W \in \mathrm{mod}(k[G])$.*
4) *$\mathscr{A}(V \otimes_k W) \cong \mathscr{A}(V) \boxtimes_R \mathscr{A}(W)$ for every $V, W \in \mathrm{mod}(k[G])$.*
5) *If $V_0$ is the trivial representation then $\mathscr{A}(V_0) = R$.*
6) *$\mathrm{rank}_R \mathscr{A}(V) = \dim_k V$ for every $k[G]$-module $V$.*

*Proof:* Properties *1)* and *2)* follow from Proposition 3.4.3 and Theorem 3.4.13. Properties *3), 4)*, and *6)* follow from Proposition 3.4.7 and Proposition 3.4.14, while *5)* is a straightforward computation. □

The Auslander functor gives a one-one correspondence between

- irreducible $k$-representations of $G$;
- indecomposable $R$-direct summands of $S$.

In the two-dimensional regular case, Herzog [Her78b] proved that the latter are exactly the indecomposable reflexive $R$-modules.

**Theorem 3.4.19** (Herzog)**.** *Let $S$, $G$, $k$ and $R$ be as in Situation 3.1, and assume that $S$ is regular of Krull dimension* 2*. Then*

$$\mathrm{Ref}(R) = \mathrm{Add}_R(S).$$

*Proof:* We observed already that $\mathrm{Add}_R(S) \subseteq \mathrm{Ref}(R)$ holds in general. In fact, $S$ is $R$-reflexive since it satisfies Serre's condition $(S_2)$, then also a finite sum of copies of $S$ is reflexive. Then reflexive summands of reflexive modules are reflexive, so the inclusion holds.

Conversely, let $M$ be a reflexive $R$-module. We apply the functor $\mathrm{Hom}_R(M^*, -)$ to the inclusion $R \hookrightarrow S$, we get another inclusion

$$M^{**} = \mathrm{Hom}_R(M^*, R) \hookrightarrow \mathrm{Hom}_R(M^*, S).$$

The Reynolds operator gives a splitting for the previous inclusion, hence $M \cong M^{**}$ is a direct summand of $N := \mathrm{Hom}_R(M^*, S)$. From Lemma 1.4.13 we have that $N$ satisfies Serre's condition $(S_2)$ as $R$-module, hence it satisfies it also as $S$-module. Since $S$ is normal, Theorem 1.4.16 implies that $N$ is reflexive over $S$. It follows that $N$ is $S$-free, because $S$ is regular. Therefore $M$ is a direct sum of a finite free $S$-module, that is an object of $\mathrm{Add}_R(S)$. □

**Corollary 3.4.20.** *Let $S$, $G$, $k$ and $R$ be as in Situation 3.1, and assume that $S$ is regular of dimension* 2*. Then the following facts hold.*

*1) The category* $\mathrm{Ref}(R)$ *has the Krull-Remak-Schmidt property.*

*2) Let $P_0, P_1, \dots, P_r$ be a complete set of non-isomorphic indecomposable projective $S * G$-modules and let $N_i := \mathscr{G}(P_i) = P_i^G$. Then $N_0, \dots, N_r$ is a complete set of non-isomorphic indecomposable reflexive $R$-modules. In particular these are exactly the indecomposable $R$-direct summands of $S$.*

*3) $R$ is a Cohen-Macaulay ring of finite CM type.*

*Proof:* By Herzog's Theorem above we have $\mathrm{Ref}(R) = \mathrm{Add}_R(S)$, and by Proposition 3.4.15 $\mathrm{Add}_R(S)$ is equivalent to $\mathrm{proj}(S * G)$, which is a Krull-Schmidt category by Corollary 3.4.5. This shows part *1)*.

Part *2)* follows again from the equivalence of Proposition 3.4.15. The fact that an indecomposable summand of a free *S*-module must be an indecomposable summand of *S* is a consequence of the KRS property.

Finally, *R* is CM by Theorem 1.5.19, and MCM(*R*) = Ref(*R*) in dimension 2, so part *3)* holds. □

**Corollary 3.4.21.** *Let S, G, k and R be as in Situation 3.1, and assume that S is regular of dimension* 2. *Let* $V_0, V_1, \ldots, V_r$ *be a complete set of non-isomorphic irreducible k-representations of G, and fix* $N_i = \mathscr{A}(V_i) = (S \otimes_k V_i)^G$. *Then* $N_0, \ldots, N_r$ *is a complete set of non-isomorphic indecomposable MCM R-modules. Moreover let V be a k-representation, which decomposes as*

$$V = V_0^{n_0} \oplus \cdots \oplus V_r^{n_r},$$

*for some natural numbers* $n_i$. *Then the MCM R-module* $N := \mathscr{A}(V) = (S \otimes_k V)^G$ *decomposes as*

$$N = N_0^{n_0} \oplus \cdots \oplus N_r^{n_r}.$$

A converse to the previous Corollary 3.4.20 is given by the following theorem of Auslander [Aus86b, Theorem 4.9], which says that the invariant rings $S = \mathbb{C}[\![u, v]\!]^G$ are the unique rings of finite CM type over the complex numbers.

**Theorem 3.4.22** (Auslander)**.** *Let R be a two-dimensional normal complete local* $\mathbb{C}$-*domain such that the number of isomorphism classes of indecomposable MCM R-modules is finite. Then there exists a finite group G acting linearly on* $S = \mathbb{C}[\![u, v]\!]$ *such that* $R \cong S^G$.

In the following diagram we summarize the main results of this chapter. We recall that the functors $\mathscr{F}$ and $\mathscr{F}'$ are an adjoint pair, while $\mathscr{G}$ and $\mathscr{G}'$ are equivalence of categories.

$$
\begin{array}{ccccc}
\mathrm{Ref}(S * G) & \underset{\mathscr{G}'}{\overset{\mathscr{G}}{\rightleftarrows}} & \mathrm{Ref}(R) & \overset{\mathrm{dim}2}{=\!=} & \mathrm{MCM}(R) \\
\big\uparrow & & \big\uparrow & & \\
\mathrm{mod}(k[G]) \;\underset{\mathscr{F}'}{\overset{\mathscr{F}}{\rightleftarrows}}\; \mathrm{proj}(S * G) & \underset{\mathscr{G}'}{\overset{\mathscr{G}}{\rightleftarrows}} & \mathrm{Add}_R(S) & &
\end{array}
$$

## 3.5. Quotient singularities

From now on we fix $S = k[\![x_1, \ldots, x_n]\!]$, a power series ring over an algebraically closed field *k* and we consider a finite group *G* acting on *S* and the corresponding invariant ring $R = S^G$. Assume that the order of *G* is invertible in *k*. From Theorem 3.4.1 we know that *R* is a complete, Noetherian, normal, Cohen-Macaulay, local domain of dimension *n*. In particular since *R* is complete, the category of finitely generated *R*-modules mod(*R*) has the KRS property by Theorem 1.1.10.

We would like to apply the theory developed in the previous sections to this situation. To do this, one should assume the following hypothesis:

- the group $G$ acts $k$-linearly on $S$, i.e. $G \subseteq \mathrm{GL}(n, k)$;
- the ring extension $R \hookrightarrow S$ is unramified in codimension one.

In fact these conditions are not so restrictive as it could seem, on the contrary they appear quite naturally in invariant theory as the next results will illustrate.

The following old theorem of Cartan [Car57] allows us to assume that $G$ is acting by linear changes of variables on $S$, in other words $G \subseteq \mathrm{GL}(n, k)$.

**Theorem 3.5.1** (Cartan). *Let $k$ be a field and let $S = k[\![x_1, \dots, x_n]\!]$. Let $G$ be a finite group of $k$-algebras automorphism of $S$ with $|G|$ invertible in $k$. Then there exists a finite group $G_1 \subseteq \mathrm{GL}(n, k)$, acting on $S$ via linear changes of variables, such that $S^G \cong S^{G_1}$.*

*Proof:* Let $\mathfrak{m} = (x_1, \dots, x_n)$ and consider the vector space $V = \mathfrak{m}/\mathfrak{m}^2$, consisting of linear forms of $S$. The action of $G$ on $V$ is linear, so we have a group homomorphism $\varphi : G \to \mathrm{GL}(V)$. We set $G_1 = \varphi(G)$, and we extend the action of $G_1$ to $S$ by linearity. We define a ring homomorphism $\theta : S \to S$

$$\theta(s) := \frac{1}{|G|} \sum_{g \in G} \varphi(g)^{-1} g(s).$$

The restriction of $\theta$ to $V$ is the identity, hence $\theta$ is an automorphism of $S$. For every $h \in G$ we have $\varphi(h) \circ \theta = \theta \circ h$, so the actions of $G$ and $G_1$ are conjugate. It follows that they have isomorphic rings of invariants. $\square$

A consequence of the Theorem of Cartan is that, dealing with invariant rings we may assume without loss of generality that the group $G$ acts linearly on $S$.

**Definition 3.5.2.** Let $k$ be a field and let $G \subseteq \mathrm{GL}(n, k)$ be a finite group. The invariant ring $R = k[\![x_1, \dots, x_n]\!]^G$ is called *quotient singularity*. If $G \subseteq \mathrm{SL}(n, k)$ then $R$ is called *special quotient singularity*. If $G$ is cyclic then $R$ is called *cyclic quotient singularity*.

For a quotient singularity the requirement that the prime ideals of height one are unramified can be translated into a condition on the group $G$, it is equivalent to the group being small.

**Definition 3.5.3.** Let $G \subseteq \mathrm{GL}(n, k)$ be a finite subgroup. An element $g \in G$ is called a *pseudo-reflection* if $\mathrm{rank}(I_n - g) = 1$, where $I_n$ denotes the identity matrix in $\mathrm{GL}(n, k)$. A group $G$ is called *small* if it contains no pseudo-reflections.

**Example 3.5.4.** Subgroups of $\mathrm{SL}(n, k)$ are small. To see this, consider $g \in \mathrm{SL}(n, k)$ such that $\mathrm{rank}(I_n - g) = 1$. This forces $g$ to be a diagonal matrix with diagonal entries given by all 1's and one element $\xi \neq 1$. However this is not possible, since $\det g = 1$.

We present the following characterization without proof, the interested reader shoud consult [LW12, Theorem B.29].

**Theorem 3.5.5.** *Let $k$ be a field, let $S = k[\![x_1, \ldots, x_n]\!]$ and let $G \subseteq \mathrm{GL}(n, k)$ be a finite group acting linearly on $S$ with invariant ring $R = S^G$. Then $R \hookrightarrow S$ is unramified in codimension one if and only if $G$ is small.*

The opposite of a small group is a *reflection group*, that is a group which is generated by pseudo-reflections. From the point of view of invariant theory the action of these groups is irrelevant, as the following classical result of Chevalley, Shephard, and Todd shows.

**Theorem 3.5.6** (Chevalley-Shephard-Todd). *Let $k$ be an algebraically closed field, $S = k[\![x_1, \ldots, x_n]\!]$ and $G$ a finite subgroup of $\mathrm{GL}(n, k)$ such that $|G|$ is invertible in $k$. Then the following facts are equivalent.*

- *$R = S^G$ is a power series ring.*
- *$G$ is generated by pseudo-reflections.*
- *$S$ is free as an $R$-module.*

**Example 3.5.7.** We consider the cyclic group $G$ generated by the pseudo-reflection

$$g = \begin{pmatrix} 1 & 0 & \ldots & 0 & 0 \\ 0 & 1 & \ldots & 0 & 0 \\ & & \ddots & & \\ & & & 1 & 0 \\ & & & 0 & \xi \end{pmatrix},$$

where $\xi$ is a primitive $m$-th root of 1 in $k$. Then we have $k[\![x_1, \ldots, x_n]\!]^G = k[\![x_1, \ldots, x_{n-1}, x_n^m]\!] \cong k[\![y_1, \ldots, y_n]\!]$. So the invariant ring is regular.

As for the condition $G \subseteq \mathrm{GL}(n, k)$, also in this case we may assume without loss of generality that the group is small. This is a result of Prill [Pri67, Proposition 6].

**Theorem 3.5.8** (Prill). *Let $k$ be a field, $S = k[\![x_1, \ldots, x_n]\!]$ and $G$ a finite subgroup of $\mathrm{GL}(n, k)$. Then there exists a small finite subgroup $G_1 \subseteq \mathrm{GL}(n, k)$ such that $S^G \cong S^{G_1}$.*

As we have already noticed, quotient singularities are Cohen-Macaulay rings. Watanabe [Wat74a, Wat74b] proved that the Gorenstein ones are exactly the special quotient singularities.

**Theorem 3.5.9** (Watanabe). *Let $k$ be a field and $G \subseteq \mathrm{GL}(n, k)$ a finite group acting linearly on $S = k[\![x_1, \ldots, x_n]\!]$. Assume that $G$ is small and the order of $G$ is invertible in $k$. Then $S^G$ is Gorenstein if and only if $G \subseteq \mathrm{SL}(n, k)$.*

We can apply the results of the previous Section 3.4 to quotient singularities. For ease of reference we collect them in the following theorem.

**Theorem 3.5.10** (Auslander correspondence). *Let $k$ be an algebraically closed field, let $G$ be a small subgroup of $\mathrm{GL}(n, k)$ acting on $S = k[\![x_1, \ldots, x_n]\!]$ such that $\mathrm{char}\,k$ does not divide $|G|$ and let $R = S^G$ be the corresponding quotient singularity. There is a one-one correspondence between*

*1) indecomposable objects in* $\mathrm{Add}_R(S)$*;*
*2) indecomposable objects in* $\mathrm{proj}(S * G)$*;*
*3) indecomposable objects in* $\mathrm{proj}(\mathrm{End}_R(S))$*;*
*4) indecomposable objects in* $\mathrm{mod}(k[G])$*;*
*5) irreducible $k$-representations of $G$.*

*Explicitly, the correspondence 5) $\to$ 1) is given by the Auslander functor*

$$V \mapsto \mathscr{A}(V) = \left(S \otimes_k V\right)^G;$$

*and the correspondence 1) $\to$ 5) is given by its right-adjoint*

$$N \mapsto \mathscr{A}'(N) = \left(S \otimes_R N\right)^{**} \otimes_S k.$$

*Moreover, if $n = 2$ we have* $\mathrm{Add}_R(S) = \mathrm{Ref}(R) = \mathrm{MCM}(R)$ *and $R$ has finite Cohen-Macaulay type.*

**Remark 3.5.11.** There is also a graded version of Theorem 3.5.10. Instead of a power series ring, one may consider a polynomial ring $S = k[x_1, \ldots, x_n]$ with a $\mathbb{Z}$-grading and restrict to finitely generated $\mathbb{Z}$-graded modules. For a proof and further details in the graded setting the reader may consult the remarkable paper of Iyama and Takahashi [IT13].

**Example 3.5.12.** Let $C_n = \langle g \rangle$ be the cyclic group of Example 3.1.14, with irreducible representations $(V_t, \rho_t)$ over an algebraically closed field $k$ given by $\rho_t(g) = \xi^t$, for a fixed primitive $n$-th root of unity $\xi \in k$. Assume that $\mathrm{char}\, k$ does not divide $|G|$. We can embed $G$ into $\mathrm{GL}(2, k)$ via the representation $V_1 \oplus V_a$, where $a$ is a natural number such that $(a, n) = 1$, otherwise the representation is not faithful. In other words, we consider the cyclic group generated by

$$\begin{pmatrix} \xi & 0 \\ 0 & \xi^a \end{pmatrix}.$$

This group acts linearly on $S = k[\![u, v]\!]$ and the invariant subring $R$ is generated by monomials $u^i v^j$ such that $i + aj \cong 0 \mod n$.

For each irreducible representation $V_t$, we have an indecomposable MCM $R$-module $M_t = (S \otimes_k V_t)^G$. A straightforward computation shows that this is given by

$$M_t = R\left(u^i v^j : i + aj \cong -t \mod n\right).$$

**Example 3.5.13.** We consider the group $BD_2$ of Example 3.1.15 and we embed it into $\mathrm{GL}(2, \mathbb{C})$ through the two-dimensional irreducible faithful representation $(V_1, \rho_1)$. We identify $BD_2$ with its image, generated by the matrices

$$\rho_1(a) = \begin{pmatrix} i & 0 \\ 0 & -i \end{pmatrix}, \ \rho_1(b) = \begin{pmatrix} 0 & i \\ i & 0 \end{pmatrix}.$$

This group acts on $S = \mathbb{C}[\![u, v]\!]$ and the invariant ring is

$$R = \mathbb{C}[\![u^4 + v^4, u^2 v^2, uv(u^4 - v^4)]\!].$$

Following González-Sprinberg and Verdier [GSV81], we can represent the indecomposable MCM $R$-modules as $A$-modules inside $S$, where $A = \mathbb{C}[\![u^4 + v^4, u^2 v^2]\!]$ is the ring of invariants of the group generated by the matrices $\rho(a)$, $\rho(b)$, and $\begin{pmatrix} 0 & 1 \\ 1 & 0 \end{pmatrix}$.

For every irreducible representation $V_t$ (labelled as in Example 3.1.15) the corresponding indecomposable MCM $R$-module $M_t = \mathscr{A}(V_t) = (S \otimes_k V_t)^G$ is given by

$$
\begin{aligned}
M_0 &= R, \\
M_1 &= (u, uv^4, u^2 v, v^3) A, \\
M_2 &= (uv, u^4 - v^4) A, \\
M_3 &= (uv(u^2 - v^2), u^2 + v^2) A \\
M_4 &= (uv(u^2 + v^2), u^2 - v^2) A.
\end{aligned}
$$

One should check that these are actually $R$-modules. To do so, it is enough to multiply the element $g = uv(u^4 - v^4) \in R \setminus A$ with every element in a system of generators for the module $M_t$ and check that the product is again in $M_t$. We do this explicitly for the module $M_2$, and we leave the rest as an exercise for the reader.

A system of generators for $M_2$ as $A$-module is given by $uv$ and $u^4 - v^4$. Thus, we have

$$g \cdot uv = u^2 v^2 (u^4 - v^4) \in M_2,$$

since $u^2 v^2 \in A$, and $u^4 - v^4 \in M_2$. Then, we have

$$g \cdot (u^4 - v^4) = uv(u^4 - v^4)^2 = uv(u^8 + v^8 - 2u^4 v^4) \in M_2,$$

because $uv \in M_2$ and $u^8 + v^8 - 2u^4 v^4 \in A$, since $-2u^4 v^4 \in A$, and $u^8 + v^8 = (u^4 + v^4)^2 - 2u^4 v^4 \in A$.

The following example is due to Auslander and Reiten [AR89, Theorem 4.1], the reader may consult also Yoshino's book [Yos90, Proposition 16.10]. It shows that Theorem 3.4.19 may fail in dimension $\geq 3$.

**Example 3.5.14** (Auslander-Reiten)**.** Let $k$ be an algebraically closed field of characteristic $\neq 2$, and let $G$ be the cyclic group of order 2 acting on $S = k[\![x, y, z]\!]$ by negating each variable. Then the invariant ring $R = k[\![x^2, y^2, z^2, xy, xz, yz]\!]$ is a Cohen-Macaulay ring of finite CM type. Its only indecomposable modules are the ring $R$ itself, the canonical module $K_R = (xy, y^2, yz) R$, and $U = K_R^3 / < (xy, 0, 0), (0, y^2, 0), (0, 0, y^2) >$.

Notice that the modules $R$ and $K_R$ are coming from irreducible representations of $G$, while $U$ is not. We have $R = \mathscr{A}(V_0)$ and $K_R = \mathscr{A}(V_1)$, where $V_0$ is the trivial representation and $V_1$ is the unique non-trivial irreducible representation of $G$.

# 4. Symmetric signature of quotient singularities

The main goal of this chapter is to compute the symmetric signature of two classes of two-dimensional quotient singularites: cyclic quotient singularites and Kleinian singularities. The latter are a particular class of special quotient singularites, which in fact coincide with them if the base field has characteristic zero. Since the strategy of the proof is similar in the two cases, we will give an outline of the main ideas before going on with the results.

Let $k$ be an algebraically closed field, and let $G$ be a small finite subgroup of $\mathrm{GL}(2,k)$ such that $\mathrm{char}\,k$ does not divide $|G|$. The group $G$ acts linearly on $S = k[\![u,v]\!]$ and we denote by $R$ the invariant ring. $R$ is a Cohen-Macaulay normal ring of dimension 2 of finite CM type, so we may consider and compute also the generalized symmetric signature of $R$. Given an indecomposable MCM $R$-module $M_t$, we need to count how many copies of $M_t$ appear as direct summands of the MCM $R$-module

$$\mathscr{S}^q = \left(\mathrm{Sym}_R^q\left(\mathrm{Syz}_R^2(k)\right)\right)^{**}$$

for $q \to +\infty$. We will use the Auslander correspondence and the theory developed in Chapter 3 to translate this algebraic problem into the language of representation theory.

Let $\mathfrak{m}_R$ be the maximal ideal of $R$ and let $p_1,\ldots,p_\mu$ be a minimal system of generators for it. We consider $p_1,\ldots,p_\mu$ as elements of $S$ and we write the beginning of a free resolution of the $S$-module $S/(p_1,\ldots,p_\mu)$

$$0 \to \mathrm{Syz}_S^1(p_1,\ldots,p_\mu) \to S^\mu \xrightarrow{p_1,\ldots,p_\mu} S \to S/(p_1,\ldots,p_\mu) \to 0. \tag{4.1}$$

We observe that we have a natural action of $G$ on this sequence. The group $G$ acts linearly on $S$, and this action is expressed by the fundamental representation $V_1$. This action extends naturally to the free $S$-module $S^\mu$, and therefore also to its submodule $\mathrm{Syz}_S^1(p_1,\ldots,p_\mu)$. We denote by $V$ the representation of $G$ which express its action on $\mathrm{Syz}_S^1(p_1,\ldots,p_\mu)$. From (4.1) it is clear that $\mathrm{rank}_S\,\mathrm{Syz}_S^1(p_1,\ldots,p_\mu) = \mu - 1$, therefore $V$ is a representation of dimension $\mu - 1$.

Actually, we have just observed that the sequence (4.1) is an exact sequence of $S * G$-modules. We apply the exact $G$-invariant functor $\mathscr{G}$ to it and we get

$$0 \to M \to R^\mu \to R \to k \to 0.$$

Since $p_1, \ldots, p_\mu$ is a system of generators of the maximal ideal of $R$, we obtain a copy of the residue field $k$ in the last position. It follows that the module $M$ is the second syzygy of the residue field, in other words $M = \mathrm{Syz}_R^2(k)$. From our construction it is clear that $\mathrm{Syz}_S^1(p_1, \ldots, p_\mu)^G = \mathrm{Syz}_R^2(k)$, and consequently that $M = \mathscr{A}(V) = (S \otimes_k V)^G$.

We will show in Section 4.1 that the Auslander functor $\mathscr{A}$ commutes with reflexive symmetric powers, in particular we have

$$\mathrm{Sym}_R^q(M)^{**} = \mathscr{A}(\mathrm{Sym}_k^q(V)),$$

where $\mathrm{Sym}_k^q(-)$ are the usual symmetric powers of representations. Now let $V_t$ be an irreducible $k$-representation and let $M_t = \mathscr{A}(V_t)$ be the corresponding indecomposable MCM $R$-module. From Corollary 3.4.21 we have that the multiplicity $a_{t,q}$ of the module $M_t$ into $\mathrm{Sym}_R^q(M)^{**}$ is equal to the multiplicity of the representation $V_t$ into $\mathrm{Sym}_k^q(V)$. The latter is much easier to compute, since we may use characters from Section 3.2 to do that.

To procede in this way, one needs to determine the representation $V$ and its character. For the Kleinian singularities it turns out that the ideal $\mathfrak{m}_R$ is minimally generated by 3 elements, so $M = \mathrm{Syz}_R^2(k)$ is an $R$-module of rank 2, and $V$ is a two-dimensional representation. In Section 4.3 we will prove that $V$ is exactly the fundamental representation $V_1$. For cyclic singularities one cannot expect this, since in general $\mu$ may be bigger than 3. However, we will prove in Section 4.5 that $V$ is always a faithful representation in this case, and this is enough to compute the symmetric signature.

The outline of this chapter is the following. As already mentioned, in Section 4.1 we will prove that the Auslander functor commutes with reflexive symmetric powers (Theorem 4.1.8). In Section 4.2 we will present the Kleinian or ADE singularities. First, we will introduce them over $\mathbb{C}$, since this is the context where they appear naturally, and then we will define them over any algebraically closed field. Section 4.3 is dedicated to the study of the syzygy module $\mathrm{Syz}_R^2(k)$ and its associated representation for the Kleinian singularities. Finally, in Section 4.4 and Section 4.5 we will compute the symmetric signature for Kleinian and cyclic singularites respectively.

## 4.1. Auslander functor and symmetric powers

We recall the definition and some properties of the symmetric powers of sheaves. Let $(X, \mathscr{O}_X)$ be a scheme, in the sense of Hartshorne ([Har77, Chapter II, 2.]), and let $\mathscr{F}$ be a sheaf of $\mathscr{O}_X$-modules. We construct a presheaf by assigning to each open set $U$ of $X$ the $\mathscr{O}_X(U)$-module $\mathrm{Sym}^q \mathscr{F}(U)$. The sheaf associated to this presheaf is the *$q$-th symmetric power* of $\mathscr{F}$ and it is denoted by $\mathrm{Sym}_X^q \mathscr{F}$.

**Remark 4.1.1.** Let $R$ be a commutative Noetherian ring and let $M$ be an $R$-module. Then over $X = \operatorname{Spec} R$, we have the following equality of sheaves

$$\operatorname{Sym}_X^q(\widetilde{M}) = \widetilde{\operatorname{Sym}_R^q(M)}.$$

In other words, symmetric powers commute with sheafification.

**Lemma 4.1.2.** *Let $f : X \to Y$ be a morphism of schemes and let $\mathscr{F}$ be a sheaf on $Y$. Then*

$$\operatorname{Sym}_X^q(f^*\mathscr{F}) \cong f^*\left(\operatorname{Sym}_Y^q(\mathscr{F})\right),$$

*as sheaves on $X$.*

*Proof:* See [Har77, Ex. II 5.16]. $\square$

Now, we restrict to the situation of two-dimensional quotient singularites. We fix an algebraically closed field $k$, and a small finite subgroup $G$ of $\operatorname{GL}(2, k)$ such that the characteristic of $k$ does not divide the order of $|G|$. $G$ acts on the power series ring $S = k[\![u, v]\!]$, and we denote by $R$ the corresponding invariant ring. We know that $R$ is a two-dimensional isolated singularity of finite CM type and that $\operatorname{MCM}(R) = \operatorname{Ref}(R) = \operatorname{Add}_R(S)$. We denote by $U = \operatorname{Spec} R \setminus \{\mathfrak{m}\}$ the punctured spectrum of $R$.

**Remark 4.1.3.** Consider the following commutative diagram

$$
\begin{array}{ccc}
U' & \longrightarrow & \operatorname{Spec} S \\
{\scriptstyle \pi}\downarrow & & \downarrow{\scriptstyle \pi} \\
U & \longrightarrow & \operatorname{Spec} R
\end{array}
$$

where the map $\pi$ is induced by the inclusion $R \hookrightarrow S$, and $U' = \pi^{-1}(U)$ is the pull-back of $U$ to $\operatorname{Spec} S$, which is actually the punctured spectrum of $S$. For every $R$-module $M$ this gives us the following identification of sheaves on $U'$

$$\pi^*(\widetilde{M}|_U) \cong \widetilde{S \otimes_R M}\Big|_{U'}. \tag{4.2}$$

We recall that the Auslander functor $\mathscr{A} : \operatorname{mod}(k[G]) \to \operatorname{MCM}(R)$ and its right-adjoint $\mathscr{A}' : \operatorname{MCM}(R) \to \operatorname{mod}(k[G])$ are defined as $\mathscr{A}(V) = (S \otimes_k V)^G$ and $\mathscr{A}'(M) = (S \otimes_R M)^{**} \otimes_S k$.

We want to prove that the Auslander functor commutes with reflexive symmetric powers $\operatorname{Sym}^q(-)^{**}$, that is

$$\mathscr{A}(\operatorname{Sym}_k^q(V)) \cong (\operatorname{Sym}_R^q(\mathscr{A}(V)))^{**}$$

for every $k[G]$-module $V$. Since the Auslander functor is the composition of the functors $\mathscr{F} : \operatorname{mod}(k[G]) \to \operatorname{proj}(S * G)$, $\mathscr{F}(V) = S \otimes_k V$ and $\mathscr{G} : \operatorname{proj}(S * G) \to \operatorname{MCM}(R)$, $\mathscr{G}(N) = N^G$, we will split the proof of this fact in two propositions.

**Remark 4.1.4.** Observe that in the categories mod($k[G]$) and proj($S * G$) the reflexive symmetric powers coincide with the usual symmetric powers. In fact, every finitely generated $k[G]$-module $V$ is reflexive, so it is canonically isomorphic to its double dual $V^{**}$. Since $S$ is regular, for every finitely generated projective $S * G$-module $N$ the symmetric powers $\mathrm{Sym}_S^q(N)$ are $S$-free, hence reflexive.

**Proposition 4.1.5.** *For every $k[G]$-module $V$ we have an isomorphism of $S * G$-modules*

$$\mathrm{Sym}_S^q(S \otimes_k V) \cong S \otimes_k \mathrm{Sym}_k^q(V).$$

*Proof:* We consider the following $S$-multilinear map

$$\psi \colon (S \times V) \times \cdots \times (S \times V) \longrightarrow S \otimes_k \mathrm{Sym}_k^q(V)$$
$$\big((a_1, v_1), \ldots, (a_q, v_q)\big) \mapsto (a_1 \ldots a_q) \otimes (v_1 \circ \cdots \circ v_q).$$

From universal properties of tensor product and symmetric powers functor we obtain a homomorphism $\overline{\psi} \colon \mathrm{Sym}_S^q(S \otimes_k V) \longrightarrow S \otimes_k \mathrm{Sym}_k^q(V)$, given by the following commutative diagram

$$
\begin{array}{ccc}
(S \times V) \times \cdots \times (S \times V) & \xrightarrow{\;\psi\;} & S \otimes_k \mathrm{Sym}_k^q(V) \\
\downarrow & \nearrow & \nearrow \\
(S \otimes_k V) \times \cdots \times (S \otimes_k V) & \overline{\psi} & \\
\downarrow & \nearrow & \\
\mathrm{Sym}_S^q(S \otimes_k V). &  &
\end{array}
$$

The homomorphism $\overline{\psi}$ is in fact an isomorphism (see [Eis94, Proposition A2.2]). It remains to check that $\overline{\psi}$ compatible with the action of $G$, and this is shown by the following commutative diagram

$$
\begin{array}{ccc}
(a_1 \otimes v_1) \circ \cdots \circ (a_q \otimes v_q) & \xrightarrow{\;\overline{\psi}\;} & (a_1 \cdots a_q) \otimes (v_1 \circ \cdots \circ v_q) \\
\downarrow{\scriptstyle g} & & \downarrow{\scriptstyle g} \\
\big(g(a_1) \otimes g(v_1)\big) \circ \cdots \circ \big(g(a_q) \otimes g(v_q)\big) & \xrightarrow{\;\overline{\psi}\;} & \begin{array}{l} g(a_1 \ldots a_q) \otimes g(v_1 \circ \cdots \circ v_q) = \\ \big(g(a_1) \cdots g(a_q)\big) \otimes \big(g(v_1) \circ \cdots \circ g(v_q)\big). \end{array}
\end{array}
$$

$\square$

**Remark 4.1.6.** The statement and the proof of Proposition 4.1.5 are true also if the dimension of $S$ is greater than 2.

**Proposition 4.1.7.** *Let $N$ be a projective $S * G$-module, then we have an isomorphism of $R$-modules*

$$(\mathrm{Sym}_S^q(N))^G \cong \big(\mathrm{Sym}_R^q(N^G)\big)^{**}.$$

*Proof:* Let $\mathscr{G}'(-) = (S \otimes_R -)^{**}$ be the right adjoint of the $G$-invariants functor $\mathscr{G}(-) = (-)^G$. Since $\mathscr{G}$ and $\mathscr{G}'$ are an equivalence of categories (Theorem 3.4.13), it is enough to show that

$$\mathscr{G}'\left((\mathrm{Sym}_S^q(N))^G\right) \cong \mathscr{G}'\left(\left(\mathrm{Sym}_R^q(N^G)\right)^{**}\right),$$

that is

$$\left(S \otimes_R \left(\mathrm{Sym}_S^q(N)\right)^G\right)^{**} \cong \left(S \otimes_R \left(\mathrm{Sym}_R^q(N^G)\right)^{**}\right)^{**}, \tag{4.3}$$

as $S * G$-modules.

Notice that the two double duals on the right hand side of (4.3) are different: the first one is the double dual in $\mathrm{mod}(R)$ and the second is the double dual in $\mathrm{mod}(S)$.

In order to prove (4.3), we consider the left hand side and the right hand side separately. From Lemma 3.4.8, the left hand side of (4.3) is $\left(S \otimes_R \left(\mathrm{Sym}_S^q(N)\right)^G\right)^{**} \cong \mathrm{Sym}_S^q(N)$.

For the $S$-module on the right hand side of (4.3) we use Lemma 1.5.23, and we interprete it as the evaluation of the sheaf

$$\left.\widetilde{S \otimes_R \left(\mathrm{Sym}_R^q(N^G)\right)^{**}}\right|_{U'}$$

on the punctured spectrum $U'$ of $S$.

From the commutative diagram of Remark 4.1.3 we get the isomorphism

$$\left.\widetilde{S \otimes_R \left(\mathrm{Sym}_R^q(N^G)\right)^{**}}\right|_{U'} \cong \pi^*\left(\left.\widetilde{\mathrm{Sym}_R^q(N^G)^{**}}\right|_U\right) = \pi^*\left(\left.\widetilde{\mathrm{Sym}_R^q(N^G)}\right|_U\right),$$

where we can remove the double dual over $U$, thanks to Remark 1.5.24.

Since taking symmetric powers commute with sheafification and with the restriction map of sheaves we get

$$\pi^*\left(\left.\widetilde{\mathrm{Sym}_R^q(N^G)}\right|_U\right) \cong \pi^*\left(\left.\left(\mathrm{Sym}_X^q(\widetilde{N^G})\right)\right|_U\right) \cong \pi^*\left(\mathrm{Sym}_X^q\left(\widetilde{N^G}|_U\right)\right),$$

where the second and the third symmetric powers are sheaf symmetric powers taken over $X = \mathrm{Spec}\,R$. We set $Y = \mathrm{Spec}\,S$, then by Lemma 4.1.2 we have

$$\pi^*\left(\mathrm{Sym}_X^q\left(\widetilde{N^G}|_U\right)\right) \cong \mathrm{Sym}_Y^q\left(\pi^*(\widetilde{N^G}|_U)\right).$$

We apply again (4.2) and Remark 1.5.24 to obtain

$$\mathrm{Sym}_Y^q\left(\pi^*(\widetilde{N^G}|_U)\right) \cong \mathrm{Sym}_Y^q\left(\left.\widetilde{S \otimes_R N^G}\right|_{U'}\right) = \mathrm{Sym}_Y^q\left(\left.\widetilde{(S \otimes_R N^G)^{**}}\right|_{U'}\right).$$

Since taking symmetric powers commutes with sheafification and with the restriction map of sheaves we get

$$\mathrm{Sym}_Y^q\left(\left.\widetilde{\left(S \otimes_R N^G\right)^{**}}\right|_{U'}\right) \cong \left.\widetilde{\mathrm{Sym}_S^q\left(\left(S \otimes_R N^G\right)^{**}\right)}\right|_{U'} \cong \left.\widetilde{\mathrm{Sym}_S^q(N)}\right|_{U'},$$

where the last isomorphism follows from Lemma 3.4.8. Taking global sections on $U'$ we obtain that the right hand side of (4.3) is also isomorphic to $\mathrm{Sym}_S^q(N)$. Therefore we have an isomorphism of $S$-modules as in (4.3). This is actually an isomorphism of $S * G$-modules, because it is obtained from an isomorphism of sheaves on the punctured spectrum $U'$, where the action of $G$ is free. □

From Proposition 4.1.5 and Proposition 4.1.7 we immediately get the following.

**Theorem 4.1.8.** *Let $V$ be a $k[G]$-module, and let $M = (S \otimes_k V)^G$ be the corresponding MCM $R$-module via Auslander functor. Then we have*

$$\mathrm{Sym}_R^q(M)^{**} \cong \left(S \otimes_k \mathrm{Sym}_k^q(V)\right)^G.$$

*In other words $\mathrm{Sym}_R^q(\mathscr{A}(V))^{**} \cong \mathscr{A}(\mathrm{Sym}_k^q(V))$.*

## 4.2. Kleinian singularities

We want to compute the symmetric signature of certain special quotient singularities, called Kleinian singularities. In the literature these are also called with different names: *Du Val singularities, two-dimensional rational double points, ADE surface singularites, simple surface singularities*, to mention some of them. We will use these names indifferently. First, we will work over $\mathbb{C}$, and then we extend the definition to algebraically closed fields whose characteristic does not divide the order of the group.

One of the advantages of working over $\mathbb{C}$ is that two-dimensional special quotient singularites are completely classified. This is a consequence of the classical result that all finite subgroups of $\mathrm{SL}(2, \mathbb{C})$ are known up to isomorphism.

**Theorem 4.2.1** (Klein). *Let $G$ be a finite subgroup of $\mathrm{SL}(2, \mathbb{C})$. Then $G$ is isomorphic to one of the following groups:*

1. *the cyclic group $C_n$ of order $n \geq 1$;*
2. *the binary dihedral group $BD_n$ of order $4n$, for $n \geq 2$;*
3. *the binary tetrahedral group $BT$ of order $24$;*
4. *the binary octahedral group $BO$ of order $48$;*
5. *the binary icosahedral group $BI$ of order $120$.*

**Definition 4.2.2.** The groups of Theorem 4.2.1 are called *Klein groups*, and the corresponding special quotient singularites $\mathbb{C}[\![u, v]\!]^G$ are called *complete complex Kleinian singularities* or *complex ADE singularities*.

We will give an outline of the ideas behind the proof of Theorem 4.2.1, and then we will describe the Klein groups and their quotient singularities. For a more detailed exposition we recommend the notes of Brenner [Bre12].

The classification of Theorem 4.2.1 rely on the classification of finite subgroups of SO(3, $\mathbb{R}$). Traditionally, this is attributed to Hessel, Bravais, and Möbius at the beginning of 19th century. However Jordan was the first one to describe them with the modern language of group theory.

The symmetry groups of the five platonic solids in $\mathbb{R}^3$ (tetrahedron, cube, octahedron, dodecahedron, and icosahedron) are finite subgroups of SO(3, $\mathbb{R}$). Since the middle points of the faces of a octahedron (respectively icosahedron) are a cube (resp. dodecahedron), we get only three distinct groups: the tetrahedral group $T$ of order 12, the octahedral group $O$ of order 24 and the icosahedral group $I$ of order 60. The tetrahedral group is isomorphic to the alternating group $A_4$, the octahedral group is isomorphic to the symmetric group $S_4$, and the icosahedral group is isomorphic to $A_5$.

These are not all finite subgroups of SO(3, $\mathbb{R}$). We have also the cyclic group $C_n$ of order $n$, which is the rotation group of a regular $n$-agon in the plane, and the dihedral group $D_n$ of order $2n$. This group can be understood as the symmetry group of the hosohedron, also called beach ball.

We collect these groups in the following theorem.

**Theorem 4.2.3.** *Let $G$ be a finite subgroup of* SO(3, $\mathbb{R}$). *Then $G$ is isomorphic to one of the following groups:*

1. *the cyclic group $C_n$ of order $n \geq 1$;*
2. *the dihedral group $D_n$ of order $2n$, for $n \geq 2$;*
3. *the tetrahedral group $T \cong A_4$;*
4. *the octahedral group $O \cong S_4$;*
5. *the icosahedral group $I \cong A_5$.*

For a proof of the following two lemmas, see [LW12, Chapter 6] or [Bre12, Vorlesung 23 und 24].

**Lemma 4.2.4.** *Every finite subgroup of* SL(2, $\mathbb{C}$) *is conjugate to a subgroup of* SU(2).

**Lemma 4.2.5.** *There exists a surjective group homomorphism $\pi :$ SU(2) $\rightarrow$ SO(3, $\mathbb{R}$) with kernel $\{\pm I\}$.*

**Remark 4.2.6.** Notice that there is only one element of order 2 in SU(2), namely $-I$. So every matrix of order 2 in SL(2, $\mathbb{C}$) is conjugated to $-I$.

From Lemma 4.2.4 and Lemma 4.2.5, we get the following proposition that allows us to describe the groups of Theorem 4.2.1.

**Proposition 4.2.7.** *Let $G$ be a finite subgroup of* SU(2). *Then either $G$ is cyclic of odd order, or $|G|$ is even and $G = \pi^{-1}(\pi(G))$ is the preimage of a finite subgroup $H = \pi(G)$ of* SO(3).

We can rephrase Proposition 4.2.7 by saying that every finite subgroup $G$ of SU(2) is either cyclic of odd order, or it fits in a non-split short exact sequence

$$0 \rightarrow C_2 \rightarrow G \rightarrow H \rightarrow 0,$$

where $C_2 = \{\pm I\}$ and $H$ is one of the groups of Theorem 4.2.3. In this case we say that $G$ is an extension of $H$ by the cyclic group $C_2$.

Now we are going to describe briefly the Klein groups and the corresponding quotient singularities. In particular, for each group we will describe the one and two dimensional representations and present a minimal set of generators of the maximal ideal of the invariant ring, since these will play an important role in Section 4.3. This information is taken from [GAP15] and [LW12].

### 4.2.1. The singularity $A_{n-1}$

Let $n \in \mathbb{Z}$, $n \geq 2$ and fix $\xi := \exp(\frac{2\pi i}{n})$, a primitive $n$-th root of unity. The cyclic group $C_n$ of order $n$ can be realized as subgroup of $\mathrm{SL}(2,\mathbb{C})$, it is the group generated by the matrix

$$A = \begin{pmatrix} \xi & 0 \\ 0 & \xi^{-1} \end{pmatrix}.$$

The group $C_n$ is abelian, so its irreducible representations $(V_j, \rho_j)$ are all one dimensional. They are given by

$$\rho_j : C_n \to \mathrm{GL}(1,\mathbb{C})$$
$$A \mapsto \xi^j,$$

for $j = 0, \ldots, n-1$ (see also Example 3.1.14, Example 3.2.5, and Example 3.5.12). In particular the fundamental representation given by the matrix $A$ is the representation $V_1 \oplus V_{n-1}$.

The ring of invariants $R = \mathbb{C}[\![u,v]\!]^{C_n}$ is generated by elements

$$X = u^n, \ Y = uv, \ \text{and} \ Z = v^n,$$

which satisfy the relation $Y^n = XZ$. Thus, the ring $R$ is isomorphic to the surface singularity

$$\mathbb{C}[\![x,y,z]\!]/(y^n - xz),$$

which is called $A_{n-1}$-*singularity*.

### 4.2.2. The singularity $D_{n+2}$

The binary tetrahedral group $BD_n$ is an extension of a cyclic group of order $2n$ by the cyclic group $C_2$. It has order $4n$ and $n+3$ conjugacy classes. As subgroup of $\mathrm{SL}(2,\mathbb{C})$ it is the group generated by the following matrices

$$A = \begin{pmatrix} \xi & 0 \\ 0 & \xi^{-1} \end{pmatrix}, \ \text{and} \ B = \begin{pmatrix} 0 & i \\ i & 0 \end{pmatrix},$$

where $\xi := \exp(\frac{2\pi i}{2n})$ and $n \geq 2$.

The group $BD_n$ has 4 one-dimensional irreducible representations, and $n-1$ irreducible two-dimensional representations. The one-dimensional representations $(W_j, \rho'_j)$ are

$$\rho'_0 \colon A \mapsto 1,\ B \mapsto 1,\quad \rho'_1 \colon A \mapsto 1,\ B \mapsto -1,$$
$$\rho'_{n-1} \colon A \mapsto -1,\ B \mapsto i,\quad \rho'_n \colon A \mapsto -1,\ B \mapsto -i.$$

The two-dimensional irreducible representations $(V_j, \rho_j)$ are

$$\rho_j \colon a \to \begin{pmatrix} \xi^j & 0 \\ 0 & \xi^{-j} \end{pmatrix},\ b \mapsto \begin{pmatrix} 0 & i^j \\ i^j & 0 \end{pmatrix} \ \text{ for } j = 1,\dots, n-1.$$

Notice that the representation $V_1$ is the fundamental representation, which gives the action of $BD_n$ on $\mathbb{C}[\![u, v]\!]$. The corresponding invariant ring $R = \mathbb{C}[\![u, v]\!]^{BD_n}$ is generated by polynomials

$$u^{2n} + v^{2n},\ u^2 v^2,\ \text{ and } \ uv(u^{2n} - v^{2n}),$$

if $n$ is even, and polynomials

$$u^{2n} - v^{2n},\ u^2 v^2,\ \text{ and } \ uv(u^{2n} + v^{2n}),$$

if $n$ is odd.

In both cases $R$ is isomorphic to the surface singularity

$$\mathbb{C}[\![x, y, z]\!]/(x^2 + y^{n+1} + yz^2),$$

called $D_{n+2}$-*singularity*. For $n = 2$ we obtain the group $BD_2$ of Example 3.1.15, Example 3.2.6, and Example 3.5.13.

### 4.2.3. The singularity $E_6$

The binary tetrahedral group $BT$ is an extension of the tetrahedral group $T$ by the cyclic group $C_2$. It has order 24 and 7 conjugacy class. Let $\xi := \exp(\frac{2\pi i}{8})$ be a primitive 8-th root of unity, then $BT$ is generated by the matrices

$$A^2 = \begin{pmatrix} i & 0 \\ 0 & -i \end{pmatrix},\ B = \begin{pmatrix} 0 & i \\ i & 0 \end{pmatrix},\ \text{ and } C = \frac{1}{\sqrt{2}} \begin{pmatrix} \xi & \xi^3 \\ \xi & \xi^7 \end{pmatrix}.$$

$BT$ has only two non-trivial normal subgroups:

- the center $C_2 = \{\pm I\}$ of order 2,
- the binary dihedral group $BD_2$ of order 8, generated by $A^2$ and $B$.

The character table for $BT$ is the following

| representative | $I$ | $-I$ | $B$ | $C$ | $C^2$ | $C^4$ | $C^5$ |
|---|---|---|---|---|---|---|---|
| \| class \| | 1 | 1 | 6 | 4 | 4 | 4 | 4 |
| order | 1 | 2 | 4 | 6 | 3 | 3 | 6 |
| $V_0$ | 1 | 1 | 1 | 1 | 1 | 1 | 1 |
| $V_1$ | 2 | -2 | 0 | 1 | -1 | -1 | 1 |
| $V_2$ | 3 | 3 | -1 | 0 | 0 | 0 | 0 |
| $V_3$ | 2 | -2 | 0 | $\zeta_3$ | $-\zeta_3$ | $-\zeta_3^2$ | $\zeta_3^2$ |
| $V_3^\vee$ | 2 | -2 | 0 | $\zeta_3^2$ | $-\zeta_3^2$ | $-\zeta_3$ | $\zeta_3$ |
| $V_4$ | 1 | 1 | 1 | $\zeta_3$ | $\zeta_3$ | $\zeta_3^2$ | $\zeta_3^2$ |
| $V_4^\vee$ | 1 | 1 | 1 | $\zeta_3^2$ | $\zeta_3^2$ | $\zeta_3$ | $\zeta_3$ |

where $\zeta_3 = \frac{-1+\sqrt{3}i}{2}$ is a third root of unity.

We point out that the one-dimensional irreducible representations $V_0$, $V_4$, $V_4^\vee$ form a cyclic group of order 3, which acts on the set of the two-dimensional representations via tensor product. Thus, we can obtain the irreducible representations $V_3$ and $V_3^\vee$ from the fundamental representation $V_1$ via tensorization:

$$V_3 = V_1 \otimes V_4 \ \text{ and } \ V_3^\vee = V_1 \otimes V_4^\vee.$$

The invariant ring $R = \mathbb{C}[\![u, v]\!]^{BT}$ is generated by invariants

$$Z := uv(u^4 - v^4),$$
$$Y := u^8 + 14u^4 v^4 + v^8,$$
$$X := u^{12} - 33u^8 v^4 - 33u^4 v^8 + v^{12},$$

which satisfy the relation $X^2 = Y^3 + 108Z^4$. Adjusting the polynomials by appropriate roots of unity we obtain the isomorphism

$$R \cong \mathbb{C}[\![x, y, z]\!]/(x^2 + y^3 + z^4).$$

This is the *singularity* $E_6$.

### 4.2.4. The singularity $E_7$

The binary octahedral group $BO$ is an extension of the octahedral group $O$ by the cyclic group $C_2$. It is generated by the matrices

$$D = A^3 = \begin{pmatrix} \xi^3 & 0 \\ 0 & \xi^5 \end{pmatrix}, \ B = \begin{pmatrix} 0 & i \\ i & 0 \end{pmatrix}, \ \text{ and } \ C = \frac{1}{\sqrt{2}} \begin{pmatrix} \xi & \xi^3 \\ \xi & \xi^7 \end{pmatrix},$$

where $\xi := \exp(\frac{2\pi i}{8})$. $BO$ has order 48 and 8 conjugacy classes. The normal subgroups of $BO$ are the following

$$C_2 \lhd BD_2 \lhd BT \lhd BO.$$

As usual we denote by $V_1$ the fundamental representation, then the character table of $BO$ is

| representative | $I$ | $-I$ | $B$ | $C$ | $C^2$ | $D$ | $BD$ | $D^3$ |
|---|---|---|---|---|---|---|---|---|
| \|class\| | 1 | 1 | 6 | 8 | 8 | 6 | 12 | 6 |
| order | 1 | 2 | 4 | 6 | 3 | 8 | 4 | 8 |
| $V_0$ | 1 | 1 | 1 | 1 | 1 | 1 | 1 | 1 |
| $V_1$ | 2 | -2 | 0 | 1 | -1 | $-\sqrt{2}$ | 0 | $\sqrt{2}$ |
| $V_2$ | 3 | 3 | -1 | 0 | 0 | 1 | -1 | 1 |
| $V_3$ | 4 | -4 | 0 | -1 | 1 | 0 | 0 | 0 |
| $V_4$ | 3 | 3 | -1 | 0 | 0 | -1 | 1 | -1 |
| $V_5$ | 2 | -2 | 0 | 1 | -1 | $\sqrt{2}$ | 0 | $-\sqrt{2}$ |
| $V_6$ | 1 | 1 | 1 | 1 | 1 | -1 | -1 | -1 |
| $V_7$ | 2 | 2 | 2 | -1 | -1 | 0 | 0 | 0 |

We point out that the irreducible one-dimensional representations form a cyclic group of order 2. The two-dimensional representation $V_5$ is obtained from the fundamental representation via tensorization: $V_5 = V_1 \otimes V_6$. The representation $V_7$ is not faithful, and is given by matrices

$$D \mapsto \begin{pmatrix} 0 & 1 \\ 1 & 0 \end{pmatrix}, \quad B \mapsto \begin{pmatrix} 1 & 0 \\ 0 & 1 \end{pmatrix}, \quad \text{and} \quad C \mapsto \begin{pmatrix} \zeta_3 & 0 \\ 0 & \zeta_3^{-1} \end{pmatrix},$$

with $\zeta_3 = \frac{-1+\sqrt{3}i}{2}$.

Invariant polynomials for the action of $BO$ on $\mathbb{C}[\![u, v]\!]$ are

$$Y = u^2 v^2 (u^4 - v^4)^2,$$
$$Z = u^8 + 14 u^4 v^4 + v^8,$$
$$X = uv(u^4 - v^4)(u^{12} - 33 u^8 v^4 - 33 u^4 v^8 + v^{12}).$$

The corresponding invariant ring is the $E_7$-singularity and is isomorphic to

$$\mathbb{C}[\![x, y, z]\!]/(x^2 + y^3 + yz^3).$$

### 4.2.5. The singularity $E_8$

The binary icosahedral group $BI$ is an extension of the icosahedral group $I$ by the cyclic group $C_2$. The group $BI$ has order 120, has 9 conjugacy classes, and has only one normal non-trivial subgroup, its center $C_2$.

$BI$ is generated by the matrices

$$F = \frac{1}{\sqrt{5}} \begin{pmatrix} \zeta_5^4 - \zeta_5 & \zeta_5^2 - \zeta_5^3 \\ \zeta_5^2 - \zeta_5^3 & \zeta_5 - \zeta_5^4 \end{pmatrix}, \quad \text{and} \quad E = \frac{1}{\sqrt{5}} \begin{pmatrix} \zeta_5^2 - \zeta_5^4 & \zeta_5^4 - 1 \\ 1 - \zeta_5 & \zeta_5^3 - \zeta_5 \end{pmatrix},$$

where $\zeta_5 := \exp(\frac{2\pi i}{5})$ is a primitive 5-th root of unity.

The character table for $BI$ is the following

| representative | $I$ | $-I$ | $F$ | $E$ | $E^2$ | $FE$ | $(FE)^2$ | $(FE)^3$ | $(FE)^4$ |
|---|---|---|---|---|---|---|---|---|---|
| \| class \| | 1 | 1 | 30 | 20 | 20 | 12 | 12 | 12 | 12 |
| order | 1 | 2 | 4 | 6 | 3 | 10 | 5 | 10 | 5 |
| $V_0$ | 1 | 1 | 1 | 1 | 1 | 1 | 1 | 1 | 1 |
| $V_1$ | 2 | -2 | 0 | 1 | -1 | $\varphi^+$ | $-\varphi^-$ | $\varphi^-$ | $-\varphi^+$ |
| $V_2$ | 3 | 3 | -1 | 0 | 0 | $\varphi^+$ | $\varphi^-$ | $\varphi^-$ | $\varphi^+$ |
| $V_3$ | 4 | -4 | 0 | -1 | 1 | 1 | -1 | 1 | -1 |
| $V_4$ | 5 | 5 | 1 | -1 | -1 | 0 | 0 | 0 | 0 |
| $V_5$ | 6 | -6 | 0 | 0 | 0 | -1 | 1 | -1 | 1 |
| $V_6$ | 4 | 4 | 0 | 1 | 1 | -1 | -1 | -1 | -1 |
| $V_7$ | 2 | -2 | 0 | 1 | -1 | $\varphi^-$ | $-\varphi^+$ | $\varphi^+$ | $-\varphi^-$ |
| $V_8$ | 3 | 3 | -1 | 0 | 0 | $\varphi^-$ | $\varphi^+$ | $\varphi^+$ | $\varphi^-$ |

Here $\varphi^+ = \dfrac{1+\sqrt{5}}{2}$ is the golden ratio, $\varphi^- = \dfrac{1-\sqrt{5}}{2}$, and $V_1$ is the fundamental representation.

The two-dimensional representations $V_1$ and $V_7$ cannot be obtained by tensoring with a one-dimensional representation (as it happened in the cases of $BT$ and $BO$), rather they are connected by an automorphism of $BI$.

The invariant ring is generated by polynomials

$$Z = u^{11}v + 11u^6v^6 - uv^{11},$$
$$Y = u^{20} - 228u^{15}v^5 + 494u^{10}v^{10} + 228u^5v^{15} + v^{20},$$
$$X = u^{30} + 522u^{25}v^5 - 10005u^{20}v^{10} - 10005u^{10}v^{20} - 522u^5v^{25} + v^{30},$$

with the relation $X^2 = Y^3 + 1728Z^5$. By adjusting the polynomials by appropriate roots of unity, $\mathbb{C}[\![u,v]\!]^{BI}$ is isomorphic to the surface singularity

$$\mathbb{C}[\![x,y,z]\!]/(x^2 + y^3 + z^5),$$

which is called $E_8$-*singularity*.

### 4.2.6. Kleinian singularities over arbitrary fields

Now let $k$ be an algebraically closed field, we want to define the Kleinian singularities over $k$. The matrices used to generate the Klein groups of Theorem 4.2.1 can be defined over $\mathrm{SL}(2,k)$ as well, provided that the characteristic of $k$ does not divide the order of the group. More precisely, if $|G| = n$ we fix an isomorphism $\psi : \mu_n(\mathbb{C}) \to \mu_n(k)$, and we use $\psi$ to lift the entries of the matrices generating the Klein groups in $\mathbb{C}$ to $k$. The numbers $\sqrt{2}$ and $\sqrt{5}$ appearing in the generators of the groups $BT$, $BO$ and $BI$ can also be lifted to $k$ using the following relations

$$\exp\left(\frac{2\pi i}{8}\right) = \frac{1+i}{\sqrt{2}}, \ \exp\left(\frac{4\pi i}{5}\right) + \exp\left(\frac{6\pi i}{5}\right) = \frac{-1+\sqrt{5}}{2}.$$

**Theorem 4.2.8.** *Let $k$ be an algebraically closed field of characteristic* $0$. *Then the unique finite subgroups of* $\mathrm{SL}(2,k)$ *up to isomorphism are the Klein groups.*

*Proof:* We have already seen that the Klein groups can be realized as subgroups of $\mathrm{SL}(2,k)$. Conversely, let $G$ be a finite subgroup of $\mathrm{SL}(2,k)$. We prove that $G$ is isomorphic to a subgroup of $\mathrm{SL}(2,\mathbb{C})$, and this will force $G$ to be isomorphic to one of the Klein groups by Theorem 4.2.1.

Since $k$ has characteristic $0$, it contains $\mathbb{Q}$ as subfield. Let $L$ be the smallest field extension $\mathbb{Q} \subseteq L \subseteq k$, such that the elements of $G$ are defined over $L$, in other words $G \subseteq \mathrm{SL}(2,L)$. The group $G$ is finite, it follows that $\mathbb{Q} \subseteq L$ is a finite field extension. So there exists an isomorphism $\phi$ between $L$ and a subfield $\mathbb{Q}(\zeta)$ of $\mathbb{C}$. Then we can extend $\phi$ to $\mathrm{SL}(2,L)$ in the natural way, and obtain the required isomorphism between $G$ and a subgroup of $\mathrm{SL}(2,\mathbb{C})$. $\square$

**Definition 4.2.9.** Let $k$ be an algebraically closed field, and let $G$ be a Klein group in $\mathrm{SL}(2,k)$. Assume that $|G|$ is invertible in $k$. A *Kleinian singularity or ADE singularity over $k$* is the invariant ring $k[\![u,v]\!]^G$.

**Remark 4.2.10.** If $k$ has prime characteristic it may not be true that the Klein groups are the unique finite subgroups of $\mathrm{SL}(2,k)$. Thus, there may be other two-dimensional special quotient singularites different from the Kleinian singularities.

The invariant polynomials are the same as the complex case, so the Kleinian singularities over $k$ are hypersurface rings $k[\![x,y,z]\!]/(f)$ as well.

| singularity name | $G$ | $|G|$ | $f$ |
|---|---|---|---|
| $A_{n-1}$ | cyclic | $n$ | $y^n - xz$ |
| $D_{n+2}$ | binary dihedral | $4n$ | $x^2 + y^{n+1} + yz^2$ |
| $E_6$ | binary tetrahedral | $24$ | $x^2 + y^3 + z^4$ |
| $E_7$ | binary octahedral | $48$ | $x^2 + y^3 + yz^3$ |
| $E_8$ | binary icosahedral | $120$ | $x^2 + y^3 + z^5$ |

**Remark 4.2.11.** One could use the hypersurface representation $k[\![x,y,z]\!]/(f)$ to define the Kleinian singularities also in the case where $|G|$ and the characteristic of $k$ are not coprime. However one loses the ability to define the Reynolds operator and the Auslander functor. Since these are very important for our purposes, we will not consider this situation.

**Remark 4.2.12.** The Klein groups are subgroups of SL(2, $k$), so they are small. Therefore the Auslander theory of Chapter 3 applies to the Kleinian singularities.

In the following theorem we summarize the properties of the Kleinian singularities.

**Theorem 4.2.13.** *Let $R = k[\![u, v]\!]^G$ be a Kleinian singularity. Then the following facts hold.*

1. *$R$ is a Gorenstein hypersurface normal domain of dimension 2.*
2. *$R$ has finite Cohen-Macaulay type.*
3. *$R$ is a simple rational isolated singularity of multiplicity $e(R) = 2$.*

## 4.3. The second syzygy of the residue field

In this section we study the maximal Cohen-Macaulay module $\mathrm{Syz}_R^2(k)$ of the Kleinian singularities. In particular we want to determine its associated representation in the Auslander correspondence.

The syzygy modules of the residue field of two-dimensional hypersurfaces have been already investigated by Kawamoto and Yoshino [YK88], and Takahashi [Tak08]. We briefly recall their results and notations.

Let $(R, \mathfrak{m}, k)$ be a complete two-dimensional normal non-regular domain with canonical module $K_R$. In this setting Auslander proved in [Aus86b] that in the category Ref($R$) there exists a unique non-split short exact sequence of the form

$$0 \to K_R \to E \to \mathfrak{m} \to 0,$$

which is called the *fundamental sequence* of $R$.

The module $E$ appearing in the middle term is also unique up to isomorphism and is called the *fundamental module* or *Auslander module* of $R$ (cf. [Yos90, Chapter 11]). It is a reflexive (hence MCM) module of rank 2. Moreover we have an isomorphism of $R$-modules $\left(\bigwedge^2 E\right)^{**} \cong K_R$.

The following example clarifies the name fundamental module.

**Example 4.3.1.** Let $V$ be a $k$-vector space of dimension 2 with basis $u$, $v$ and let $G \subseteq$ GL(2, $k$) be a finite subgroup. Then, $G$ acts on the power series ring $S = k[\![u, v]\!]$ and we consider the invariant ring $R = S^G$, which is a two-dimensional quotient singularity. The fundamental module of $R$ is the image via Auslander functor of the fundamental representation $V$ of $G$, that is

$$E = (S \otimes_k V)^G.$$

Now assume in addition that $R$ is a hypersurface ring over a field $k$, that is $R = T/(f)$, where $T$ is a regular local ring of dimension 3 with a regular system of parameters $\{x, y, z\}$ and $f \in T$. Observe that the Kleinian singularities satisfy these conditions, and that in this setting we have a short exact sequence

$$0 \to \mathrm{Syz}_R^2(k) \to R^3 \xrightarrow{x,y,z} R \to k \to 0,$$

which tells us that $\mathrm{Syz}_R^2(k)$ has rank 2.

Under the previous assumptions, Yoshino and Kawamoto [YK88] proved the following two results.

**Theorem 4.3.2** (Yoshino-Kawamoto)**.** *The fundamental module $E$ is isomorphic to the third syzygy of $k$, i.e.*

$$E \cong \mathrm{Syz}_R^3(k).$$

**Theorem 4.3.3** (Yoshino-Kawamoto)**.** *The following facts are equivalent.*

  1. *The fundamental module $E$ is decomposable.*
  2. *$R$ is a cyclic quotient singularity.*

Now, take $f_x, f_y$ and $f_z$ in $T$ so that they satisfy $f = x f_x + y f_y + z f_z$. The minimal free resolution of the residue field $k$ over $R = T/(f)$ is given as follows

$$\cdots \to R^4 \xrightarrow{C} R^4 \xrightarrow{D} R^4 \xrightarrow{C} R^4 \xrightarrow{B} R^3 \xrightarrow{A} R \to k \to 0,$$

where

$$A = (x, y, z)^t, \; B = \begin{pmatrix} 0 & -z & y \\ z & 0 & -x \\ -y & x & 0 \\ f_x & f_y & f_z \end{pmatrix},$$

$$C = \begin{pmatrix} 0 & f_z & -f_y & x \\ -f_z & 0 & f_x & y \\ f_y & -f_x & 0 & z \\ -x & -y & -z & 0 \end{pmatrix}, \; D = \begin{pmatrix} 0 & -z & y & -f_x \\ z & 0 & -x & -f_y \\ -y & x & 0 & -f_z \\ f_x & f_y & f_z & 0 \end{pmatrix}.$$

In a certain sense, only two modules appear in the previous exact sequence

$$\mathrm{Syz}_R^3(k) = \ker B = \mathrm{im}\, C = \ker D = R^4/\mathrm{im}\, D = \mathrm{Coker}\, D, \text{and}$$
$$\mathrm{Syz}_R^2(k) = \mathrm{im}\, B = R^4/\ker B = R^4/\mathrm{im}\, C = \mathrm{Coker}\, C = \ker C = \mathrm{im}\, D.$$

In fact, the couple of matrices $(C, D)$ is a matrix factorization of the MCM module $\mathrm{Syz}_R^2(k)$, which is equal to the cokernel of $C$, and the couple $(D, C)$ is a matrix factorization of the third syzygy of $k$, $\mathrm{Syz}_R^3(k)$, which is equal to the cokernel of $D$. For the definition and basic properties of matrix factorizations, see [Yos90, Chapter 7].

We also point out that $\mathrm{Syz}_R^1(\mathrm{Syz}_R^2(k)) \cong \mathrm{Syz}_R^3(k)$, or equivalently the fundamental module $E$ is the first cosyzygy of $\mathrm{Syz}_R^2(k)$. Using this relation, Takahashi [Tak08] proved an analogous of Theorem 4.3.3 for the second syzygy of the field.

**Theorem 4.3.4** (Takahashi). *Let $(T, \mathfrak{n}, k)$ be a regular local ring, and $I$ an ideal of $T$ contained in $\mathfrak{n}^2$, and $R = T/I$ the residue class ring. Suppose that $R$ is a Henselian Gorenstein ring of dimension two. Then the following facts are equivalent.*

1. $\mathrm{Syz}^2_R(k)$ *is decomposable.*
2. $\dim T = 3$ *and* $I = (xy - zf)$ *for some regular system of parameters $x, y, z$ of $T$ and $f \in \mathfrak{n}$.*

As a consequence of Theorem 4.3.4 we obtain that the second syzygy module of the field of the Kleinian singularities $D_n$, $E_6$, $E_7$ and $E_8$ is indecomposable, while for the singularity $A_n$ is decomposable. In all cases $\mathrm{Syz}^2_R(k)$ is a MCM $R$-module by Proposition 1.5.16, therefore it corresponds to a $k$-representation of the acting group via Auslander functor (Theorem 3.5.10). We determine explicitely which representation it is.

**Theorem 4.3.5.** *Let $(R, \mathfrak{m}, k)$ be a two-dimensional Kleinian singularity. Then the second syzygy of the residue field $k$ is isomorphic to the fundamental module, that is*

$$\mathrm{Syz}^2_R(k) \cong E.$$

From Example 4.3.1 we immediately get the following corollary.

**Corollary 4.3.6.** *Let $(R, \mathfrak{m}, k)$ be a two-dimensional Kleinian singularity, and let $V_1$ be the two-dimensional fundamental representation of the acting group $G$ of $R$. Then, the second syzygy of the residue field $k$ is isomorphic to the image of $V_1$ via Auslander functor, that is*

$$\mathrm{Syz}^2_R(k) \cong \mathscr{A}(V_1)$$

We prove Theorem 4.3.5 by case considerations. We illustrate the general strategy and prove it for the singularities $A_n$ and $BD_n$. The computations for $E_6$, $E_7$ and $E_8$ are done with Macaulay2 [GS] and left to Appendix A.

Let $S = k[\![u, v]\!]$ be the powers series ring over an algebraically closed field $k$, and let $G$ be one of the Klein subgroups of $\mathrm{SL}(2, k)$ acting on $S$ through a faithful representation $V_1$. Assume that $\mathrm{char}\,k$ and $|G|$ are coprime. We denote by $R = S^G$ the invariant ring, and by $\mathfrak{m}_R$ its maximal ideal. Let $p_1, p_2, p_3$ be a minimal system of generators of $\mathfrak{m}_R$ as ideal in $S$. In fact, these are the polynomials that we denoted by $X$, $Y$, and $Z$ in the previous Section 4.2.

We consider the following short exact sequence of $S$-modules, which is the beginning of an $S$-free resolution of $S/\mathfrak{m}_R S$

$$0 \to \mathrm{Syz}^1_S(p_1, p_2, p_3) \to S^3 \xrightarrow{p_1, p_2, p_3} S \to S/\mathfrak{m}_R S \to 0. \tag{4.4}$$

The group $G$ acts linearly on $S$ through the fundamental representation $V_1$, and this action extends naturally to $S^3$, and its submodule $\mathrm{Syz}^1_S(p_1, p_2, p_3)$. In other words, the sequence (4.4) is an exact sequence of $S * G$-modules. We apply the functor $\mathscr{G}(N) = N^G$, which is exact, to the sequence (4.4) and we obtain an exact sequence

$$0 \to M \to R^3 \to R \to R/\mathfrak{m}_R \to 0.$$

Since $p_1, p_2, p_3$ are a system of generators for the maximal ideal $\mathfrak{m}_R$ of $R$, the last module on the right is a copy of the residue field. It follows that the module $M$ appearing on the left of the last sequence is just the second syzygy of $k$, that is $M = \mathrm{Syz}_R^2(k)$. In other words we have that

$$\mathrm{Syz}_S^1(p_1, p_2, p_3)^G = \mathrm{Syz}_R^2(k).$$

So, to understand which $R$-module $\mathrm{Syz}_R^2(k)$ is, we need to understand which is the action of $G$ on $\mathrm{Syz}_S^1(p_1, p_2, p_3)$, that its $S * G$-module structure.

We know that as $S$-module $\mathrm{Syz}_S^1(p_1, p_2, p_3) \cong S^2$ and it is therefore generated by two elements $s_1(u, v), s_2(u, v) \in S^3$. We need to keep track of the action of $G$ through this isomorphism. The action of $G$ on $\mathrm{Syz}_S^1(p_1, p_2, p_3)$ is inherited by the action on $S$, which is linear and given by matrices $M = M_g$ in $SL(2, k)$. In other words, the $M$'s are the matrices of the fundamental representation of $G$, which we listed for each Klein group in Section 4.2.

In order to understand how these matrices act on the generators $s_1(u, v)$ and $s_2(u, v)$ we procede as follows. For each matrix $M$ we apply the linear transformation $(u, v) \mapsto M(u, v)^T$ to $s_1(u, v)$ and $s_2(u, v)$. We obtain two elements $s_1'(u, v)$ and $s_2'(u, v)$ in $S^3$ which belong to $\mathrm{Syz}_S^1(p_1, p_2, p_3)$. Therefore we can write them as linear combination of $s_1$ and $s_2$

$$\begin{pmatrix} s_1' \\ s_2' \end{pmatrix} = N \begin{pmatrix} s_1 \\ s_2 \end{pmatrix}$$

for some matrix $N$. In this way we obtain a collection of matrices $N = N_g$, which give us the representation corresponding to $\mathrm{Syz}_S^1(p_1, p_2, p_3)$.

We will show that in each case we get a representation isomorphic to the fundamental representation, therefore by Example 4.3.1 we conclude that $\mathrm{Syz}_R^2(k) \cong E$.

### 4.3.1. Computation of $\mathrm{Syz}_R^2(k)$ for the singularity $A_{n-1}$

Let $\xi$ be a primitive $n$-th root of unity in $k$, and consider the cyclic group $C_n$ generated by

$$A = \begin{pmatrix} \xi & 0 \\ 0 & \xi^{-1} \end{pmatrix}.$$

The maximal ideal $\mathfrak{m}_R$ of the invariant ring $R = k[\![u, v]\!]^{C_n}$ is generated by polynomials $p_1 = u^n$, $p_2 = v^n$, $p_3 = uv$. Their syzygy module $\mathrm{Syz}_S^1(p_1, p_2, p_3)$ is generated by

$$s_1 = \begin{pmatrix} 0 \\ -u \\ v^{n-1} \end{pmatrix} \quad \text{and} \quad s_2 = \begin{pmatrix} -v \\ 0 \\ u^{n-1} \end{pmatrix}.$$

We apply the linear transformation given by $A$: $u \mapsto \xi u$, $v \mapsto \xi^{-1} v$, to $s_1$ and $s_2$

$$s_1 \mapsto s_1' = \begin{pmatrix} 0 \\ -\xi u \\ (\xi^{-1})^{n-1} v^{n-1} \end{pmatrix} = \xi \begin{pmatrix} 0 \\ -u \\ v^{n-1} \end{pmatrix} = \xi s_1;$$

$$s_2 \mapsto s_2' = \begin{pmatrix} -\xi^{-1} v \\ 0 \\ \xi^{n-1} u^{n-1} \end{pmatrix} = \xi^{-1} \begin{pmatrix} -v \\ 0 \\ u^{n-1} \end{pmatrix} = \xi^{-1} s_2.$$

Thus, the representation corresponding to $\mathrm{Syz}_S^1(p_1, p_2, p_3)$ is exactly the fundamental representation.

### 4.3.2. Computation of $\mathrm{Syz}_R^2(k)$ for the singularity $D_{n+2}$

Let $\xi \in k$ be a primitive $2n$-th root of unity for $n \geq 2$, and let $i \in k$ be a primitive fourth root of unity. The binary dihedral group $G = BD_n \subseteq \mathrm{SL}(2, k)$ is generated by the matrices

$$A = \begin{pmatrix} \xi & 0 \\ 0 & \xi^{-1} \end{pmatrix} \text{ and } B = \begin{pmatrix} 0 & i \\ i & 0 \end{pmatrix}.$$

If $n$ is even, the maximal ideal $\mathfrak{m}_R$ of the invariant subring $R = k[\![u, v]\!]^{BD_n}$ is generated by polynomials $p_1 = u^{2n} + v^{2n}$, $p_2 = u^2 v^2$ and $p_3 = uv(u^{2n} - v^{2n})$. The syzygy module $\mathrm{Syz}_S^1(p_1, p_2, p_3)$ is generated by

$$s_1 = \begin{pmatrix} -u^2 v \\ 2v^{2n-1} \\ u \end{pmatrix} \text{ and } s_2 = \begin{pmatrix} -uv^2 \\ 2u^{2n-1} \\ -v \end{pmatrix}.$$

We apply the linear transformation $A$: $u \mapsto \xi u$, $v \mapsto \xi^{-1} v$, and we obtain

$$s_1 \mapsto s_1' = \begin{pmatrix} -\xi u^2 v \\ 2(\xi^{-1})^{2n-1} v^{2n-1} \\ \xi u \end{pmatrix} = \xi \begin{pmatrix} -u^2 v \\ 2v^{2n-1} \\ u \end{pmatrix} = \xi s_1;$$

$$s_2 \mapsto s_2' = \begin{pmatrix} -\xi^{-1} u v^2 \\ 2\xi^{2n-1} u^{2n-1} \\ -\xi^{-1} v \end{pmatrix} = \begin{pmatrix} -\xi^{-1} u v^2 \\ 2\xi^{-1} u^{2n-1} \\ -\xi^{-1} v \end{pmatrix} = \xi^{-1} s_2.$$

Thus, we have

$$\begin{pmatrix} s_1' \\ s_2' \end{pmatrix} = \begin{pmatrix} \xi & 0 \\ 0 & \xi^{-1} \end{pmatrix} \begin{pmatrix} s_1 \\ s_2 \end{pmatrix}.$$

For the transformation $B\colon u \mapsto iv$, $v \mapsto iu$, we have

$$s_1 \mapsto s_1' = \begin{pmatrix} -i^3 u v^2 \\ 2i^{2n-1} u^{2n-1} \\ iv \end{pmatrix} = \xi \begin{pmatrix} i u v^2 \\ -2i u^{2n-1} \\ iv \end{pmatrix} = i^3 s_2;$$

$$s_2 \mapsto s_2' = \begin{pmatrix} -i^3 u^2 v \\ 2i^{2n-1} v^{2n-1} \\ -iu \end{pmatrix} = \begin{pmatrix} i u^2 v \\ -2i v^{2n-1} \\ -iu \end{pmatrix} = i^3 s_2.$$

Thus, we have

$$\begin{pmatrix} s_1' \\ s_2' \end{pmatrix} = \begin{pmatrix} 0 & i^3 \\ i^3 & 0 \end{pmatrix} \begin{pmatrix} s_1 \\ s_2 \end{pmatrix}.$$

So the representation $(V', \rho')$ acting on $\mathrm{Syz}_S^1(p_1, p_2, p_3)$ is given by matrices

$$A' = \begin{pmatrix} \xi & 0 \\ 0 & \xi^{-1} \end{pmatrix}, \quad B' = \begin{pmatrix} 0 & i^3 \\ i^3 & 0 \end{pmatrix}.$$

This representation is isomorphic to the fundamental representation $(V_1, \rho_1)$. An isomorphism between them is given by

$$\varphi\colon (V_1, \rho_1) \to (V', \rho')$$
$$(x, y) \mapsto (x, -y),$$

with associated matrix

$$P = \begin{pmatrix} 1 & 0 \\ 0 & -1 \end{pmatrix}.$$

In fact, it is easy to check that $PA' = AP$ and $PB' = BP$.

Now, we consider the case where $n$ is odd. Generators for $\mathfrak{m}_R$ are given by polynomials $p_1 = u^{2n} - v^{2n}$, $p_2 = u^2 v^2$, $p_3 = uv(u^{2n} + v^{2n})$. The syzygy module $\mathrm{Syz}_S^1(p_1, p_2, p_3)$ is generated by

$$s_1 = \begin{pmatrix} u^2 v \\ 2v^{2n-1} \\ -u \end{pmatrix} \quad \text{and} \quad s_2 = \begin{pmatrix} -uv^2 \\ 2u^{2n-1} \\ -v \end{pmatrix}.$$

We apply the linear transformation $A\colon u \mapsto \xi u$, $v \mapsto \xi^{-1} v$, and we obtain

$$s_1 \mapsto s_1' = \begin{pmatrix} \xi u^2 v \\ 2(\xi^{-1})^{2n-1} v^{2n-1} \\ -\xi u \end{pmatrix} = \xi \begin{pmatrix} u^2 v \\ 2v^{2n-1} \\ -u \end{pmatrix} = \xi s_1;$$

$$s_2 \mapsto s_2' = \begin{pmatrix} -\xi^{-1} u v^2 \\ 2\xi^{2n-1} u^{2n-1} \\ -\xi^{-1} v \end{pmatrix} = \begin{pmatrix} -\xi^{-1} u v^2 \\ 2\xi^{-1} u^{2n-1} \\ -\xi^{-1} v \end{pmatrix} = \xi^{-1} s_2.$$

Thus, we have

$$\begin{pmatrix} s_1' \\ s_2' \end{pmatrix} = \begin{pmatrix} \xi & 0 \\ 0 & \xi^{-1} \end{pmatrix} \begin{pmatrix} s_1 \\ s_2 \end{pmatrix} = A \begin{pmatrix} s_1 \\ s_2 \end{pmatrix}.$$

For the transformation $B: u \mapsto iv, \, v \mapsto iu$, we have

$$s_1 \mapsto s_1' = \begin{pmatrix} i^3 u v^2 \\ 2i^{2n-1} u^{2n-1} \\ -iv \end{pmatrix} = \xi \begin{pmatrix} i^3 u v^2 \\ -2i^3 u^{2n-1} \\ i^3 v \end{pmatrix} = -i^3 s_2 = i s_2;$$

$$s_2 \mapsto s_2' = \begin{pmatrix} -i^3 u^2 v \\ 2i^{2n-1} v^{2n-1} \\ -iu \end{pmatrix} = \begin{pmatrix} -i^3 u^2 v \\ -2i^3 v^{2n-1} \\ i^3 u \end{pmatrix} = -i^3 s_2 = i s_2.$$

Thus, we have

$$\begin{pmatrix} s_1' \\ s_2' \end{pmatrix} = \begin{pmatrix} 0 & i \\ i & 0 \end{pmatrix} \begin{pmatrix} s_1 \\ s_2 \end{pmatrix} = B \begin{pmatrix} s_1 \\ s_2 \end{pmatrix}.$$

Therefore the representation acting on $\mathrm{Syz}^1_S(p_1, p_2, p_3)$ is the fundamental representation.

## 4.4. Symmetric signature of Kleinian singularities

In this section we compute the symmetric signature, and the generalized symmetric signature of the Kleinian singularities over an algebraically closed field $k$. As already mentioned, we will translate this into a problem in representation theory using the Auslander correspondence of Chapter 3. For this reason, let us begin with some results in representation theory.

**Lemma 4.4.1.** *Let $\xi \neq \pm 1$ be a root of unity in $\mathbb{C}$. Then the function $f : \mathbb{N} \to \mathbb{C}$*

$$f(q) := \sum_{t=0}^{q} \xi^{2t-q}$$

*is bounded.*

*Proof:* We have

$$\sum_{t=0}^{q} \xi^{2t-q} = \xi^{-q} \sum_{t=0}^{q} (\xi^2)^t.$$

Since $\xi \neq \pm 1$, $\xi^2$ is a root of unity and $\xi^2 \neq 1$. Let $m \geq 2$ be the order of $\xi^2$ and write $q = dm + r$, with $0 \leq r < m$, then

$$\xi^{-q} \sum_{t=0}^{q} (\xi^2)^t = \xi^{-r} \sum_{t=0}^{r} (\xi^2)^t,$$

which is periodic modulo $m$, hence bounded. □

**Example 4.4.2.** Let $(V, \rho)$ be a two-dimensional representation of a finite group $G$ over an algebraically closed field $k$. Since $k$ is algebraically closed, for every element $g \in G$ the matrix $\rho(g)$ can be diagonalized

$$\rho(g) = \begin{pmatrix} \lambda & 0 \\ 0 & \mu \end{pmatrix},$$

with $\lambda, \mu \in k$. The representation $\mathrm{Sym}^q(V)$ evaluated at the element $g$ is given by the matrix

$$\begin{pmatrix} \lambda^q & 0 & \cdots & \cdots & 0 \\ 0 & \lambda^{q-1}\mu & 0 & \cdots & 0 \\ \vdots & & \ddots & & \vdots \\ & \cdots & \cdots & \lambda\mu^{q-1} & 0 \\ 0 & \cdots & \cdots & 0 & \mu^q \end{pmatrix}.$$

In other words, if $\lambda$ and $\mu$ are the eigenvalues of $V$ at the element $g$, then the eigenvalues of $\mathrm{Sym}^q(V)$ at $g$ are $\{\lambda^t \mu^{q-t} : t = 0, \ldots, q\}$.

**Lemma 4.4.3.** *Let $k$ be an algebraically closed field and let $G$ be finite group such that char $k$ and $|G|$ are coprime. Let $(V, \rho)$ be a two-dimensional representation of $G$ whose image is contained in $\mathrm{SL}(2, k)$, then the (Brauer) character of the representation $\mathrm{Sym}^q(V)$ is given by*

$$\chi_{\mathrm{Sym}^q(V)}(g) = \sum_{t=0}^{q} \lambda_g^{2t-q},$$

*where $\lambda_g$ is the lift to $\mathbb{C}$ of an eigenvalue of the matrix $\rho(g)$.*

*Proof:* Let $g \in G$ and let $\lambda_g$ and $\mu_g$ be the lift to $\mathbb{C}$ of the eigenvalues of $\rho(g)$. Since $\rho(g) \in \mathrm{SL}(2, k)$, we have $\mu_g = \lambda_g^{-1}$. Then the lift of the eigenvalues of $g$ in the representation $\mathrm{Sym}^q(V)$ are

$$\{\lambda_g^t \cdot \mu_g^{q-t} : t = 0, \ldots, q\} = \{\lambda_g^{2t-q} : t = 0, \ldots, q\},$$

so the formula for the character of $\mathrm{Sym}^q(V)$ follows immediately. □

**Remark 4.4.4.** The character of the symmetric representation $\mathrm{Sym}^q(V)$ can be computed also using *Molien's formula* (cf. [Web14, Chapter 4, Ex.14])

$$\sum_{q=0}^{+\infty} \chi_{\mathrm{Sym}^q(V)}(g) t^q = \frac{1}{\det(I - t\rho(g))}. \tag{4.5}$$

Here $t$ is an indeterminate, and the determinant on the right hand side is of a matrix with entries conveniently lifted to the polynomial ring $\mathbb{C}[t]$. Expanding the rational function on the right, we obtain a formal power series which is equal to the formal power series on the left, and we can use this equality to compute the character $\chi_{\mathrm{Sym}^q(V)}$. Notice that Molien's formula holds also if $\dim_k V > 2$.

**Theorem 4.4.5.** *Let $k$ be an algebraically closed field and let $G$ be a Klein subgroup of* $\mathrm{SL}(2,k)$, *such that* char$k$ *does not divide* $|G|$. *Let* $(V,\rho)$ *be a two-dimensional faithful $k$-representation of $G$ with image contained in* $\mathrm{SL}(2,k)$ *and let* $(V_i,\rho_i)$ *be an irreducible $k$-representation of $G$. We denote by $\alpha_{i,q}(V)$ the multiplicity of $V_i$ in* $\mathrm{Sym}^q(V)$ *and by $\beta_q(V) = \dim_k \mathrm{Sym}^q(V)$. Then*

$$\lim_{N \to +\infty} \frac{\sum_{q=0}^{N} \alpha_{i,q}(V)}{\sum_{q=0}^{N} \beta_q(V)} = \frac{\dim_k V_i}{|G|}.$$

*Proof:* Since $\dim_k V = 2$, the dimension of $\mathrm{Sym}^q(V)$ is $q+1$. So the denominator in the limit above is $\sum_{q=0}^{N} \beta_q(V) = \frac{1}{2}(N+1)(N+2)$.

From Corollary 3.2.10 and Lemma 4.4.3 we have

$$|G|\alpha_{i,q}(V) = \langle \chi_{\mathrm{Sym}^q(V)}, \chi_{V_i} \rangle = \sum_{g \in G} \chi_{\mathrm{Sym}^q(V)}(g) \cdot \overline{\chi_{V_i}(g)}$$

$$= \sum_{g \in G} \overline{\chi_{V_i}(g)} \left( \sum_{t=0}^{q} \lambda_g^{2t-q} \right),$$

where $\lambda_g$ is the lift to $\mathbb{C}$ of an eigenvalue of $\rho(g)$. The order $m_g$ of the root of unity $\lambda_g$ coincides with the order of $\rho(g)$ in $\mathrm{SL}(2,k)$, since $V$ is faithful this coincides also with the order of the group element $g$ in $G$. If $m_g > 2$ then the sum $\sum_{t=0}^{q} \lambda_g^{2t-q}$ is bounded by Lemma 4.4.1.

If $G$ is cyclic of even order or $G$ is $BD_n$, $BT$, $BO$ or $BI$, then there are only two elements of order $\leq 2$ in $G$, namely $I$ and $-I$. Thus, we can write the previous sum as

$$\sum_{g \in G} \overline{\chi_{V_i}(g)} \left( \sum_{t=0}^{q} \lambda_g^{2t-q} \right) = \overline{\chi_{V_i}(I)} \left( \sum_{t=0}^{q} \lambda_I^{2t-q} \right) + \overline{\chi_{V_i}(-I)} \left( \sum_{t=0}^{q} \lambda_{-I}^{2t-q} \right) + O(1)$$

We have that $\chi_{V_i}(I) = \dim_k V_i$ and $\chi_{V_i}(-I) = \pm \dim_k V_i$, where the sign $\pm$ depends only on the irreducible representation $V_i$. Moreover we have $\lambda_I = 1$ and $\lambda_{-I} = -1$, since $V$ is faithful. Therefore the previous sum is equal to

$$\dim_k(V_i) \left( \sum_{t=0}^{q} 1 \right) \pm \dim_k(V_i) \left( \sum_{t=0}^{q} (-1)^q \right) + O(1)$$

$$= \dim_k(V_i)\big( (q+1)(1 \pm (-1)^q) \big) + O(1),$$

where the term $(1 \pm (-1)^q)$ is equal to 0 or to 2 depending on the parity of $q$ and on the irreducible representation $V_i$. We sum for $q$ running from 0 to a fixed natural number $N$ and we obtain

$$|G| \sum_{q=0}^{N} \alpha_{i,q}(V) = \sum_{q=0}^{N} \big( \dim_k(V_i)\big( (q+1)(1 \pm (-1)^q) \big) + O(1) \big)$$

$$= \frac{1}{2} \dim_k(V_i)(N+1)(N+2) + O(N).$$

Therefore, the numerator of the limit is

$$\sum_{q=0}^{N} \alpha_{i,q}(V) = \frac{(N+1)(N+2)\dim_k V_i}{2|G|} + O(N)$$

$$= \frac{\dim_k V_i}{|G|} \sum_{q=0}^{N} \beta_q(V) + O(N),$$

so the limit is equal to $\frac{\dim_k V_i}{|G|}$ as desired.

If $G$ is a cyclic group of odd order, then there is only one element of order $\leq 2$, the identity $I$, and the computations are analogous to the previous case. $\square$

**Remark 4.4.6.** The multiplicity $\alpha_{i,q}(V)$ can be computed also using Molien's formula (4.5), as showed by Springer in [Spr87]. One obtains

$$\sum_{q=0}^{+\infty} \alpha_{i,q}(V)\, t^q = \frac{1}{|G|} \sum_{g \in G} \frac{\chi_{V_i}(g^{-1})}{\det(I - t\rho(g))}.$$

**Theorem 4.4.7.** *Let $k$ be an algebraically closed field and let $G$ be a Klein subgroup of* SL$(2,k)$ *such that* char$k$ *does not divide* $|G|$. *Let $R$ be the corresponding Kleinian singularity and let $M_i$ be an indecomposable MCM $R$-module. Then the generalized symmetric signature of $R$ with respect to $M_i$ is*

$$s_\sigma(R, M_i) = \frac{\operatorname{rank}_R M_i}{|G|}.$$

*Proof:* We apply the Auslander correspondence of Theorem 3.5.10, and we fix $M = \operatorname{Syz}_R^2(k)$. Let $V_i$ be the irreducible $k$-representation of $G$ such that $M_i = \mathscr{A}(V_i)$, and let $V$ be the $k$-representation such that $\mathscr{A}(V) = M$. Let $\alpha_{i,q}(M)$ be the multiplicity of $M_i$ into $\operatorname{Sym}_R^q(M)^{**}$ and let $\beta_q(M) = \operatorname{rank}_R \operatorname{Sym}_R^q(M)^{**}$.

From Theorem 4.1.8 we have

$$\mathscr{A}'\big(\operatorname{Sym}_R^q(M)^{**}\big) = \operatorname{Sym}^q(V).$$

It follows that $\beta_q(M) = \dim_k \operatorname{Sym}^q(V)$ and $\alpha_{i,q}(M)$ equals the multiplicity of the representation $V_i$ in $\operatorname{Sym}^q(V)$. Since $M$ is the fundamental module by Theorem 4.3.5, $V$ is the fundamental representation. In particular $V$ is faithful and its image is contained in SL$(2,k)$, so we can apply Theorem 4.4.5 to conclude the proof. $\square$

**Corollary 4.4.8.** *Let $k$ be an algebraically closed field and let $G$ be a Klein subgroup of* SL$(2,k)$ *such that* char$k$ *does not divide* $|G|$. *The symmetric signature of the corresponding Kleinian singularity $R$ is*

$$s_\sigma(R) = \frac{1}{|G|}.$$

**Example 4.4.9.** Consider the $A_{n-1}$-singularity $R = k[\![u, v]\!]^{C_n}$ over an algebraically closed field $k$, and let $\mathscr{S}^q = \big(\mathrm{Sym}_R^q(\mathrm{Syz}_R^2(k))\big)^{**}$ as in the definition of symmetric signature. From the proofs of Theorem 4.4.5 and Theorem 4.4.7, we have that $\mathrm{rank}_R \mathscr{S}^q = q + 1$, and

$$\mathrm{frk}_R \mathscr{S}^q = (q+1)(1 + (-1)^q) + O(1) = \begin{cases} 2(q+1) + O(1) & \text{for } q \text{ even} \\ O(1) & \text{for } q \text{ odd.} \end{cases}$$

This shows that the limit

$$\lim_{q \to +\infty} \frac{\mathrm{frk}_R \mathscr{S}^q}{\mathrm{rank}_R \mathscr{S}^q}$$

does not exist. Specifically, in this setting the syzygy module splits as $\mathrm{Syz}_R^2(k) \cong M_1 \oplus M_{n-1}$ (with notation as in Example 3.5.12). Therefore we have

$$\mathscr{S}^q \cong \bigoplus_{t=0}^q M_1^{\otimes t} \otimes M_{n-1}^{\otimes q-t} \cong \bigoplus_{t=0}^q M_1^{\otimes 2t - q}.$$

### 4.4.1. Differential symmetric signature of Kleinian singularities

For the differential symmetric signature of Definition 2.2.15, the situation is very nice and easy. In fact, Martsinkovsky proved in [Mar87] and [Mar90] that for two-dimensional quotient singularities (not necessarily Kleinian or cyclic) the module of Zariski differentials $\Omega_{R/k}^{**}$ is isomorphic to the fundamental module $E$. Therefore with the same proof of Theorem 4.4.7, one obtains that also $s_{d\sigma}(R) = \frac{1}{|G|}$ for a two-dimensional Kleinian singularity $R$.

**Theorem 4.4.10** (Martsinkovsky)**.** *Let $(R, \mathfrak{m}, k)$ be a two-dimensional quotient singularity over an algebraically closed field $k$, and assume that the characteristic of $k$ does not divide the order of the acting group. Then the module of Zariski differentials $\Omega_{R/k}^{**}$ of $R$ over $k$ is isomorphic to the fundamental module $E$.*

**Theorem 4.4.11.** *Let $k$ be an algebraically closed field and let $G$ be a Klein subgroup of $\mathrm{SL}(2, k)$ such that $\mathrm{char}\, k$ does not divide $|G|$. Then the differential symmetric signature of the corresponding Kleinian singularity $R$ is*

$$s_{d\sigma}(R) = \frac{1}{|G|}.$$

## 4.5. Symmetric signature of cyclic quotient singularites

The content of this section is part of a joint work with Lukas Katthän.

We fix a cyclic group $C_n$ of order $n \geq 1$ with generator $g$, and an algebraically closed field $k$, such that $\mathrm{char}\, k$ does not divide $n$. Let $1 \leq a \leq n - 1$, $(a, n) = 1$, and consider the

$k$-representation $V_1 \oplus V_a$, given by

$$\rho(g) = \begin{pmatrix} \xi & 0 \\ 0 & \xi^a \end{pmatrix},$$

where $\xi$ is a fixed primitive $n$-th root of unity in $k$. We denote the cyclic group which is the image of this representation by $\frac{1}{n}(1, a)$. Since $a$ is coprime with $n$, this is a small subgroup of $GL(2, k)$ of order $n$.

Let $S = k[\![u, v]\!]$, then the invariant ring $R$ is generated by the monomials $u^i v^j$ satisfying $i + aj \cong 0 \mod n$. The indecomposable MCM modules have all rank one, and are given by

$$M_t = R\left(u^i v^j : i + aj \cong -t \mod n\right),$$

for $t = 0, \dots, n-1$. If $(V_t, \rho_t)$ is the irreducible representation $\rho_t(g) = \xi^t$, then $M_t = \mathcal{A}(V_t)$. See also Examples 3.1.14, 3.2.5 and 3.5.12.

**Example 4.5.1.** For $a = n-1$, the invariant ring $R$ is the $A_{n-1}$-singularity of Section 4.2.1.

If $a = 1$, then $R = k[\![u^n, u^{n-1}v, \dots, uv^{n-1}, v^n]\!]$ is the Veronese ring.

We want to compute the generalized symmetric signature $s_\sigma(R, M_t)$ for any indecomposable MCM $R$-module $M_t$. The strategy is similar to the one used in Theorem 4.4.7 for the Kleinian singularities, but some differences occur.

First of all, if $M = \mathrm{Syz}_R^2(k)$ and $V$ is the $k$-representation associated to $M$, i.e. $M = \mathcal{A}(V)$, then it is no longer true that $M$ has always rank 2 and $V$ is two-dimensional. We have that $\mathrm{rank}_R \mathrm{Syz}_R^2(k) = \mu - 1$, where $\mu$ is the minimal number of generators of the maximal ideal of $R$, but this number can be bigger than 3. For example for the Veronese ring of Example 4.5.1, we have $\mu = n+1$, therefore $\mathrm{Syz}_R^2(k)$ has rank $n$.

For this reason, $\mathrm{Syz}_R^2(k)$ is not isomorphic in general to the fundamental module $E$, as it happens for Kleinian singularities, since $E$ has always rank 2. Moreover, we cannot hope $\mathrm{Syz}_R^2(k)$ being indecomposable, because the indecomposable modules over $R$ have rank one.

Anyway, Theorem 4.1.8 remains true also in this setting, so the multiplicity of an indecomposable MCM $R$-module $M_t$ in the decomposition of $\mathrm{Sym}_R^q(M)^{**}$ is the same as the multiplicity of $V_t$ in $\mathrm{Sym}^q(V)$, and we can compute the last one using representation theory of finite groups. In particular, we will prove the following theorem.

**Theorem 4.5.2.** *Let $G$ be an abelian group of order $n$, let $\rho : G \to GL(V)$ be a faithful $k$-representation of $G$, and let $\rho_t : G \to GL(V_t)$ be an irreducible $k$-representation of $G$. Assume that the characteristic of $k$ does not divide $n$. We denote by $\alpha_{t,q}(V)$ the multiplicity of $V_t$ in $\mathrm{Sym}^q(V)$ and $\beta_q(V) := \dim_k \mathrm{Sym}^q(V)$. Then we have*

$$\lim_{N \to +\infty} \frac{\sum_{q=0}^N \alpha_{t,q}(V)}{\sum_{q=0}^N \beta_q(V)} = \frac{1}{|G|}.$$

*4. Symmetric signature of quotient singularities*

For the proof of the theorem, we first prove two lemmas. We start with an elementary result from group theory.

**Lemma 4.5.3.** *Let $H \subseteq (\mathbb{Z}/n)^v$ be a subgroup and let $L$ be the kernel of the map*

$$\ell : \mathbb{Z}^v \longrightarrow \operatorname{Hom}_{\mathbb{Z}}(H, \mathbb{Z}/n)$$

*given by $x \mapsto \left( y \mapsto \sum_i x_i y_i \right)$. Then $[\mathbb{Z}^v : L] = \#H$.*

*Proof:* We first show that $\#H = \#\operatorname{Hom}_{\mathbb{Z}}(H, \mathbb{Z}/n)$. $H$ is a finitely generated abelian group, so it splits into a direct sum $H = \bigoplus_i H_i$ of cyclic groups. Moreover, every $H_i$ has order divisible by $n$. The generator of a $H_i$ of order $u$ can be mapped to any element whose order divides $u$ in $\mathbb{Z}/n$, and there are exactly $u$ such elements. Hence

$$\#H = \prod_i \#H_i = \prod_i \#\operatorname{Hom}_{\mathbb{Z}}(H_i, \mathbb{Z}/n) = \#\operatorname{Hom}_{\mathbb{Z}}(H, \mathbb{Z}/n).$$

Next we show that $\ell$ is surjective. We denote by $\langle -, - \rangle$ the standard scalar product, i.e. $\ell(x)(y) = \sum_i x_i y_i = \langle x, y \rangle$. Let $e_1, \ldots, e_v$ be the unit vectors of $(\mathbb{Z}/n)^v$. By the elementary divisor theorem, there exists an invertible matrix $A \in (\mathbb{Z}/n)^{v \times v}$, elements $\alpha_1, \ldots, \alpha_r \in \mathbb{Z}/n$, $r \leq v$, and generators $h_1, \ldots, h_r$ of $H$, such that $Ah_i := \alpha_i e_i$ for $1 \leq i \leq r$. It is clear that $\operatorname{Hom}_{\mathbb{Z}}(H, \mathbb{Z}/n)$ is generated by the maps $\varphi_i$ sending $h_j$ to $\alpha_j$ for $j = i$ and all other generators to zero. But the map $\varphi_i$ is the image under $\ell$ of (an arbitrary lifting to $\mathbb{Z}^v$ of) the vector $A^t e_i$, where $A^t$ denotes the transpose of $A$:

$$\ell(A^t e_i)(h_j) = \langle A^t e_i, h_j \rangle = \langle e_i, A h_j \rangle$$
$$= \langle e_i, \alpha_j e_j \rangle = \alpha_j \delta_{ij} = \varphi_i(h_j),$$

where $\delta_{ij}$ is the Kronecker symbol. So the claim follows. $\square$

**Lemma 4.5.4.** *Let $G$, $(V, \rho)$, and $(V_t, \rho_t)$ be as in Theorem 4.5.2, and let $v = \dim_k V$. Then there exists a lattice $L \subseteq \mathbb{Z}^v$ and $a_0 \in \mathbb{Z}^v$ such that*

$$\alpha_{t,q}(V) = \#\{x \in \mathbb{N}^v \cap (a_0 + L) : |x| = q\},$$
$$\beta_q(V) = \#\{x \in \mathbb{N}^v : |x| = q\},$$

*and $[\mathbb{Z}^v : L] = \#G$. Here $|x| := \sum_i x_i$.*

*Proof:* The representation $V$ splits into one-dimensional representations $V = \bigoplus_i V_i$. Let $x_i$ be a basis vector for $V_i$, then $\operatorname{Sym} V$ can be identified with the polynomial ring $k[x_1, \ldots, x_v]$. It is clear that $\beta_q(V) = \dim_k \operatorname{Sym}^q V$ equals the number of monomials of degree $q$, so the claimed formula holds.

$G$ acts diagonally on $\bigoplus_i V_i$, so we can identify the image of $G$ under $\rho$ with a subgroup $G_1 \subseteq (k^*)^v$ which is isomorphic to $G$, since $(V, \rho)$ is faithful. The action of $G_1$ is now given by

$$g \cdot \underline{x}^{\underline{a}} = g^{\underline{a}} \underline{x}^{\underline{a}},$$

where $\underline{x}^{\underline{a}} = \prod_i x_i^{a_i}$ is a monomial and $g^{\underline{a}} := \prod_i g_i^{a_i}$, $g = (g_1, \ldots, g_\nu)$.

For each point $\underline{a} \in \mathbb{Z}^\nu$, we have a one-dimensional representation $\rho_{\underline{a}} : G_1 \to k^*$, $\rho_{\underline{a}}(g) = g^{\underline{a}}$. So, we obtain a map

$$\ell : \mathbb{Z}^\nu \to \mathrm{Hom}_{\mathbb{Z}}(G_1, k^*)$$
$$\underline{a} \mapsto \rho_{\underline{a}}.$$

Two points $\underline{a}$ and $\underline{b}$ are mapped to the same representation if and only if $\underline{a} - \underline{b} \in \ker \ell$. Let $L$ be the kernel of $\ell$. We claim that $[\mathbb{Z}^\nu : L] = \#G_1$. Let $\mu_n(k) \subseteq k^*$ be the group of $n$-th roots of unity. The order of every element of $G_1$ is divisible by $n$, so $G_1$ is in fact a subgroup of $\mu_n(k)^\nu$. Further, every $\rho_{\underline{a}}$ takes values in $\mu_n(k)$. Note that $\mu_n(k)$ is (not canonically) isomorphic to $\mathbb{Z}/n$ and let $G_2 \subset (\mathbb{Z}/n)^\nu$ be the image of $G_1$ under such an isomorphism. Then the map $\ell$ corresponds to the map

$$\bar{\ell} : \mathbb{Z}^\nu \to \mathrm{Hom}_{\mathbb{Z}}(G_2, \mathbb{Z}/n),$$
$$\underline{a} \mapsto \Big( g \mapsto \sum_i a_i g_i \Big),$$

and the claim follows from the preceding Lemma 4.5.3.

Now, since $[\mathbb{Z}^\nu : L] = \#G_1$ the lattice $L$ divides $\mathbb{Z}^\nu$ into $\#G_1$ cosets. On the other hand, there are exactly $\#G_1$ irreducible representations, so there exists $a_0 \in \mathbb{Z}^\nu$ such that $\ell(a_0)$ is the representation $V'$ which completes the proof. $\square$

Now we are ready to prove Theorem 4.5.2.

*Proof of Theorem 4.5.2.* Let $\nu := \dim_k V$ and let $L \subseteq \mathbb{Z}^\nu$ and $a_0 \in \mathbb{Z}^\nu$ be as in Lemma 4.5.4. Further, let $\Delta \subset \mathbb{R}^\nu$ denote the simplex spanned by the origin and the unit vectors. By Lemma 4.5.4, $\sum_{q=0}^N \beta_q(V)$ equals the number of lattice points in $N\Delta = \{Np : p \in \Delta\}$, and $\sum_{q=0}^N \alpha_{t,q}(V)$ equals $\#(N\Delta \cap (a_0 + L))$.

By choosing a basis for $L$ we find a $\nu \times \nu$ matrix $A$ of full rank, such that $L = A\mathbb{Z}^\nu$. Note that by the elementary divisor theorem, $A$ can be diagonalized and we see that $\det A = [\mathbb{Z}^\nu : L]$. Now we compute:

$$\lim_{N \to \infty} \frac{\sum_{q=0}^N \alpha_{t,q}(V)}{\sum_{q=0}^N \beta_q(V)} = \lim_{N \to \infty} \frac{\#(N\Delta \cap (a_0 + L))}{\#(N\Delta \cap \mathbb{Z}^\nu)}$$

$$= \lim_{N \to \infty} \frac{1/N^\nu \#\big(NA^{-1}\Delta \cap (A^{-1}a_0 + \mathbb{Z}^\nu)\big)}{1/N^\nu \#(N\Delta \cap \mathbb{Z}^\nu)}$$

$$= \frac{\lim_{N \to \infty} 1/N^\nu \#\big(NA^{-1}\Delta \cap (A^{-1}a_0 + \mathbb{Z}^\nu)\big)}{\lim_{N \to \infty} 1/N^\nu \#(N\Delta \cap \mathbb{Z}^\nu)}$$

$$= \frac{\mathrm{Vol}(A^{-1}\Delta)}{\mathrm{Vol}(\Delta)} = \frac{\det(A^{-1})\mathrm{Vol}(\Delta)}{\mathrm{Vol}(\Delta)} = \frac{1}{\det A}$$

$$= \frac{1}{[\mathbb{Z}^\nu : L]} = \frac{1}{\#G}. \qquad\qquad \square$$

*4. Symmetric signature of quotient singularities*

Now we come back to the cyclic group $\frac{1}{n}(1, a)$.

**Lemma 4.5.5.** *Let $G$ be the cyclic group $\frac{1}{n}(1, a)$ over the algebraically closed field $k$ of characteristic coprime with $n$, and let $R = k[\![u, v]\!]^G$ be the invariant ring. If $M$ is the MCM $R$-module $\mathrm{Syz}_R^2(k)$, then the corresponding $k$-representation via Auslander $V = \mathscr{A}'(M)$ is a faithful representation.*

*Proof:* Since $G$ has order $n$, it is enough to show that there exists an element of order $n$ in $V$.

Let $\mathfrak{m}_R$ be the maximal ideal of $R$, and let $\{p_1, \ldots, p_\mu\}$ be a minimal system of generators for it. In particular we choose the $p_i$'s among the monomials $u^i v^j$ such that $i + aj \cong 0$ mod $n$, and we may assume without loss of generality that $p_1 = u^n$ and $p_2 = u^{n-a} v$.

We consider the $S * G$-module $\mathrm{Syz}_S^1(p_1, \ldots, p_\mu)$, and we recall that $\mathrm{Syz}_S^1(p_1, \ldots, p_\mu)^G = \mathrm{Syz}_R^2(k)$. We want to understand the action of $G$ on this module. We consider the element

$$s = (v, u^a, 0, \ldots, 0) \in \mathrm{Syz}_S^1(p_1, \ldots, p_\mu) \subseteq S^\mu$$

Since $(p_1, \ldots, p_\mu)$ is a monomial ideal in $S$, it follows from [MS05, Proposition 3.1] that $s$ is a minimal generator for $\mathrm{Syz}_S^1(p_1, \ldots, p_\mu)$, so it corresponds to a non-zero element of $V$. Let $g = \begin{pmatrix} \xi & 0 \\ 0 & \xi^a \end{pmatrix}$ be a generator of $G$. Then the action of $g$ on $s$ is given by $g \cdot s = \xi^a s$. Since $a$ is coprime with $n$, this shows that $s$ has order $n$, hence the representation $V$ is faithful. $\square$

**Theorem 4.5.6.** *Let $G$ be the group $\frac{1}{n}(1, a)$ over an algebraically closed field $k$. Assume that $(\mathrm{char}\, k, n) = 1$ and let $R$ be the cyclic quotient singularity $k[\![u, v]\!]^G$. Then for any indecomposable MCM $R$-module $M_t$ we have*

$$s_\sigma(R, M_t) = \frac{1}{|G|}.$$

*Proof:* Let $M = \mathrm{Syz}_R^2(k)$, let $\alpha_{t,q}(M)$ be the multiplicity of $M_t$ in $\left(\mathrm{Sym}_R^q(M)\right)^{**}$, and $\beta_q(M) = \mathrm{rank}_R \left(\mathrm{Sym}_R^q(M)\right)^{**}$. If $V = \mathscr{A}'(M)$ and $V_t = \mathscr{A}'(M_t)$, then $V_t$ is an irreducible representation by Theorem 3.5.10 and $V$ is a faithful representation by Lemma 4.5.5.

From Theorem 4.1.8 we have

$$\mathscr{A}'\left(\mathrm{Sym}_R^q(M)^{**}\right) = \mathrm{Sym}^q(V).$$

It follows that $\beta_q(M) = \dim_k \mathrm{Sym}^q(V)$ and $\alpha_{t,q}(M)$ equals the multiplicity of the representation $V_t$ in $\mathrm{Sym}^q(V)$. Thus the application of Theorem 4.5.2 concludes the proof. $\square$

If we choose $M_t = R$ in the previous theorem, we immediately get the following corollary.

**Corollary 4.5.7.** *Let $G = \frac{1}{n}(1, a)$ over an algebraically closed field $k$. Assume that $(\text{char } k, n) = 1$ and let $R$ be the cyclic quotient singularity $k[\![u, v]\!]^G$. Then the symmetric signature of $R$ is*

$$s_\sigma(R) = \frac{1}{|G|}$$

# 5. Symmetric signature of cones over elliptic curves

The main goal of this chapter is to investigate the symmetric signature of the coordinate ring of a plane elliptic curve over an algebraically closed field. Our intention is motivated by the following observation.

**Remark 5.0.8.** We consider $R = k[x, y, z]/(f)$, with $k$ an algebraically closed field of positive characteristic, and $f$ a homogeneous non-singular polynomial of degree 3. In other words, $R$ is the coordinate ring of a plane elliptic curve over $k$. The ring $R$ is not strongly F-regular, and in particular its F-signature $s(R)$ is 0 by the result of Aberbach and Leuschke (Theorem 2.1.10). To see this, assume for simplicity that $f$ is in Weierstrass normal form, that is $f = y^2 z - x^3 - axz^2 - bz^3$ for some $a, b \in k$. This is always possible if the characteristic of $k$ is different from 2 and 3. Then the ideal $I = (x, z)$ is not tightly closed in $R$, since $y^2 \notin I$, but $y^2 \in I^*$. Therefore $R$ is not strongly F-regular and $s(R) = 0$.

Inspired by the previous remark, we would like to prove that for such rings $R$ also the symmetric signature $s_\sigma(R)$ and the differential symmetric signature $s_{d\sigma}(R)$ are 0. The methods we use in this situation are different from those of previous chapters, and are of geometric nature. We will use the correspondence between graded MCM $R$-modules and vector bundles over the smooth projective curve $Y = \operatorname{Proj} R$ to translate the problem into geometric language. We will take advantage of this, and use Atiyah's classification of vector bundles over an elliptic curve to obtain results, which can then be pulled back to the algebraic setting. Since we will use also the multiplicative structure of vector bundles, and in particular Theorem 5.1.22 which holds in characteristic zero, we will focus mainly on elliptic curves over the field of the complex numbers.

We prove that the differential symmetric signature of the coordinate ring $R$ of a plane elliptic curve is 0 (Theorem 5.2.1), but unfortunately, at the time of this writing we are not able to compute $s_\sigma(R)$. However, we present two possible strategies (Sections 5.2.1 and 5.2.2) which allow us to obtain partial results and give an upper bound for the symmetric signature, namely $s_\sigma(R) \leq \frac{1}{2}$. Hopefully these strategies can be completed to prove that $s_\sigma(R) = 0$.

## 5.1. Vector bundles over curves

We briefly recall some definitions and results concerning vector bundles over curves. We assume that the reader is familiar with the concepts of varieties and sheaves, and we will use the notations of Hartshorne's book [Har77]. By *variety* we mean a complete separated integral scheme of finite type over a field, and by *curve* a variety of dimension one. All sheaves are assumed to be coherent. We will skip the proofs of most results, and we refer to the books of Le Potier [LeP97] and Mukai [Muk03, Chapter 10], or to the notes of Teixidor i Bigas [Tei].

Let $Y$ be a curve over an algebraically closed field $k$. A *vector bundle of rank $r$* over $Y$ is a couple $(E, \pi)$, where $E$ is a variety and $\pi$ a morphism $\pi : E \to Y$ called *projection*, such that there exists an open covering $(U_i)_{i \in I}$ of $Y$ and isomorphisms

$$\varphi_i : \pi^{-1}(U_i) \to U_i \times \mathbb{A}_k^r$$

such that on the intersection $U_i \cap U_j$ the composition

$$\varphi_i \circ \varphi_j^{-1} : (U_i \cap U_j) \times \mathbb{A}_k^r \to (U_i \cap U_j) \times \mathbb{A}_k^r$$

is given by $\varphi_i \circ \varphi_j^{-1}(x, v) = (x, g_{i,j}(x)v)$, where $g_{i,j}(x) \in \mathrm{GL}(r, k)$. The maps $\varphi_i$ are called *trivializations* or *charts* and the functions $g_{i,j} : Y \to \mathrm{GL}(r, k)$ are called *transition functions*. The vector bundle is called *algebraic* if the transition functions are morphisms. All vector bundles we deal with are assumed to be algebraic, and sometimes we will refer to them simply as *bundles*. A vector bundle of rank one is also called *line bundle*.

A *morphism* of vector bundles over $Y$ is a morphism of varieties that commutes with the projections to $Y$ and restricts to a linear map on each fiber. The set of all morphisms between two vector bundles $E$ and $F$ is denoted by $\mathrm{Hom}(E, F)$. The category whose objects are vector bundles over $Y$, and whose maps are morphisms between them is denoted by $\mathrm{VB}(Y)$. If $Y$ is a smooth projective connected curve, then $\mathrm{VB}(Y)$ is a Krull-Schmidt category (cf. [Ati56, Theorem 3]).

A *subbundle* of a vector bundle is a closed subvariety which is itself a bundle and such that the inclusion is a morphism. The usual operations on vector spaces like direct sum, tensor product, symmetric and wedge powers extend naturally to vector bundles.

A *section* of a vector bundle $E$ over an open set $U \subseteq Y$ is a map $s : U \to E$ such that $\pi \circ s = \mathrm{Id}_U$. A *global section* is a section over $U = Y$. For any vector bundle $E$ we can define its sheaf of sections

$$U \mapsto \Gamma(U, E) := \{s : U \to E \text{ such that } \pi \circ s = \mathrm{Id}_U\},$$

which is a locally free sheaf.

**Theorem 5.1.1.** *Let $Y$ be a curve and let $r$ be a positive integer. The correspondence $E \mapsto \Gamma(-, E)$ induces an equivalence of categories between the category of vector bundles of rank $r$ over $Y$ and the category of locally free sheaves of rank $r$ on $Y$.*

*Proof:* See [Har77, Ex. II 5.18]. □

Having this identification in mind, we will use indifferently the words vector bundle and locally-free sheaf. For example, we will use the notation for the structure sheaf $\mathscr{O}_Y$ of $Y$ also to denote the trivial bundle of rank one corresponding to it. Nevertheless, we point out that the notions of subbundle and (locally-free) subsheaf are different. For example, the structure sheaf $\mathscr{O}_Y$ has in general many subsheaves, but no non-trivial subbundle.

In the equivalence of Theorem 5.1.1 line bundles correspond to invertible sheaves. If $Y$ is a smooth projective curve, these correspond also to (Cartier) divisors via the map which sends every divisor $D$ to the sheaf $\mathscr{O}_Y(D)$. Finally, we mention the fact that the usual constructions on vector bundles such as tensor, symmetric and wedge products correspond to the same constructions in the category of locally free sheaves.

**Lemma 5.1.2.** *Let $0 \to \mathscr{E} \to \mathscr{F} \to \mathscr{G} \to 0$ be a short exact sequence of coherent sheaves over $Y$. If $\mathscr{F}$ and $\mathscr{G}$ are locally free, then $\mathscr{E}$ is also locally free.*

*Proof:* Let $P$ be a point of $Y$ and consider the corresponding short exact sequence $0 \to E \to F \to G \to 0$ of finitely generated modules over the stalk $\mathscr{O}_P$. Since $\mathscr{F}$ and $\mathscr{G}$ are locally free, the modules $F$ and $G$ are free. So the sequence splits, and $E$ is a direct summand of a free module $F$, hence it is projective. Since $\mathscr{O}_P$ is a Noetherian local ring, $E$ is also free by Theorem 1.4.4. □

The *dual* of a vector bundle $E$ is the bundle $E^\vee := \mathrm{Hom}(E, \mathscr{O}_Y)$. From ordinary properties of vector spaces one has that $E^{\vee\vee}$ is canonically isomorphic to $E$. If $E$ has rank $r$, then the *determinant* of $E$ is the line bundle

$$\det E := \bigwedge^r E.$$

**Remark 5.1.3.** Since every vector bundle $E$ is canonically isomorphic to its double dual $E^{\vee\vee}$, reflexive symmetric powers $\mathrm{Sym}^q(-)^{\vee\vee}$ coincide with ordinary symmetric powers $\mathrm{Sym}^q(-)$ in the category $\mathrm{VB}(Y)$.

Apart from direct sums, tensor, symmetric and wedge products, one of the most useful ways to construct a vector bundle from other bundles is via extensions. We say that a vector bundle $E$ is an *extension* of the bundle $E_2$ by the bundle $E_1$ if there is an exact sequence

$$0 \to E_1 \to E \to E_2 \to 0. \tag{5.1}$$

If the previous sequence splits, then the extension is said to be *trivial* and $E$ is isomorphic to the direct sum $E_1 \oplus E_2$. Two extensions $E$ and $E'$ coming from sequences of the form (5.1) are called *equivalent* if there exists an isomorphism $\varphi : E \to E'$ such that the following diagram commutes

$$\begin{array}{ccccccccc}
0 & \longrightarrow & E_1 & \longrightarrow & E & \longrightarrow & E_2 & \longrightarrow & 0 \\
& & \downarrow{\scriptstyle\mathrm{id}} & & \downarrow{\scriptstyle\varphi} & & \downarrow{\scriptstyle\mathrm{id}} & & \\
0 & \longrightarrow & E_1 & \longrightarrow & E' & \longrightarrow & E_2 & \longrightarrow & 0.
\end{array}$$

**Theorem 5.1.4.** *Let $E_1$ and $E_2$ be two vector bundles, then the extensions of $E_2$ by $E_1$ modulo the equivalence defined above correspond bijectively to the elements of $\mathrm{Ext}^1(E_2, E_1)$.*

*Proof:* See [Rot09, Theorem 7.30]. □

From now on, we assume that $Y$ is a smooth projective curve. We want to define the *degree* of a vector bundle $E$. If $E$ is a line bundle, then there exists a unique (up to linear equivalence) divisor $D$ which corresponds to $E$, then we define $\deg E := \deg D$. If the rank of $E$ is greater than 1, we define $\deg E := \deg(\det E)$, the degree of the determinant.

**Lemma 5.1.5.** *Let $Y$ be a smooth projective curve. Then the degree map $\deg: \mathrm{VB}(Y) \to \mathbb{Z}$ satisfies the following properties.*

1. *If $0 \to E' \to E \to E'' \to 0$ is a short exact sequence of vector bundles over $Y$ then $\deg E = \deg E' + \deg E''$.*

2. *If $E_1$ and $E_2$ are two vector bundles of ranks $r_1$ and $r_2$ then*

$$\deg(E_1 \otimes E_2) = r_1 \deg(E_2) + r_2 \deg(E_1).$$

*Proof:* See [Har77, Ex. II 6.12] □

**Theorem 5.1.6** (Riemann-Roch for vector bundles). *Let $E$ be a vector bundle of rank $r$ over a smooth projective curve $Y$, then*

$$\deg(E) = \chi(E) - r\chi(\mathscr{O}_Y),$$

*where $\chi(-)$ denotes the Euler-Poincaré characteristic, i.e. $\chi(E) := h^0(E) - h^1(E)$.*

**Definition 5.1.7.** *The* slope *of a vector bundle $E$ on a curve $Y$ is defined as $\mu(E) := \frac{\deg E}{\mathrm{rank}\,E}$. A bundle $E$ is said to be* semistable *(resp.* stable*) if for every proper subbundle $F \subset E$ one has $\mu(F) \leq \mu(E)$ (resp. $\mu(F) < \mu(E)$).*

**Proposition 5.1.8.** *Let $E$ be a vector bundle on $Y$ of slope $\mu$. Then the following facts are equivalent.*

1. *$E$ is semistable.*

2. *For every non-zero coherent sheaf $\mathscr{F} \subseteq E$ we have $\mu(\mathscr{F}) \leq \mu$.*

3. *For every non-zero quotient coherent sheaf $\mathscr{F}'$ of $E$ we have $\mu(\mathscr{F}') \geq \mu$.*

From Lemma 5.1.5 we immediately obtain the following lemma.

**Lemma 5.1.9.** *Let $Y$ be a smooth projective curve. Then the following facts hold.*

1. *If $0 \to E' \to E \to E'' \to 0$ is a short exact sequence of vector bundles over $Y$ of ranks $r'$, $r$ and $r''$ respectively, then*

$$\mu(E) = \frac{\mu(E')r' + \mu(E'')r''}{r}.$$

2. *If $E_1$ and $E_2$ are two vector bundles of ranks $r_1$ and $r_2$ then*

$$\mu(E_1 \otimes E_2) = \mu(E_1) + \mu(E_2).$$

We apply the definitions of simple and indecomposable objects of an abelian category from Chapter 1 to the category of vector bundles VB($Y$). A vector bundle $E$ is *simple* if its only proper subbundle is the zero bundle, or equivalently $\operatorname{Hom}(E, E) = k$. $E$ is *indecomposable* if it cannot be written as direct sum of two non-zero subbundles. We have the following implications:

$$\text{indecomposable} \iff \text{simple} \iff \text{stable} \implies \text{semistable.}$$

The fact that a simple object is indecomposable holds in every abelian category and the implication stable $\implies$ semistable follows directly from the definition. For the fact that stable vector bundles are simple see [Muk03, Corollary 10.25]. We point out that the reverse arrows do not hold in general.

**Proposition 5.1.10.** *Let $Y$ be a smooth projective curve, let $E$ be a vector bundle and $\mathscr{L}$ a line bundle over $Y$. Then the following facts hold.*

*1)* $\operatorname{Sym}^q(E \otimes \mathscr{L}) \cong \operatorname{Sym}^q(E) \otimes \mathscr{L}^{\otimes q}$.

*2)* $\mu(\operatorname{Sym}^q(E)) = q \cdot \mu(E)$.

*3) If the field $k$ has characteristic $0$ and $E$ is semistable, then $\operatorname{Sym}^q(E)$ is also semistable.*

*Proof:* Property *1)* is a standard fact of symmetric powers, and for *2)* and *3)* see [LeP97, Theorem 10.2.1]. □

## 5.1.1. Syzygies of vector bundles

Let $R$ be a normal standard graded domain of dimension 2 over a field $k$, that is $R_0 = k$ and $R$ is generated by finitely many elements of degree 1. The normal assumption on $R$ implies that $R_{\mathfrak{p}}$ is a regular ring for every prime ideal $\mathfrak{p} \neq R_+$. It follows that the projective variety $Y = \operatorname{Proj} R$ is a smooth projective curve over $k$. For every graded module $M$ we denote by $\widetilde{M}$ the corresponding coherent sheaf on $Y$. The sheaves $\widetilde{R(n)}$ are invertible (cf. [Har77, Proposition II.5.12]) and denoted by $\mathscr{O}_Y(n)$. Moreover, every MCM graded $R$-module $M$ is locally free on the punctured spectrum, so the associated coherent sheaf $\widetilde{M}$ is in fact a vector bundle over $Y$. If we denote by $\operatorname{MCM}_{\mathbb{Z}}(R)$ the category of finitely generated graded MCM $R$-modules, then we have a functor

$$\widetilde{\phantom{m}} : \operatorname{MCM}_{\mathbb{Z}}(R) \to \operatorname{VB}(Y),$$

whose properties we collect in the following proposition.

**Proposition 5.1.11.** *Let $R$ be a normal standard graded domain of dimension $2$ over a field $k$, let $Y = \operatorname{Proj} R$, and let $M$ and $N$ be finitely generated graded MCM $R$-modules. Then the following facts hold.*

*1)* $\widetilde{M \oplus N} \cong \widetilde{M} \oplus \widetilde{N}$;

*2)* $\widetilde{M \boxtimes_R N} \cong \widetilde{M} \otimes_{\mathscr{O}_Y} \widetilde{N}$;

*3)* $\operatorname{Sym}_R^q(M)^{**} \cong \operatorname{Sym}_Y^q(\widetilde{M})$.

*Proof:* These statements follow from the corresponding statements on the punctured spectrum $U = D(R_+)$. Then, property *1)* is straightforward, and *2)* is a consequence of Corollary 1.5.25. We check *3)*.

Since symmetric powers commute with sheafification and from Remark 1.5.24, we have

$$\operatorname{Sym}_U^q(\widetilde{M}) \cong \widetilde{\operatorname{Sym}_R^q(M)}\Big|_U = \widetilde{\operatorname{Sym}_R^q(M)^{**}}\Big|_U$$

as required. □

Now we give an appropriate definition of free rank for the category of vector bundles over $Y$.

**Definition 5.1.12.** Let $E$ be a vector bundle over $Y$ with a fixed very ample invertible sheaf $\mathscr{O}_Y(1)$. We define the *free rank* of $E$ as

$$\operatorname{frk}_{\mathscr{O}_Y(1)}(E) := \max\Big\{n : \exists \text{ a split surjection } \varphi : E \twoheadrightarrow F, \text{ with } F = \bigoplus_{i=1}^n \mathscr{O}_Y(-d_i)$$

$$\text{a splitting vector bundle of rank } n\Big\}.$$

Since the category $\operatorname{VB}(Y)$ has the KRS property, to compute the free rank of a bundle $E$, one should count how many copies of twisted structure sheaves $\mathscr{O}_Y(-d_i)$ appear in the decomposition of $E$ into indecomposable bundles. For this reason, the free rank of vector bundles $\operatorname{frk}_{\mathscr{O}_Y(1)}$ agrees with the graded free rank $\operatorname{frk}_R^{\mathrm{gr}}$ of Definition 1.1.19, as explained in the following corollary.

**Corollary 5.1.13.** *Let $R$ be a normal standard graded domain of dimension $2$ over a field $k$, let $Y = \operatorname{Proj} R$, and let $M$ be a finitely generated graded MCM $R$-module. Then the following facts hold.*

*1.* $\operatorname{rank}_R M = \operatorname{rank}_{\mathscr{O}_Y} \widetilde{M}$.

*2.* $\operatorname{frk}_R^{\mathrm{gr}} M = \operatorname{frk}_{\mathscr{O}_Y(1)} \widetilde{M}$.

**Remark 5.1.14.** Let $I$ be an $R_+$-primary homogeneous ideal with homogeneous generators $f_1, \ldots, f_n$ of degrees $d_i$. The following presentation of $R/I$

$$\bigoplus_{i=1}^n R(-d_i) \xrightarrow{f_1, \ldots, f_n} R \to R/I \to 0$$

induces a short exact sequence of sheaves on $Y$

$$0 \to \mathscr{S} \to \bigoplus_{i=1}^{n} \mathscr{O}_Y(-d_i) \to \mathscr{O}_Y \to 0,$$

where $R/I$ corresponds to 0 because $I$ is $R_+$-primary. Since $\mathscr{O}_Y$ and $\bigoplus_{i=1}^{n} \mathscr{O}_Y(-d_i)$ are locally free, the sheaf $\mathscr{S}$ is also locally free by Lemma 5.1.2. The sheaf $\mathscr{S}$ is the kernel of the sheaf morphism $f_1, \ldots, f_n$ and is denoted by $\mathrm{Syz}(f_1, \ldots, f_n)$ and called *syzygy bundle*. In fact, it is nothing but the vector bundle corresponding to the MCM graded $R$-module $\mathrm{Syz}_R^2(R/I) = \mathrm{Syz}_R^1(f_1, \ldots, f_n)$.

Syzygies of the irrelevant maximal ideal $R_+$ play a special role, they correspond to the restriction of the cotangent bundle $\Omega_{\mathbb{P}^n}$ of $\mathbb{P}^n$ to the curve. In fact, it follows from the Euler sequence on $\mathbb{P}^n$ [Har77, Theorem 8.13]

$$0 \to \Omega_{\mathbb{P}^n} \to \bigoplus_{i=0}^{n} \mathscr{O}_{\mathbb{P}^n}(-1) \to \mathscr{O}_{\mathbb{P}^n} \to 0,$$

that $\Omega_{\mathbb{P}^n}$ is a syzygy bundle. So if $x_0, \ldots, x_n$ is a system of generators for $\mathscr{O}_{\mathbb{P}^n}(1)$, we have the isomorphism

$$\mathrm{Syz}(x_0, \ldots, x_n) \cong \Omega_{\mathbb{P}^n}.$$

If we restrict these bundles to the curve $Y$, we obtain isomorphic bundles. In fact, we have just proved the following proposition.

**Proposition 5.1.15.** *Let $Y$ be a smooth projective curve with an embedding in $\mathbb{P}^n$ given by $\varphi : Y \to \mathbb{P}^n$. Let $x_0, \ldots, x_n$ be a system of generators of $\mathscr{O}_{\mathbb{P}^n}(1)$ and denote by $\Omega_{\mathbb{P}^n}$ the cotangent bundle of $\mathbb{P}^n$ and by $\Omega_{\mathbb{P}^n}|_Y$ its restriction to $Y$. Then, we have the following isomorphism of vector bundles on $Y$*

$$\mathrm{Syz}(x_0, \ldots, x_n)|_Y \cong \Omega_{\mathbb{P}^n}|_Y.$$

### 5.1.2. Vector bundles over elliptic curves

Since the category of vector bundles over a smooth projective curve has the KRS property, the problem of finding its indecomposable objects arises naturally. This may be a difficult task in general, however for a curve of genus 0 the situation is quite simple.

**Theorem 5.1.16** (Grothendieck)**.** *The indecomposable objects of* $\mathrm{VB}(\mathbb{P}^1)$ *are those of the form* $\mathscr{O}_{\mathbb{P}^1}(a)$ *for some integer $a$. Then, every vector bundle on $\mathbb{P}^1$ can be uniquely decomposed as*

$$\mathscr{O}(a_1) \oplus \cdots \oplus \mathscr{O}(a_r)$$

*for some integers $a_1 \geq \cdots \geq a_r$.*

*5. Symmetric signature of cones over elliptic curves*

The next interesting case is a curve of genus 1, that is an elliptic curve. The situation here is more complicated, but still it is possible to give a complete description of the indecomposable objects of VB($Y$) for an elliptic curve $Y$ over an algebraically closed field. This was done by Atiyah in his wonderful paper [Ati57]. Moreover, he was able to understand the multiplicative structure of the monoid VB($Y$) with the tensor product operation if the field has characteristic zero.

We will illustrate Atiyah' s results, but we will omit proofs. The interested reader may consult the original paper [Ati57] or the notes of Teixidor i Bigas [Tei].

Let $Y$ be an elliptic curve over an algebraically closed field $k$ with structure sheaf $\mathscr{O}$. We fix a line bundle $A$ on $Y$ of degree 1. This is equivalent to fix a point $P_0$ as neutral element for the group operation on $Y$, in fact $A = \mathscr{O}(P_0)$. For every positive integer $r$ and every $d \in \mathbb{Z}$, we denote by $\mathscr{E}(r, d)$ the set of isomorphism classes of indecomposable vector bundles on $Y$ of rank $r$ and degree $d$. First, we are going to describe $\mathscr{E}(r, 0)$.

**Theorem 5.1.17** (Atiyah). *For every positive integer $r$ there exists a unique indecomposable vector bundle $F_r$ in $\mathscr{E}(r, 0)$ with $\Gamma(Y, F_r) \neq 0$. Moreover, the following facts hold.*

1. *$F_r$ has only one section up to scalar multiplication, i.e. $\Gamma(Y, F_r) \cong k$.*
2. *$F_r \cong F_r^\vee$.*
3. *There is a non-split exact sequence*

$$0 \to \mathscr{O} \to F_r \to F_{r-1} \to 0.$$

4. *$\mathrm{Sym}^q(F_2) \cong F_{q+1}$.*
5. *If $E$ is an indecomposable vector bundle of rank $r$ and degree $0$ then $E \cong L \otimes F_r$, where $L$ is a line bundle of degree zero, unique up to isomorphism and such that $L \cong \det E$.*

*Proof:* These are Theorem 5 (and subsequent corollaries), and Theorem 9 of [Ati57]. □

We call the bundle $F_r$ *Atiyah bundle of rank $r$*. Clearly we have that $F_1 = \mathscr{O}$, the structure sheaf is the unique line bundle of degree 0 with non-zero sections. The bundle $F_2$ is given by the unique non-trivial extension

$$0 \to \mathscr{O} \to F_2 \to \mathscr{O} \to 0.$$

In other words, $F_2$ is the non-zero element of $\mathrm{Ext}^1(\mathscr{O}, \mathscr{O})$.

For a plane elliptic curve we can give another description of $F_2$ using syzygy bundles.

**Proposition 5.1.18.** *Let $Y$ be a plane elliptic curve over a field $k$ of characteristic $\geq 5$ with coordinate ring $k[x, y, z]/(f)$, where $f$ is a homogeneous polynomial of degree 3. Then we have*

$$\mathrm{Syz}\left(\frac{\partial f}{\partial x}, \frac{\partial f}{\partial y}, \frac{\partial f}{\partial z}\right)(3) \cong F_2.$$

*Proof:* It is enough to show that $\mathrm{Syz}\left(\frac{\partial f}{\partial x}, \frac{\partial f}{\partial y}, \frac{\partial f}{\partial z}\right)(3)$ corresponds to a non-zero element of $\mathrm{Ext}^1(\mathcal{O}, \mathcal{O})$ with global sections, then the uniqueness of $F_2$ will imply the desired isomorphism.

We consider the following injective map of sheaves, $\varphi : \mathcal{O} \to \widetilde{\mathcal{O}(+1)}^{\oplus 3}$, given by $\varphi(1) = (x, y, z)$. From the Euler formula

$$\frac{\partial f}{\partial x} x + \frac{\partial f}{\partial y} y + \frac{\partial f}{\partial z} z = 3f,$$

which vanishes on $Y$, we obtain that the image of $\varphi$ is contained in the syzygy bundle $\mathrm{Syz}\left(\frac{\partial f}{\partial x}, \frac{\partial f}{\partial y}, \frac{\partial f}{\partial z}\right)(3)$. Actually, since this image does not vanish anywhere, it defines a subbundle and hence also a quotient bundle, which is by rank and degree reasons the structure sheaf. In other words, we have the following short exact sequence of vector bundles

$$0 \to \mathcal{O} \xrightarrow{\varphi} \mathrm{Syz}\left(\frac{\partial f}{\partial x}, \frac{\partial f}{\partial y}, \frac{\partial f}{\partial z}\right)(3) \to \mathcal{O} \to 0, \tag{5.2}$$

It remains to prove that the sequence (5.2) is non-split, equivalently $(x, y, z)$ is the unique non-zero global section. This can be done easily by explicit computations, for example assuming that $f$ is in Weierstrass normal form. $\square$

**Corollary 5.1.19.** *Let $\widetilde{\Omega_{R/k}}$ be the sheaf version of the cotangent module $\Omega_{R/k}$ of the cone $R$ of $Y$, then we have an isomorphism of vector bundles*

$$\widetilde{\Omega_{R/k}}(-1) \cong F_2.$$

*Proof:* We recall that the cotangent module $\Omega_{R/k}$ is the graded $R$-module

$$\Omega_{R/k} = < \mathrm{d}x, \mathrm{d}y, \mathrm{d}z > / \left(\frac{\partial f}{\partial x} \mathrm{d}x + \frac{\partial f}{\partial y} \mathrm{d}y + \frac{\partial f}{\partial z} \mathrm{d}z\right).$$

In other words, $\Omega_{R/k}$ can be defined by the following short exact sequence of graded $R$-modules

$$0 \to R(-2) \xrightarrow{\psi} R^{\oplus 3} \xrightarrow{\varphi} \Omega_{R/k} \to 0, \tag{5.3}$$

where $\psi(1) = \left(\frac{\partial f}{\partial x}, \frac{\partial f}{\partial y}, \frac{\partial f}{\partial z}\right)$, and $\varphi$ sends the canonical basis to $\mathrm{d}x, \mathrm{d}y, \mathrm{d}z$.

Sequence (5.3) induces a short exact sequence of locally free sheaves on $Y$

$$0 \to \mathcal{O}_Y(-2) \to \mathcal{O}_Y^{\oplus 3} \to \widetilde{\Omega_{R/k}} \to 0.$$

We dualize this sequence and we get

$$0 \to \widetilde{\Omega_{R/k}}^\vee \to \mathcal{O}_Y^{\oplus 3} \to \mathcal{O}_Y(2) \to 0,$$

where the last map is given by $\psi$. We obtain that

$$\widetilde{\Omega_{R/k}}^\vee \cong \mathrm{Syz}\left(\frac{\partial f}{\partial x}, \frac{\partial f}{\partial y}, \frac{\partial f}{\partial z}\right)(2).$$

So by the previous Proposition 5.1.18 we have $\widetilde{\Omega_{R/k}}^\vee(1) \cong F_2$ on $Y$, but the Atiyah bundle $F_2$ is self-dual, so we have $\widetilde{\Omega_{R/k}}(-1) \cong F_2^\vee \cong F_2$. $\quad\square$

Now that the elements of $\mathscr{E}(r,0)$ have been described, we have two fundamental ways to produce the elements of $\mathscr{E}(r,d)$: tensoring with the line bundle $A$ and taking extensions. Given a positive integer $r$ and $d \in \mathbb{Z}$, Atiyah constructs a map

$$\alpha_{r,d} : \mathscr{E}(h,0) \to \mathscr{E}(r,d)$$

where $h = (r,d)$ is the highest common factor of $r$ and $d$. The map $\alpha_{r,d}$ is defined by the following properties:

1. $\alpha_{r,0}$ is the identity;
2. $\alpha_{r,d+r}(E) \cong \alpha_{r,d}(E) \otimes A$ for every $E \in \mathscr{E}(h,0)$;
3. if $0 < d < r$ and $E \in \mathscr{E}(h,0)$ then $\alpha_{r,d}(E)$ is the non-trivial extension

$$0 \to \mathscr{O}^d \to \alpha_{r,d}(E) \to \alpha_{r-d,d}(E) \to 0.$$

**Theorem 5.1.20** (Atiyah). *The map $\alpha_{r,d} : \mathscr{E}(h,0) \to \mathscr{E}(r,d)$ is well-defined and gives a one-one correspondence between $\mathscr{E}(h,0)$ and $\mathscr{E}(r,d)$. Moreover we have $\det \alpha_{r,d}(E) \cong \det E \otimes A^d$.*

*Proof:* See Theorem 7 of [Ati57]. $\quad\square$

For every rank $r$ and degree $d$ we define an indecomposable vector bundle $E_A(r,d)$ as the image of the Atiyah bundle $F_h$ via the map $\alpha_{r,d}$, that is

$$E_A(r,d) := \alpha_{r,d}(F_h).$$

Like Atiyah bundles, also bundles $E_A(r,d)$ have some special properties.

**Corollary 5.1.21.** *Let $h = (r,d) = 1$ and let $E \in \mathscr{E}(r,d)$. Then the following facts hold.*

1. *There exists a line bundle $\mathscr{L}$ of degree $0$ such that $E \cong E_A(r,d) \otimes \mathscr{L}$.*
2. *$E_A(r,d) \otimes \mathscr{L} \cong E_A(r,d)$ if and only if $\mathscr{L}^{\otimes r} \cong \mathscr{O}$.*
3. *$E_A(r,d)^\vee \cong E_A(r,-d)$.*
4. *$\det E_A(r,d) \cong A^{\otimes d}$.*

Now we look at the multiplicative structure of vector bundles over $Y$ with tensor product.

**Theorem 5.1.22** (Atiyah). *If the characteristic of the field $k$ is $0$, then the following relations hold.*

1) *If $r \geq s$, then $F_r \otimes F_s \cong F_{r-s+1} \oplus F_{r-s+3} \oplus \cdots \oplus F_{r+s-1}$.*
2) *If $(r,d) = 1$, then $E_A(r,d) \otimes F_h \cong E_A(rh,dh)$.*
3) *If $(r,r') = (r,d) = (r',d') = 1$, then $E_A(r,d) \otimes E_A(r',d') \cong E_A(rr', rd' + r'd)$.*

*4) If $a_1 > a_2$, $(d_1, p) = (d_2, p) = 1$, then*

$$E_A(p^{a_1}, d_1) \otimes E_A(p^{a_2}, d_2) \cong \mathcal{O}^{p_{a_2}} \otimes E_A(p^{a_1}, d_1 + p^{a_2 - a_1} d_2).$$

*5) If $(d_1, p) = (d_2, p) = 1$, then*

$$E_A(p^a, d_1) \otimes E_A(p^a, d_2) \cong \mathcal{O}^{p^b} \otimes \left( \bigoplus L_i \right) \otimes E_A(p^b, d_3),$$

*where $p^{a-b} = (p^a, d_1 + d_2)$, $d_3 = (d_1 + d_2)/p^{a-b}$ and the $L_i$'s are line bundles of degree 0 and order $p^a$ taken modulo the class of elements of order $p^b$.*

*6) If $\mathcal{L}$ is a line bundle of degree 0 and $h = (r, d)$, then $E_A(r, d) \otimes \mathcal{L} \cong E_A(r, d)$ if and only if $\mathcal{L}^{\otimes r/h} \cong \mathcal{O}$ with $h = (r, d)$.*

*Proof:* These are Theorem 10, Theorem 13, and Theorem 14 of [Ati57]. □

Theorem 5.1.22 allows us to compute the tensor product of every two indecomposable vector bundles $E_A(r_1, d_1)$ and $E_A(r_2, d_2)$. First, if $(r_i, d_i) = h_i > 1$ we apply *2)* and we write $E_A(r_i, d_i) \cong F_{h_i} \otimes E_A(\frac{r_i}{h_i}, \frac{d_i}{h_i})$. Then we reduce to the situation where the ranks of the bundles $E_A(-, -)$ are a power of a prime using *3)*. So we can multiply these bundles using *4)* and *5)* if the primes are the same, and *3)* if not. We can also multiply the Atiyah bundles using *1)*. Finally, using *2)* and *3)* we can write the final product as a sum of indecomposable bundles.

## 5.2. Symmetric signature of elliptic curves

We would like to compute the symmetric signature $s_\sigma(R)$ of the coordinate ring $R$ of a plane elliptic curve over an algebraically closed field $k$. So we fix the following setting. Let $f$ be a homogeneous non-singular polynomial $f$ of degree 3 in $k[x, y, z]$, and let $R = k[x, y, z]/(f)$. Then $R$ is a normal standard graded $k$-domain of dimension 2 and the projective curve $Y = \text{Proj } R$ is a plane elliptic curve over $k$.

To compute the symmetric signature of $R$ one should consider reflexive symmetric powers of the second syzygy of the residue field, $\text{Syz}_R^2(k)$. Thanks to Remark 5.1.14 and Proposition 5.1.11 this is equivalent to consider the syzygy bundle $\text{Syz}(x, y, z)$ over $Y$ and taking ordinary symmetric powers in the category $\text{VB}(Y)$. We recall that by Corollary 5.1.13 we have $\text{rank}_R M = \text{rank}_{\mathcal{O}_Y} \widetilde{M}$ and $\text{frk}_R^{\text{gr}} M = \text{frk}_{\mathcal{O}_{Y(1)}} \widetilde{M}$ for every MCM graded $R$-module. Therefore, we have that $s_\sigma(R)$ exists if and only if the limit

$$\lim_{N \to +\infty} \frac{\sum_{q=0}^{N} \text{frk}_{\mathcal{O}_{Y(1)}} \text{Sym}^q(\text{Syz}(x, y, z))}{\sum_{q=0}^{N} \text{rank}_{\mathcal{O}_Y} \text{Sym}^q(\text{Syz}(x, y, z))}$$

exists, and in this case they coincide.

The same reasoning applies to the differential symmetric signature. In this case, we should consider the sheaf associated to the module of Zariski differentials $\Omega_{R/k}^{**}$ on $Y$ and then take ordinary symmetric powers in the category VB($Y$). In other words, we consider the vector bundles

$$\mathrm{Sym}_Y^q\left(\widetilde{\Omega_{R/k}}\right).$$

Notice that we can forget about the double dual inside $\mathrm{Sym}_Y^q(-)$, thanks to Remark 2.2.16 and Proposition 5.1.11. The differential symmetric signature of the coordinate ring $R$ exists if and only if the limit

$$\lim_{N\to+\infty}\frac{\sum_{q=0}^N \mathrm{frk}_{\mathscr{O}_Y(1)}\mathrm{Sym}_Y^q\left(\widetilde{\Omega_{R/k}}\right)}{\sum_{q=0}^N \mathrm{rank}_{\mathscr{O}_Y}\mathrm{Sym}_Y^q\left(\widetilde{\Omega_{R/k}}\right)} \tag{5.4}$$

exists, and in this case they coincide.

We can compute the differential symmetric signature of the cone $R$..

**Theorem 5.2.1.** *Let $Y$ be a plane elliptic curve over an algebraically closed field $k$ of characteristic $\geq 5$ with coordinate ring $R$. Then the differential symmetric signature of $R$ is $s_{d\sigma}(R) = 0$.*

*Proof:* By the previous observations we should compute the limit (5.4). From Corollary 5.1.19, we have $\widetilde{\Omega_{R/k}} \cong F_2 \otimes \mathscr{O}_Y(1)$. Therefore from part *4)* of Theorem 5.1.17 and from Proposition 5.1.10 we obtain

$$\mathrm{Sym}_Y^q\left(\widetilde{\Omega_{R/k}}\right) \cong \mathrm{Sym}_Y^q(F_2) \otimes \mathscr{O}_Y(q) \cong F_{q+1} \otimes \mathscr{O}_Y(q)$$

for all $q \geq 1$. So the module $\mathrm{Sym}_Y^q\left(\widetilde{\Omega_{R/k}}\right)$ is indecomposable, and in particular it has free rank 0 and the claim follows. □

For the rest of the chapter we concentrate on the symmetric signature of the cone $R$, and we work over the field of the complex numbers $\mathbb{C}$. First, we investigate the sygygy bundle $\mathrm{Syz}(x,y,z)$.

**Proposition 5.2.2.** *The vector bundle $\mathrm{Syz}(x,y,z)$ is stable of rank 2 and degree $-9$. Moreover $\det\mathrm{Syz}(x,y,z) = \mathscr{O}_Y(-3)$.*

*Proof:* The syzygy bundle fits into a short exact sequence

$$0 \to \mathrm{Syz}(x,y,z) \to \bigoplus_{i=1}^3 \mathscr{O}_Y(-1) \xrightarrow{x,y,z} \mathscr{O}_Y \to 0. \tag{5.5}$$

So from the additivity of rank and degree, we immediately get that $\mathrm{Syz}(x,y,z)$ has rank 2 and degree

$$\deg\mathrm{Syz}(x,y,z) = \left(\deg\bigoplus_{i=1}^3 \mathscr{O}_Y(-1) - \deg\mathscr{O}_Y\right)\deg Y = -3\deg Y = -9.$$

For the indecomposable property, we have that $\mathrm{Syz}(x,y,z)$ equals the restriction of the cotangent bundle $\Omega_{\mathbb{P}^2}$ of $\mathbb{P}^2$ to $Y$ by Proposition 5.1.15. This bundle is stable, and in particular indecomposable, by the result of Brenner and Hein [BH06, Theorem 1.3]. Finally, taking determinants of sequence (5.5) yields $\det\mathrm{Syz}(x,y,z) = \mathscr{O}_Y(-3)$. □

**Remark 5.2.3.** Since $\mathrm{Syz}(x,y,z)$ is indecomposable of rank 2 one would be tempted to write $\mathrm{Syz}(x,y,z) \cong F_2 \otimes \mathcal{L}$ for some line bundle $\mathcal{L}$ of degree 0 and then apply the formula $\mathrm{Sym}^q(F_2) \cong F_{q+1}$ as in the proof of Theorem 5.2.1. However, such a decomposition does not hold. For example, observe that $\deg \mathrm{Syz}(x,y,z) = -9$ and $\deg F_2 \otimes \mathcal{L}$ is always even by Lemma 5.1.5.

**Remark 5.2.4.** Since $\mathrm{Syz}(x,y,z)$ has rank 2, we have that $\mathrm{rank}_{\mathcal{O}_Y} \mathrm{Sym}^q(\mathrm{Syz}(x,y,z)) = q+1$. The difficult part is to determine the free rank of $\mathrm{Sym}^q(\mathrm{Syz}(x,y,z))$.

**Corollary 5.2.5.** *Let $q$ be odd, then the bundle $\mathrm{Sym}^q(\mathrm{Syz}(x,y,z))$ contains no subbundle of rank one. In particular,* $\mathrm{frk}_{\mathcal{O}_Y(1)} \mathrm{Sym}^q(\mathrm{Syz}(x,y,z)) = 0$.

*Proof:* From Proposition 5.1.10 and Proposition 5.2.2 we know that $\mathrm{Sym}^q(\mathrm{Syz}(x,y,z))$ is semistable of slope $-\frac{9}{2}q$. Let $\mathcal{L}$ be a non-zero subbundle of rank one. Since $\mathrm{Sym}^q(\mathrm{Syz}(x,y,z))$ is semistable, we have $\mu(\mathcal{L}) = -\frac{9}{2}q$, but $\mu(\mathcal{L}) = \deg \mathcal{L}$ must be an integer. Since $q$ is odd we get a contradiction. $\square$

**Corollary 5.2.6.** *If the symmetric signature of $R$ exists, then $s_\sigma(R) \leq \frac{1}{2}$.*

*Proof:* For every $q \in \mathbb{N}$, let $a_q = \mathrm{frk}_{\mathcal{O}_Y(1)} \mathrm{Sym}^q(\mathrm{Syz}(x,y,z))$, and let $b_q = \mathrm{rank}_{\mathcal{O}_Y} \mathrm{Sym}^q(\mathrm{Syz}(x,y,z))$. From the previous Corollary 5.2.5 and Remark 5.2.4, we have $a_q = 0$ for $q$ odd, and $b_q = q+1$ for all $q$.

For every even $q$ we have the following inequalities

$$a_q + a_{q+1} = a_q \leq b_q = \frac{1}{2}(b_q + b_q) \leq \frac{1}{2}(b_q + b_{q+1}).$$

It follows that

$$\frac{\sum_{q=0}^{N} a_q}{\sum_{q=0}^{N} b_q} \leq \frac{\frac{1}{2}\sum_{q=0}^{N} b_q}{\sum_{q=0}^{N} b_q} + \phi(N) = \frac{1}{2} + \phi(N),$$

where $\phi(N) = \frac{b_N}{\sum_{q=0}^{N} b_q}$ if $N$ is even, and 0 otherwise. In both cases $\phi(N) \to 0$ if $N \to +\infty$, so the claim is proved. $\square$

Now we outline two possible strategies to compute the free rank of $\mathrm{Sym}^q(\mathrm{Syz}(x,y,z))$ and understand the decomposition of $\mathrm{Sym}^q(\mathrm{Syz}(x,y,z))$ into irreducible bundles also for even $q$.

### 5.2.1. Strategy 1: $\mathrm{Sym}^q$ as subbundle of $T^q$

Let $E$ be a vector bundle over $Y$ and let $q$ a positive integer. The $q$-th symmetric powers of $E$ can be seen as the invariant subvector space of the $q$-th tensor product $T^q(E) := E \otimes \cdots \otimes E$ under the action of the symmetric group $\mathfrak{S}_q$ on $q$ elements. Thus, we have an inclusion

$$\mathrm{Sym}^q(E) \hookrightarrow T^q(E).$$

Since we are in characteristic zero, the previous inclusion splits. A splitting is given by the Reynolds operator associated to the group action of $\mathfrak{S}_q$. It follows that $\mathrm{Sym}^q(E)$ is a direct summand of $T^q(E)$. In particular, direct summands of $\mathrm{Sym}^q(E)$ are precisely the direct summands of $T^q(E)$ which are invariant under the action of $\mathfrak{S}_q$.

We apply the previous observations to the bundle $E = \mathrm{Syz}(x, y, z)$, and we outline the following strategy.

1) Write down the decomposition of $T^q(E)$ into indecomposable bundles using Theorem 5.1.22.
2) Find which indecomposable summands of $T^q(E)$ are invariant under the action of $\mathfrak{S}_q$.
3) Count how many bundles of the form $\mathscr{O}_Y(a)$ appear in the previous decomposition and compute $\mathrm{frk}_{\mathscr{O}_Y(1)} \mathrm{Sym}^q(E)$.

Part *1)* will be done in the following Proposition 5.2.7, while part *2)* and *3)* are still work in progress at the time of this writing.

**Proposition 5.2.7.** *Let $E = \mathrm{Syz}(x, y, z)$ on $Y$, and let $q > 1$ be an integer. Then, the bundle $T^q(E)$ decomposes as follows.*

1) *If $q$ is even, $T^q(E)$ is a sum of line bundles, each one being of the form*

$$L_i \otimes A^{\otimes - \frac{9}{2}q} \otimes \mathscr{L}^{\otimes q}.$$

*Here $L_i$ is a line bundle of degree $0$ and order $2$, that is $L_i \otimes L_i \cong \mathscr{O}$, and $\mathscr{L}$ is the unique line bundle of degree $0$ such that $\mathrm{Syz}(x, y, z) \cong E_A(2, -9) \otimes \mathscr{L}$.*

2) *If $q$ is odd, $T^q(E)$ is a sum of indecomposable vector bundles of rank two, each one being of the form*

$$L_i \otimes E_A(2, d) \otimes \mathscr{L}^{\otimes q}.$$

*Here $L_i$ and $\mathscr{L}$ are as in 1), and $d = -9 + 2^{\frac{9}{2}(q-1)}$.*

*Proof:* We recall that by Corollary 5.1.21 there exists a unique line bundle $\mathscr{L}$ of degree $0$ such that $\mathrm{Syz}(x, y, z) \cong E_A(2, -9) \otimes \mathscr{L}$. First, we consider the case $q = 2$. From Theorem 5.1.22 we have

$$\mathrm{Syz}(x, y, z) \otimes \mathrm{Syz}(x, y, z) \cong E_A(2, -9) \otimes E_A(2, -9) \otimes \mathscr{L}^{\otimes 2}$$

$$\cong \bigoplus_{i=0}^{3} L_i \otimes E_A(1, -9) \otimes \mathscr{L}^{\otimes 2} = \bigoplus_{i=0}^{3} L_i \otimes A^{\otimes -9} \otimes \mathscr{L}^{\otimes 2},$$

where the bundles $L_0 = \mathscr{O}$, $L_1$, $L_2$ and $L_3$ are the four degree $0$ line bundles of order $2$ corresponding to the four 2-torsion points of $Y$ (cf. [Sil86, Theorem 6.1]).

For $q > 2$ even, we write

$$T^q(E) = (T^2(E))^{\otimes \frac{q}{2}}$$

$$\cong \left( \bigoplus_{i=0}^{3} L_i \otimes A^{\otimes -9} \otimes \mathscr{L}^{\otimes 2} \right)^{\otimes \frac{q}{2}}$$

$$\cong \bigoplus_{i_t \in \{0,\dots,3\}} \left( L_{i_1} \otimes \dots \otimes L_{i_{\frac{q}{2}}} \right) \otimes A^{\otimes -\frac{9}{2} q} \otimes \mathscr{L}^{\otimes q}$$

$$\cong \bigoplus L_i \otimes A^{\otimes -\frac{9}{2} q} \otimes \mathscr{L}^{\otimes q},$$

where the last isomorphism follows from the fact that every product of the form $L_{i_1} \otimes \dots \otimes L_{i_{\frac{q}{2}}}$ is isomorphic to one of the bundles $L_i$. Therefore, we have four different kinds of line bundles appearing as summand in $T^q(E)$, whose rank is equal to $2^q$.

Finally, for $q$ odd we have

$$T^q(E) = T^{q-1}(E) \otimes E$$

$$\cong \bigoplus L_i \otimes A^{\otimes -\frac{9}{2}(q-1)} \otimes \mathscr{L}^{\otimes q-1} \otimes E_A(2,-9) \otimes \mathscr{L}$$

$$\cong \bigoplus L_i \otimes \mathscr{L}^{\otimes q} \otimes E_A\left(1, -\frac{9}{2}(q-1)\right) \otimes E_A(2,-9).$$

Computing the product $E_A\left(1, -\frac{9}{2}(q-1)\right) \otimes E_A(2,-9)$ with Theorem 5.1.22 we obtain the indecomposable bundle $E_A(2,d)$ of degree $d = -9 + 2^{\frac{9}{2}(q-1)}$. $\quad\square$

**Question 5.2.8.** Can we use Proposition 5.2.7 to get an upper bound for $s_\sigma(R)$?

**Example 5.2.9.** We give a better description of the case $q = 2$. From Proposition 5.2.7 we have that

$$T^2(E) \cong (\mathscr{O}_Y \oplus L_1 \oplus L_2 \oplus L_3) \otimes A^{\otimes -9} \otimes \mathscr{L}^{\otimes 2},$$

where $E = \text{Syz}(x, y, z)$. So $T^2(E)$ decompose as the sum of four line bundles. Of those, only one can be of the form $\mathscr{O}_Y(a)$. On the other hand, we know that

$$T^2(E) = \text{Sym}^2(E) \oplus \bigwedge^2(E),$$

and $\bigwedge^2(E) \cong \mathscr{O}_Y(-3)$ by Proposition 5.2.2. It follows that $\text{frk}_{\mathscr{O}_Y(1)} \text{Sym}^2(\text{Syz}(x, y, z)) = 0$.

## 5.2.2. Strategy 2: an upper bound for $\text{frk}_{\mathscr{O}_Y(1)}$

This strategy should allow us to obtain an upper bound for $\text{frk}_{\mathscr{O}_Y(1)} \text{Sym}^q(\text{Syz}(x, y, z))$. Our hope is that this upper bound is actually 0. The main idea is to use the following lemma.

**Lemma 5.2.10.** *Let* $0 \to \mathcal{S} \to \mathcal{F} \xrightarrow{\varphi} \mathcal{L} \to 0$ *be a short exact sequence of locally free sheaves on a smooth curve, with* $\mathcal{L}$ *a line bundle. Then, for every* $q \geq 1$ *we have a short exact sequence*

$$0 \to \mathrm{Sym}^q(\mathcal{S}) \to \mathrm{Sym}^q(\mathcal{F}) \to \mathrm{Sym}^{q-1}(\mathcal{F}) \otimes \mathcal{L} \to 0.$$

*In particular, the last non-zero map is given by*

$$s_1 \circ \cdots \circ s_q \mapsto \sum_{i=1}^{q} s_1 \circ \cdots \circ \hat{s}_i \circ \cdots s_q \otimes \varphi(s_i).$$

*Proof:* This follows from the long exact sequence of divided powers (cf. [Eis94, Appendix A2]). The reader may consult also [Sak78, Proposition 3]. □

Our strategy is the following. We apply Lemma 5.2.10 to the syzygy sequence

$$0 \to \mathrm{Syz}(x,y,z) \to \bigoplus_{i=1}^{3} \mathscr{O}_Y(-1) \xrightarrow{x,y,z} \mathscr{O}_Y \to 0,$$

and obtain a short exact sequence of vector bundles on $Y$

$$0 \to \mathrm{Sym}^q(\mathrm{Syz}(x,y,z)) \to \mathrm{Sym}^q\left(\bigoplus_{i=1}^{3} \mathscr{O}_Y(-1)\right) \xrightarrow{\phi} \mathrm{Sym}^{q-1}\left(\bigoplus_{i=1}^{3} \mathscr{O}_Y(-1)\right) \to 0 \qquad (5.6)$$

of ranks $q+1$, $\binom{q+2}{2}$, and $\binom{q+1}{2}$ respectively. From the formula $\mu(\mathrm{Sym}^q \mathscr{F}) = q\mu(\mathscr{F})$ of Proposition 5.1.10, we have that the slopes of the bundles of sequence (5.6) are $-\frac{9}{2}q$, $-3q$ and $-3q+3$ respectively.

**Remark 5.2.11.** We have that

$$\mathrm{Sym}^q\left(\bigoplus_{i=1}^{3} \mathscr{O}_Y(-1)\right) \cong \bigoplus_{\substack{i_1 \leq \cdots \leq i_q \\ i_j \in \{0,\ldots,3\}}} \mathscr{O}_Y(-q). \qquad (5.7)$$

Therefore the map $\phi$ of sequence (5.6) can be interpreted as a map between two splitting vector bundles

$$\phi \colon \bigoplus_{i_1 \leq \cdots \leq i_q} \mathscr{O}_Y(-q) \to \bigoplus_{i_1 \leq \cdots \leq i_{q-1}} \mathscr{O}_Y(-q+1)$$

and represented as a matrix with homogeneous entries in the coordinate ring $R$ of $Y$.

Now, we assume that $q$ is even, the case with $q$ odd has been clarified in Corollary 5.2.5. We tensor sequence (5.6) with the line bundle $\mathscr{O}_Y(\frac{3}{2}q)$ of degree $\frac{9}{2}q$. This yields

$$0 \to \mathrm{Sym}^q(\mathrm{Syz}(x,y,z)) \otimes \mathscr{O}_Y\left(\frac{3}{2}q\right) \to \bigoplus \mathscr{O}_Y\left(\frac{1}{2}q\right) \xrightarrow{\phi} \bigoplus \mathscr{O}_Y\left(\frac{1}{2}q+1\right) \to 0. \qquad (5.8)$$

This gives a syzygy representation for the bundle $E_q := \mathrm{Sym}^q(\mathrm{Syz}(x,y,z)) \otimes \mathscr{O}_Y\left(\frac{3}{2}q\right)$, which has degree 0, by Proposition 5.1.10 and Proposition 5.2.2. Since $E_q$ is obtained from $\mathrm{Sym}^q(\mathrm{Syz}(x,y,z))$ via tensorization with a twisted structure sheaf, the splitting sub-bundles of $E_q$ and $\mathrm{Sym}^q(\mathrm{Syz}(x,y,z))$ are in bijective correspondence. In other words, $\mathrm{frk}_{\mathscr{O}_Y(1)}\,E_q = \mathrm{frk}_{\mathscr{O}_Y(1)}\,\mathrm{Sym}^q(\mathrm{Syz}(x,y,z))$. Moreover the bundle $E_q$ is semistable of slope 0, so its only possible trivial subbundles are of degree 0, that is of the form $\mathscr{O}_Y$. This can make the computation of the free rank of $E_q$ easier.

In addition, we observe that from sequence (5.8) we have that every $\mathscr{O}_Y$-summand of $E_q$ gives rise to a section on the left, so we have an upper bound

$$\mathrm{frk}_{\mathscr{O}_Y(1)}\,\mathrm{Sym}^q(\mathrm{Syz}(x,y,z)) \leq h^0(E_q).$$

The global sections $h^0(E_q)$ can be counted by looking at elements of degree $\frac{1}{2}q$ in the kernel of the map $\phi$. If $q$ is not too big, this could be done by hand or with the help of a computer software.

**Example 5.2.12.** For $q = 2$, the sequence (5.8) looks like

$$0 \to \mathrm{Sym}^2(\mathrm{Syz}(x,y,z)) \otimes \mathscr{O}_Y(3) \to \bigoplus_{i=1}^{6} \mathscr{O}_Y(1) \xrightarrow{\phi} \bigoplus_{i=1}^{3} \mathscr{O}_Y(2) \to 0.$$

The map $\phi$ is given by the matrix

$$M_\phi = \begin{pmatrix} 2x & y & z & 0 & 0 & 0 \\ 0 & x & 0 & 2y & z & 0 \\ 0 & 0 & x & 0 & y & 2z \end{pmatrix}$$

with values in $R$. Its kernel is the bundle $\mathrm{Sym}^2(\mathrm{Syz}(x,y,z)) \otimes \mathscr{O}_Y(3)$. We have to look at elements of degree $\frac{1}{2}q = 1$ in the kernel of $\phi$, that is vectors $h = (h_1, \dots, h_6) \in \bigoplus_{i=1}^{6} \mathscr{O}_Y(1)$ such that $\phi(h) = 0$. In this case, it is clear that such an element $h$ does not exist. Observe that the defining equation of the curve $f$ does not play a role here, since $f$ has degree 3 and $M_\phi \cdot h$ has degree 2. Therefore, we have that $\mathrm{frk}_{\mathscr{O}_Y(1)}\,\mathrm{Sym}^2(\mathrm{Syz}(x,y,z)) = 0$, in accordance with Example 5.2.9.

We conclude the chapter with some considerations on the case $q = 4$.

**Example 5.2.13.** For $q = 4$, the sequence (5.8) looks like

$$0 \to \mathrm{Sym}^4(\mathrm{Syz}(x,y,z)) \otimes \mathscr{O}_Y(6) \to \bigoplus_{i=1}^{15} \mathscr{O}_Y(2) \xrightarrow{\phi} \bigoplus_{i=1}^{10} \mathscr{O}_Y(3) \to 0,$$

and the map $\phi$ is given by the matrix

$$M_\phi = \begin{pmatrix}
4x & y & z & 0 & \cdots & & 0 & \cdots & & \cdots & 0 & \cdots & & \cdots & 0 \\
0 & 3x & 0 & 2y & z & 0 & 0 & \ddots & & & & & & & \vdots \\
\vdots & \ddots & 3x & 0 & y & 2z & 0 & 0 & & & & & & & \\
& & & 2x & 0 & 0 & 3y & z & 0 & \ddots & & & & & \vdots \\
\vdots & & \ddots & 2x & 0 & 0 & 2y & 2z & 0 & 0 & & & & & 0 \\
0 & & & & 0 & 2x & 0 & 0 & y & 3z & 0 & \ddots & & & \vdots \\
\vdots & & & & \ddots & x & 0 & 0 & 0 & 4y & z & & & & \\
& & & & & & x & 0 & 0 & 0 & 3y & 2z & \ddots & & \vdots \\
\vdots & & & & & & \ddots & x & 0 & 0 & 0 & 2y & 3z & 0 \\
0 & \cdots & & \cdots & 0 & \cdots & & \cdots & 0 & x & 0 & 0 & 0 & y & 4z
\end{pmatrix}$$

with values in $R$. This time we have to look at homogeneous elements $h$ of degree 2 in the kernel of $\phi$, so the equation $f$ may play a role, since $M_\phi \cdot h$ is homogeneous of degree 3. For some specific equations $f$, it can be checked with Macaulay2 [GS] that we have no such element $h$, so $\mathrm{frk}_{\mathcal{O}_Y(1)} \mathrm{Sym}^4(\mathrm{Syz}(x, y, z)) = 0$, but to prove it in general something more is needed.

# A. Macaulay2 computations

### A.0.3. Computation of $\mathrm{Syz}_R^2(k)$ for the singularity $E_6$

```
 Macaulay2, version 1.6
with packages: ConwayPolynomials, Elimination, IntegralClosure, LLLBases,
               PrimaryDecomposition, ReesAlgebra, TangentCone
i1 : R=QQ[x,y]/(x^7+x^6+x^5+x^4+x^3+x^2+x+1,y^2-2,x^4+1);

i2 : S=R[u,v];

i3 : p1=u*v*(u^4-v^4);

i4 : p2=u^8+14*u^4*v^4+v^8;

i5 : p3=u^12-33*u^8*v^4-33*u^4*v^8+v^12;

i6 : sub(p1,{u=>x^2*u, v=>x^6*v})===p1

o6 = true

i7 : sub(p2,{u=>x^2*u, v=>x^6*v})===p2

o7 = true

i8 : sub(p3,{u=>x^2*u, v=>x^6*v})===p3

o8 = true

i9 : sub(p1,{u=>x^2*v, v=>x^2*u})===p1

o9 = true

i10 : sub(p2,{u=>x^2*v, v=>x^2*u})===p2

o10 = true

i11 : sub(p3,{u=>x^2*v, v=>x^2*u})===p3

o11 = true

i12 : sub(p1,{u=>(y/2)*(x*u+x^3*v), v=>(y/2)*(x*u+x^7*v)})===p1

o12 = true

i13 : sub(p2,{u=>(y/2)*(x*u+x^3*v), v=>(y/2)*(x*u+x^7*v)})===p2

o13 = true

i14 : sub(p3,{u=>(y/2)*(x*u+x^3*v), v=>(y/2)*(x*u+x^7*v)})===p3

o14 = true

i15 : I=ideal(p1,p2,p3);

o15 : Ideal of S

i16 : syz(gens(I))

o16 = {6, 0}  | -30u7-210u3v4 210u4v3+30v7 |
      {8, 0}  | 25u4v-5v5     -5u5+25uv4    |
      {12, 0} | 5v            5u            |

i17 : s1=vector{210*u^4*v^3+30*v^7,-5*u^5+25*u*v^4,5*u};
```

```
i18 : s2=vector{-210*u^3*v^4-30*u^7,-5*v^5+25*u^4*v,5*v};

i19 : sub(s1,{u=>x^2*u, v=>x^6*v})===x^2*s1

o19 = true

i20 : sub(s2,{u=>x^2*u, v=>x^6*v})===x^6*s2

o20 = true

i21 : sub(s1,{u=>x^2*v, v=>x^2*u})===x^2*s2

o21 = true

i22 : sub(s2,{u=>x^2*v, v=>x^2*u})===x^2*s1

o22 = true

i23 : sub(s1,{u=>(y/2)*(x*u+x^3*v), v=>(y/2)*(x*u+x^7*v)})===
                (y/2)*(x*s1+x^3*s2)

o23 = true

i24 : sub(s2,{u=>(y/2)*(x*u+x^3*v), v=>(y/2)*(x*u+x^7*v)})===
                (y/2)*(x*s1+x^7*s2)

o24 = true
```

### A.0.4. Computation of $\mathrm{Syz}_R^2(k)$ for the singularity $E_7$

```
Macaulay2, version 1.6
with packages: ConwayPolynomials, Elimination, IntegralClosure, LLLBases,
               PrimaryDecomposition, ReesAlgebra, TangentCone
i1 : R=QQ[x,y]/(x^7+x^6+x^5+x^4+x^3+x^2+x+1,y^2-2,x^4+1);

i2 : S=R[u,v];

i3 : n=u*v*(u^4-v^4);

i4 : p1=u^8+14*u^4*v^4+v^8;

i5 : q=u^12-33*u^8*v^4-33*u^4*v^8+v^12;

i6 : p2=n^2;

i7 : p3=n*q;

i8 : sub(p1,{u=>x^3*u, v=>x^5*v})===p1

o8 = true

i9 : sub(p2,{u=>x^3*u, v=>x^5*v})===p2

o9 = true

i10 : sub(p3,{u=>x^3*u, v=>x^5*v})===p3

o10 = true

i11 : sub(p1,{u=>x^2*v, v=>x^2*u})===p1

o11 = true

i12 : sub(p2,{u=>x^2*v, v=>x^2*u})===p2

o12 = true

i13 : sub(p3,{u=>x^2*v, v=>x^2*u})===p3

o13 = true

i14 : sub(p1,{u=>(y/2)*(x*u+x^3*v), v=>(y/2)*(x*u+x^7*v)})===p1
```

```
o14 = true

i15 : sub(p2,{u=>(y/2)*(x*u+x^3*v), v=>(y/2)*(x*u+x^7*v)})===p2

o15 = true

i16 : sub(p3,{u=>(y/2)*(x*u+x^3*v), v=>(y/2)*(x*u+x^7*v)})===p3

o16 = true

i17 : I=ideal(p1,p2,p3);

o17 : Ideal of S

i18 : syz(gens(I))

o18 = {8, 0}  | -35u9v2+42u5v6-7uv10 7u10v-42u6v5+35u2v9 |
      {12, 0} | 42u7+294u3v4         -294u4v3-42v7        |
      {18, 0} | -7v                  -7u                  |
                    3        2
o18 : Matrix S  <--- S

i19 : s1=vector{7*u^10*v-42*u^6*v^5+35*u^2*v^9,-294*u^4*v^3-42*v^7,-7*u}

o19 = | 7u10v-42u6v5+35u2v9 |
      | -294u4v3-42v7       |
      | -7u                 |
            3
o19 : S

i20 : s2=vector{-35*u^9*v^2+42*u^5*v^6-7*u*v^10,42*u^7+294*u^3*v^4,-7*v}

o20 = | -35u9v2+42u5v6-7uv10 |
      | 42u7+294u3v4         |
      | -7v                  |
            3
o20 : S

i21 : sub(s1,{u=>x^3*u, v=>x^5*v})===x^3*s1

o21 = true

i22 : sub(s2,{u=>x^3*u, v=>x^5*v})===x^5*s2

o22 = true

i23 : sub(s1,{u=>x^2*v, v=>x^2*u})===x^2*s2

o23 = true

i24 : sub(s2,{u=>x^2*v, v=>x^2*u})===x^2*s1

o24 = true

i25 : sub(s1,{u=>(y/2)*(x*u+x^3*v), v=>(y/2)*(x*u+x^7*v)})===(y/2)*(x*s1+x^3*s2)

o25 = true

i26 : sub(s2,{u=>(y/2)*(x*u+x^3*v), v=>(y/2)*(x*u+x^7*v)})===(y/2)*(x*s1+x^7*s2)

o26 = true
```

## A.0.5. Computation of $\operatorname{Syz}^2_R(k)$ for the singularity $E_8$

```
Macaulay2, version 1.6
with packages: ConwayPolynomials, Elimination, IntegralClosure, LLLBases,
               PrimaryDecomposition, ReesAlgebra, TangentCone
i1 : R=QQ[x,y]/(1+x+x^2+x^3+x^4,y^2-5);

i2 : S=R[u,v];
```

```
i3 : p1=u^11*v+11*u^6*v^6-u*v^11;

i4 : p2=u^20-228*u^15*v^5+494*u^10*v^10+228*u^5*v^15+v^20;

i5 : p3=u^30+522*u^25*v^5-10005*u^20*v^10-10005*u^10*v^20-522*u^5*v^25+v^30;

i6 : sub(p1, {u=>(y/5)*((-x+x^4)*u+(x^2-x^3)*v),
                v=>(y/5)*((x^2-x^3)*u+(x-x^4)*v)})===p1

o6 = true

i7 : sub(p2, {u=>(y/5)*((-x+x^4)*u+(x^2-x^3)*v),
                v=>(y/5)*((x^2-x^3)*u+(x-x^4)*v)})===p2

o7 = true

i8 : sub(p3, {u=>(y/5)*((-x+x^4)*u+(x^2-x^3)*v),
                v=>(y/5)*((x^2-x^3)*u+(x-x^4)*v)})===p3

o8 = true

i9 : sub(p1, {u=>(y/5)*((x^2-x^4)*u+(x^4-1)*v),
                v=>(y/5)*((1-x)*u+(x^3-x)*v)})===p1

o9 = true

i10 : sub(p2, {u=>(y/5)*((x^2-x^4)*u+(x^4-1)*v),
                v=>(y/5)*((1-x)*u+(x^3-x)*v)})===p2

o10 = true

i11 : sub(p3, {u=>(y/5)*((x^2-x^4)*u+(x^4-1)*v),
                v=>(y/5)*((1-x)*u+(x^3-x)*v)})===p3

o11 = true

i12 : I=ideal(p1,p2,p3);

o12 : Ideal of S

i13 : syz(gens(I))

o13 = {12, 0} | -12012u19+2054052u14v5-2966964u9v10-684684u4v15
      {20, 0} | 11011u10v+66066u5v6-1001v11
      {30, 0} | 1001v
      --------------------------------------------------------------------------
      -684684u15v4+2966964u10v9+2054052u5v14+12012v19  |
      -1001u11-66066u6v5+11011uv10                     |
      1001u                                            |

                3        2
o13 : Matrix S  <--- S

i14 : s1=vector{-684684*u^15*v^4+2966964*u^10*v^9+2054052*u^5*v^14+12012*v^19,
                -1001*u^11-66066*u^6*v^5+11011*u*v^10,1001*u}

o14 = | -684684u15v4+2966964u10v9+2054052u5v14+12012v19 |
      | -1001u11-66066u6v5+11011uv10                    |
      | 1001u                                           |

           3
o14 : S

i15 : s2=vector{-12012*u^19+2054052*u^14*v^5-2966964*u^9*v^10-684684*u^4*v^15,
                11011*u^10*v+66066*u^5*v^6-1001*v^11,1001*v}

o15 = | -12012u19+2054052u14v5-2966964u9v10-684684u4v15 |
      | 11011u10v+66066u5v6-1001v11                     |
      | 1001v                                           |

           3
o15 : S

i16 : sub(s1, {u=>(y/5)*((-x+x^4)*u+(x^2-x^3)*v),
```

```
                        v=>(y/5)*((x^2-x^3)*u+(x-x^4)*v)})===(y/5)*((-x+x^4)*s1+(x^2-x^3)*s2)
o16 = true
i17 : sub(s2, {u=>(y/5)*((-x+x^4)*u+(x^2-x^3)*v),
                        v=>(y/5)*((x^2-x^3)*u+(x-x^4)*v)})===(y/5)*((x^2-x^3)*s1+(x-x^4)*s2)
o17 = true
i18 : sub(s1, {u=>(y/5)*((x^2-x^4)*u+(x^4-1)*v),
                        v=>(y/5)*((1-x)*u+(x^3-x)*v)})===(y/5)*((x^2-x^4)*s1+(x^4-1)*s2)
o18 = true
i19 : sub(s2, {u=>(y/5)*((x^2-x^4)*u+(x^4-1)*v),
                        v=>(y/5)*((1-x)*u+(x^3-x)*v)})===(y/5)*((1-x)*s1+(x^3-x)*s2)
o19 = true
```